\documentclass[12pt]{article}
\usepackage{amsmath}
\usepackage{amssymb}
\usepackage{amsthm}
\usepackage{color}
\usepackage{mathrsfs}
\usepackage{graphicx}
\usepackage{float}
\usepackage{cases}
\usepackage{cite}
\usepackage{setspace}
\usepackage{bm}
\usepackage{enumitem}
\usepackage{titlesec}
\usepackage{dsfont}

\titleformat{\subsection}
  {\bf\normalsize} % Font style: \normalsize, \itshape (italic), etc.
  {\thesubsection}                % Numbering (e.g., 1.1)
  {1 em}                           % Space between number and title
  {}                              % Formatting before the title
\titlespacing*{\subsection}{0pt}{3.25ex plus 1ex minus .2ex}{1.5ex plus .2ex}
\usepackage{hyperref}
\newtheorem{theorem}{Theorem}[section]
 \newtheorem{corollary}[theorem]{Corollary}

 \newtheorem{lemma}[theorem]{Lemma}
 \newtheorem{remark}[theorem]{Remark}

 \title{Pointwise estimates for an anisotropic product type singular integral on the sphere}
 \author{Ahmed A. Abdelhakim}
 \date{}
\allowdisplaybreaks
\begin{document}

\maketitle
\begin{abstract}
Let $\Gamma=\{1,\dots,d\}$, $d\geq 2$. For each $j\in \Gamma$, let ${\ell}_{j}\geq 1$ be an integer and let $\alpha_{j}\in (1/2,1)$. We obtain a pointwise estimate for the integral
\begin{equation*}
\int_{\mathbb{S}^{d-1}}\prod_{j\in \Gamma}
|x^{\ell_{j}}_{j}-\xi_{j}|^{-\alpha_{j}} d
\sigma(x),
\end{equation*}
in the region
\begin{equation*}
\{\xi\in \mathbb{R}^{d}: \xi_{j}\neq 0,\;j\in\Gamma,\;\; \sum\nolimits_{j\in I}|\xi_{j}|^{2/\ell_{j}}\neq 1,\; I\subseteq \Gamma\}.
\end{equation*}
These estimates are used to prove a Tomas-Stein type restriction inequality for the oscillatory integral operator $Tf(x)=\int_{\mathbb{R}^d}e^{\dot{\imath}h(x,y)}f(y)dy$, with the phase $h(x,y)=\sum\nolimits_{j\in\Gamma} a_{j} x_{j}^{\ell_{j}} y_{j}^{m_{j}}$, where, for each $j\in \Gamma$, $a_{j}\neq 0$ is a real number and $m_{j}> 1$ is an integer.
\end{abstract}
\section{Introduction}
Let $d\geq 2$ and let $\Gamma$ denote the set $\{1,\dots,d\}$.
For each $j\in \Gamma$, let
$\alpha_{j}\in (1/2,1)$ and let ${\ell}_{j}\geq 1$ be an integer. Consider the integral
\begin{equation*}
F(\xi):=\int_{\mathbb{S}^{d-1}}\prod_{j\in \Gamma}
|x^{\ell_{j}}_{j}-\xi_{j}|^{-\alpha_{j}} d
\sigma(x),\quad \xi=(\xi_{1},\dots,\xi_{d})\in \mathbb{R}^{d},
\end{equation*}
where $d\sigma$ is the Lebesgue surface measure on the sphere $\mathbb{S}^{d-1}$ in
$\mathbb{R}^{d}$. The purpose of this study is to obtain pointwise bounds
for the integral $F(\xi)$ for all values of the exponents $\alpha_{j}\in (1/2,1)$, $j\in \Gamma$, in the region
\begin{equation*}
\mathcal{E}:=
\{\xi\in \mathbb{R}^{d}: \xi_{j}\neq 0,\,j\in\Gamma,\;\; \sum\nolimits_{j\in I}|\xi_{j}|^{2/\ell_{j}}\neq 1,\, I\subseteq \Gamma\}.
\end{equation*}
The integral $F(\xi)$ is convergent for every $\xi\in \mathbb{R}^{d}$ whenever each of the following  conditions is satisfied: $\alpha_{j}\ell_{j}<1$, $j\in \Gamma$, $\sum\nolimits_{j\in I}\alpha_{j}<|I|-1$, $I\subseteq\Gamma$, $|I|\geq 2$, and $2\alpha_{k}+\sum\nolimits_{j\in \Gamma\setminus \{k\}}
\alpha_{j}\ell_{j}<d-1$, $k\in \Gamma$. We are rather interested in the pointwise behaviour of the integral $F(\xi)$ for the larger values of the exponents $\alpha_{j}$. Outside the coordinate hyperplanes $\{\xi_{j}=0\}$, $j\in \Gamma$,
the restriction $\alpha_{j}\ell_{j}<1$ is no longer necessary, for any $j\in \Gamma$, for pointwise convergence of the integral $F(\xi)$. One must also avoid the hypersurfaces $\{\xi\in \mathbb{R}^{d}:\sum_{j\in I}|\xi_{j}|^{2/\ell_{j}}=1\}$, $I\subseteq \Gamma$, on which the hypersurfaces $\{x_{j}^{\ell_{j}}=\xi_{j}\}$, $j\in I$, may meet the sphere $\mathbb{S}^{d-1}$ tangentially (non-transversely) (see e.g. \cite{lee2013smooth,stein1970singular}). A minor modification of our approach yields a pointwise bound
in the case $0<\alpha_{j}\leq 1/2$, for some $j\in \Gamma$. We decided not to address that case explicitly here to avoid an unnecessarily cumbersome estimate.
\smallskip\par
If $\alpha_{j}<1/2$, for all $j\in \Gamma$, the integral $F(\xi)$ is convergent on $\{\xi\in \mathbb{R}^{d}:\, \xi_{j}\neq 0,\,j\in \Gamma\}$ in all dimensions $d\geq 3$. In that case, one could obtain a naive pointwise estimate for $F(\xi)$ using Young's inequality on the sphere (see \cite{carlen2004sharp,barthe2006entropy}).
The best one could get, applying Theorem 1.1 in \cite{carlen2004sharp} with
$p=2$, is that
\begin{equation*}
F(\xi)\lesssim \prod\nolimits_{j\in\Gamma,\,
2\alpha_{j}\ell_{j}=1} |\log{|\xi_{j}|}|^{1/2}
\prod\nolimits_{j\in\Gamma,\;2\alpha_{j}\ell_{j}>1}
|\xi_{j}|^{1/(2\ell_{j})
-\alpha_{j}},
\end{equation*}
for small $\xi$ (see e.g. \cite{grafakos2008classical}, Appendix D, for useful formulas for integration over the sphere). Nevertheless, the pointwise bound implied by our approach remains uniformly bounded as $|\xi|\rightarrow 0$, if $\alpha_{j}<\min\{1/2,1/\ell_{j}\}$, for every $j\in \Gamma$. This is consistent with the fact that $K(0)$ exits and is finite when $\alpha_{j}\ell_{j}< 1$, $j\in \Gamma$.
 \smallskip\par
The anisotropic structure of $F(\xi)$ suggests integrating in all coordinate directions. Roughly speaking, we
start by bounding $F(\xi)$ by a sum of integrals on spherical shell zones. These in turn are decomposed into sums of integrals on regions where the singularities are easier to integrate iteratively.
The decomposition naturally depends on $\xi$. At this point, the difficulty is that a careful pointwise estimate obtained for an integral in some coordinate direction will not depend only on $\xi$, but also on the region determined by the relation between the remaining variables and $\xi$.
\smallskip \par
For $j\in \Gamma$, let $a_{j}$ be a nonzero real number and  $m_{j}> 1$ be an integer.
We encountered the integral $F(\xi)$ in \cite{abdelhakim2026sharpl2mathbbrdl2mathbbsd1restrictiontheorem} while trying to characterize the exponents
$p\in [1,2]$ for which the oscillatory
integral operator $Tf(x)=\int_{\mathbb{R}^{d}}
e^{\dot{\imath} h(x,y)}f(y)dy$, $h(x,y)=\sum_{j\in\Gamma} a_{j} x_{j}^{\ell_{j}} y_{j}^{m_{j}}$, satisfies the restriction inequality
 $\|Tf\|_{L^{2}(\mathbb{S}^{d-1})}\lesssim
\|f\|_{L^{p}(\mathbb{R}^{d})}$. Together with the results in \cite{abdelhakim2023inhomogeneous}, the pointwise estimates proved for $F(\xi)$ here are key to solving this problem. For the classical theory of the restriction of the Fourier transform to the sphere
see e.g. \cite{stein1993harmonic,tomas1975restriction}.
For recent progress on the restriction theory, see e.g. \cite{foschi2017some,demeter1945fourier,stovall2019waves,tao2011some} and the many references therein.
\section{Pointwise estimate}\label{pstm1}
Let $\Gamma_{j}$ denote the set $\Gamma\setminus \{j\}$
and ${\overline{\ell}}_{j}$ denote
the reciprocal of the integer $\ell_{j}$, $j\in \Gamma$. Introduce the sets
\begin{equation*}
E_{1}=\{j\in \Gamma: \ell_{j}=1\},\quad
E_{2}=\{j\in \Gamma: \ell_{j}\geq 2\}.
\end{equation*}
For every $\xi\in \mathcal{E}$, we have
\begin{equation}\label{mnk}
F(\xi)\leq C\sum_{k\in \Gamma} \sum_{
\substack{
\Gamma_{k}^{+}\cup \Gamma_{k}^{-}= \Gamma_{k}\\
\Gamma_{k}^{+}\cap \Gamma_{k}^{-}= \emptyset}} \prod_{j\in \Gamma_{k}^{+}} |\xi_{j}|^{-\alpha_{j}}\chi_{[2,\infty)}(|\xi_{j}|) \prod_{n\in \Gamma_{k}^{-}}
\chi_{(0,2)}(|\xi_{n}|)H_{k,\Gamma_{k}^{-}}(\xi),
\end{equation}
with some constant $C$ independent of $\xi$, where
\begin{align*}
H_{k,\Gamma_{k}^{-}}(\xi)&:=
P_{\Gamma^{-}_{k}\cap E_{2}}(\xi)f_{k}(|\xi_{k}|)+
(\mathcal{R}_{k,\Gamma_{k}^{-}}(\xi)+
\mathcal{S}_{k,\Gamma_{k}^{-}}(\xi))g_{k}(|\xi_{k}|),\\
f_{k}(t)&:=
\chi_{(0,2]}(t)+t^{-\alpha_{k}}
\chi_{[1,\infty)}(t),\\
g_{k}(t)&:=\chi_{(d^{-\ell_{k}}/2,2)}(t),\\
\mathcal{S}_{k,\Gamma_{k}^{-}}(\xi)&:=
\sum_{\substack{I\subseteq \Gamma^{-}_{k}\\
\beta_{|\Gamma^{+}_{k}|}(I)<0}}
\sum_{J\subseteq \Gamma^{-}_{k}\setminus I}
M^{\beta_{|\Gamma^{+}_{k}|}(I)+|J|}_{\Gamma^{-}_{k}\setminus J}(\xi)
\,P_{\Gamma^{-}_{k}\setminus ((I\cap E_{1})\cup J)}(\xi)
\,Q^{+}_{J}(\xi),\\
\mathcal{R}_{k,\Gamma_{k}^{-}}(\xi)&:=
\sum_{\substack{I\subseteq\Gamma_{k}^{-}\\
-2<\beta_{|\Gamma_{k}^{+}|}(I)\leq 0}}
\sum_{J\subseteq \Gamma_{k}^{-}\setminus I}
\sum_{\substack{K\subseteq J\\
-1<\theta_{K}(\beta_{|\Gamma_{k}^{+}|}(I))\leq 0}}
\mathcal{L}_{\Gamma_{k}^{-}\setminus (I\cup J)}(\xi)\,
\mathcal{N}_{I,J,K}(\xi;\theta_{K}(\beta_{ |\Gamma_{k}^{+}|}(I)))\\
&\qquad+\sum_{j\in \Gamma_{k}^{-}}\,
\sum_{\substack{I\subseteq \Gamma_{k}^{-}\setminus\{j\}\\
-2<\beta_{|\Gamma_{k}^{+}|}(I)\leq -1}}
\sum_{J\subseteq \Gamma_{k}^{-}\setminus (I\cup\{j\})}
\sum_{\substack{K\subseteq J\\
-1<\theta_{K}(\beta_{|\Gamma_{k}^{+}|}(I))\leq -1/2}}
Q_{J, K\cup \{j\}}(\xi_{j})\,\\&
\hspace{6 cm}\mathcal{N}_{I,J,K\cup \{j\},K}(\xi;\theta_{K}(\beta_{|\Gamma_{k}^{+}|}(I))+1/2),
\end{align*}
for each $I\subseteq \Gamma_{k}^{-}$,
\begin{align*}
P_{I}(\xi)&:=\,\prod_{j\in I}
|\xi_{j}|^{\overline{\ell}_{j}-\alpha_{j}},\\
\mathcal{L}_{I}(\xi)&:=1+\sum_{j\in I}|\log{|\xi_{j}|}|,\\
M_{I}(\xi)&:=\,\sum_{j\in I}
|\xi_{j}|^{\overline{\ell}_{j}}+\zeta_{k},\qquad \zeta_{k}^{2}:=|1-|\xi_{k}||,\\
Q^{+}_{I}(\xi)&:=\,\prod_{j\in I}
|\xi_{j}|^{-\alpha_{j}}\chi_{\{|\xi_{j}|^{\overline{\ell}_{j}}
\gg M_{\Gamma^{-}_{k}\setminus I}(\xi)\}}(\xi_{j}),\\
\lambda_{I}(\xi)&:=\sum_{j\in \Gamma^{-}_{k}\setminus I}|\xi_{j}|^{2{
\overline{\ell}}_{j}},\\
\Phi_{I}(\xi)&:=\,|1-\sum_{n\in I}{
|\xi_{n}|^{2\overline{\ell}_{n}}}-|\xi_{k}|^{
2{\overline{\ell}}_{k}}|,\\
X_{I}(\xi)&:=\chi_{\{\sum_{j\in I}|\xi_{j}|^{\overline{\ell}_{j}}\lesssim \zeta_{k}\}}(\xi),\\
\beta_{m}(I)&:=\,|I\cap E_{1}|-\sum_{j\in I\cap E_{1}}\alpha_{j}+m-2\alpha_{k},\quad
m\in \mathbb{Z}_{\geq 0},\\
\theta_{I}(\sigma)&:=|I|-
\sum_{j\in I}\alpha_{j}+\sigma/2,\quad -2<\sigma\leq 0,
\end{align*}
and for each $I,J,K,L\subseteq \Gamma_{k}^{-}$,
\begin{align*}
\mathcal{N}_{I,J,K,L}(\xi;\gamma)&:=
P_{\Gamma_{k}^{-}\setminus((I\cap E_{1}) \cup K)}(\xi)\,P^{+}_{J, L}(\xi)\,X_{J}(\xi)\,
\mathcal{W}_{J,K}(\xi;\gamma),\quad -1<\gamma\leq 0,\\[0.1 cm]
P^{+}_{J, L}(\xi)&:=\,
\prod\nolimits_{n\in L}
|\xi_{n}|^{2\alpha_{n}\overline{\ell}_{n}-\alpha_{n}-
\overline{\ell}_{n}}
\chi_{\{|\xi_{n}|^{2\overline{\ell}_{n}}\gtrsim
\Phi_{J}(\xi)+\lambda_{L}(\xi)\}}(\xi_{n}),\\[0.1 cm]
\mathcal{W}_{J,K}(\xi;\gamma)&:=
\mathds{1}_{\gamma=0}|\log{(\Phi_{J}(\xi)+\lambda_{K}(\xi))}|+
\mathds{1}_{-1<\gamma<0}(\Phi_{J}(\xi)+\lambda_{K}(\xi))^{\gamma},\\[0.1 cm]
\mathcal{N}_{I,J,K}(\xi;\gamma)&:=
\mathcal{N}_{I,J,K,K}(\xi;\gamma),\\[0.1 cm]
Q_{J, K}(\xi_{j})&:=\,
|\xi_{j}|^{-\alpha_{j}}
\chi_{\{|\xi_{j}|^{2\overline{\ell}_{j}}\gtrsim
\Phi_{J}(\xi)+\lambda_{K}(\xi)\}}(\xi_{j}).
\end{align*}
\section{Proof}
Assume that $1/2<\alpha_{j}<1$, $j\in \Gamma$.
We start with the elementary observation that
\begin{equation*}
F(\xi)\leq 2\sup_{\epsilon>0}{
\epsilon^{-1} \int_{\sqrt{1-\epsilon}\,\leq |x|\leq 1} \prod\nolimits_{j\in \Gamma}|x_{j}^{{\ell}_{j}}-\xi_{j}|^{
-\alpha_{j}}{d x}},
\end{equation*}
which is finite for every $\xi\in \mathcal{E}$.
By the triangle inequality and symmetry, it suffices to estimate
\begin{equation*}
L_{\epsilon}(\xi):=
\int_{\{\sqrt{1-\epsilon}\,\leq |x|\leq 1\}\cap \mathbb{R}^{d}_{+}} \prod_{j\in \Gamma}|x_{j}^{{\ell}_{j}}-\xi_{j}|^{
-\alpha_{j}}{d x}
\end{equation*}
in the region
\begin{equation*}
\{\xi\in (0,\infty)^{d}:\;\sum\nolimits_{j\in I} \xi_{j}^{2/\ell_{j}}\neq 1,\; I\subseteq \Gamma\}.
\end{equation*}
Let $\epsilon<(d-1)/d^2$ and introduce the sets
\begin{equation*}
\mathcal{Y}_{\epsilon}^{j}:=\left\{
x\in \mathbb{R}^{d}_{+}:1-\epsilon\leq |x|^{2}\leq 1,\;
|x^{\prime}_{j}|^{2}\leq 1-\epsilon-1/d^2\right\},\quad
j\in \Gamma,
\end{equation*}
where $x^{\prime}_{j}$ denotes the point
$(x_{1},\dots,x_{j-1},x_{j+1},\dots,x_{d})$
in $\mathbb{R}^{d-1}$. Since
\begin{equation*}
\{\sqrt{1-\epsilon}\,\leq |x|\leq 1\}\cap \mathbb{R}^{d}_{+}=\bigcup_{k\in \Gamma} \mathcal{Y}_{\epsilon}^{k},
\end{equation*}
one has
\begin{equation}\label{lps}
L_{\epsilon}(\xi)\leq \sum_{k\in \Gamma}
L^{k}_{\epsilon}(\xi),
\end{equation}
where
\begin{equation*}
L^{k}_{\epsilon}(\xi):=
\int_{\mathcal{Y}_{\epsilon}^{k}}
\prod_{j\in \Gamma}|x_{j}^{{\ell}_{j}}-\xi_{j}|^{
-\alpha_{j}} d x.
\end{equation*}
Changing variables $x^{\ell_{j}}_{j}\rightarrow
x_{j}$, $j\in \Gamma$, and applying Fubini's theorem, we see that
\begin{equation}\label{gy1}
L^{k}_{\epsilon}(\xi)\,\lesssim\,
\int_{x^{\prime}_{k}\in \mathbb{R}_{+}^{d-1},\,
\tilde{d}\leq \psi_{\epsilon}(x_{k}^{\prime})\leq 1}
\prod_{j\in \Gamma_{k}} x_{j}^{\overline{\ell}_{j}-1}
|x_{j}-\xi_{j}|^{
-\alpha_{j}}\,\vartheta_{\epsilon}(x_{k}^{\prime})\,
{d x^{\prime}_{k}},
\end{equation}
where
\begin{align*}
\vartheta_{\epsilon}(x^{\prime}_{k})&:=\,\int_{\tilde{d}\leq
\psi_{\epsilon}(x_{k}^{\prime})\leq
x_{k}\leq \psi(x_{k}^{\prime})}
{ x_{k}^{\overline{\ell}_{k}-1}
|x_{k}-\xi_{k}|^{-
\alpha_{k}}}{d x_{k}},\\
\psi_{\epsilon}(x_{k}^{\prime})&:=\,
(1-\epsilon-\sum_{j\in \Gamma_{k}}
x_{j}^{2\overline{\ell}_{j}})^{\ell_{k}/2},\\
\psi(x_{k}^{\prime})&:=\,
(1-\sum_{j\in \Gamma_{k}}
x_{j}^{2\overline{\ell}_{j}})^{\ell_{k}/2},
\end{align*}
and $\tilde{d}$ denotes $d^{-\ell_{k}}$. One has
\begin{align*}
\vartheta_{\epsilon}(x^{\prime}_{k})
&\approx
\int_{\tilde{d}\leq
\psi_{\epsilon}(x_{k}^{\prime})\leq
x_{k}\leq \psi(x_{k}^{\prime})}
|x_{k}-\xi_{k}|^{-\alpha_{k}}{d x_{k}} \\[0.05 cm]
&\approx
\begin{cases}
|\psi(x_{k}^{\prime})-\xi_{k}|^{1-
\alpha_{k}}+|\psi_{\epsilon}(x_{k}^{\prime})-\xi_{k}|^{1-
\alpha_{k}},& \psi_{\epsilon}(x_{k}^{\prime})\leq \xi_{k}\leq
\psi(x_{k}^{\prime});\\
||\psi(x_{k}^{\prime})-\xi_{k}|^{1-
\alpha_{k}}-|\psi_{\epsilon}(x_{k}^{\prime})-\xi_{k}|^{1-
\alpha_{k}}|,& \text{otherwise}.
\end{cases}
\end{align*}
Since $|\psi(x_{k}^{\prime})-\psi_{\epsilon}(x_{k}^{\prime})|
\lesssim\epsilon$, by the mean value theorem, it follows that
%%when $\tilde{d}\leq
%\psi_{\epsilon}(x_{k}^{\prime})\leq \psi(x_{k}^{\prime})\leq 1$. Thus, if $|\psi(x_{k}^{\prime})-\xi_{k}|\lesssim \epsilon$, we also have $|\psi_{\epsilon}(x_{k}^{\prime})-\xi_{k}|\lesssim \epsilon$,
%by the triangle inequality. It follows then from (\ref{j0}) that
%$\vartheta_{\epsilon}(x^{\prime}_{k})\lesssim \epsilon^{1-\alpha_{k}}$,
%wherever $|\psi(x_{k}^{\prime})-\xi_{k}|\lesssim \epsilon$.
%If $|\psi(x_{k}^{\prime})-\xi_{k}|\gg  \epsilon$, we may apply the mean value theorem to $t\mapsto |\psi(x_{k}^{\prime})-\xi_{k}-t|^{1-
%\alpha_{k}}$ on $[0,\psi(x_{k}^{\prime})-\psi_{\epsilon}(x_{k}^{\prime})]$,
%then use (\ref{mmu}). In either case, we have
\begin{equation*}
\vartheta_{\epsilon}(x^{\prime}_{k})\,\lesssim\,\epsilon\,
|\psi(x_{k}^{\prime})-\xi_{k}|^{-
\alpha_{k}},\quad \xi_{k}\neq\psi(x_{k}^{\prime}).
\end{equation*}
Using this in (\ref{gy1}), we deduce that
\begin{equation}\label{kl2}
L^{k}_{\epsilon}(\xi)\lesssim\,\epsilon L^{k}(\xi),
\end{equation}
where
\begin{equation*}
L^{k}(\xi):=\int_{\substack{x^{\prime}_{k}\in
\mathbb{R}_{+}^{d-1}\\
\tilde{d}\leq \psi(x_{k}^{\prime})\leq 1}}
\prod_{j\in \Gamma_{k}}
x_{j}^{\overline{\ell}_{j}-1}
|x_{j}-\xi_{j}|^{
-\alpha_{j}}
|\psi(x_{k}^{\prime})-\xi_{k}|^{-
\alpha_{k}}
{d x^{\prime}_{k}}.
\end{equation*}
We need to consider only one of the integrals $L^{k}$. Let us estimate $L^{d}$ on $\mathbb{R}_{+}^{d}$. Update
$\tilde{d}=d^{-\ell_{d}}$ and split
\begin{equation}\label{hyg0}
\begin{aligned}[b]
L^{d}(\xi)&=\int_{\substack{x^{\prime}_{d}\in
\mathbb{R}_{+}^{d-1}\\
\tilde{d}\leq \psi(x_{d}^{\prime})\leq 1\\
\psi(x_{d}^{\prime})\geq 2 \xi_{d}}}+
\int_{\substack{x^{\prime}_{d}\in
\mathbb{R}_{+}^{d-1}\\
\tilde{d}\leq \psi(x_{d}^{\prime})\leq 1\\ \xi_{d}\geq 2 \psi(x_{d}^{\prime})}}+\int_{\substack{x^{\prime}_{d}\in
\mathbb{R}_{+}^{d-1}\\
\tilde{d}\leq \psi(x_{d}^{\prime})\leq 1\\ \psi(x_{d}^{\prime})/2
< \xi_{d}< 2 \psi(x_{d}^{\prime})}}\\
&\hspace{2 cm}\prod_{j\in \Gamma_{d}}
x_{j}^{\overline{\ell}_{j}-1}
|x_{j}-\xi_{j}|^{
-\alpha_{j}}
|\psi(x_{d}^{\prime})-\xi_{d}|^{-
\alpha_{d}}
{d x^{\prime}_{d}}\\
&=:G_{1}(\xi)+G_{2}(\xi)+W(\xi).
\end{aligned}
\end{equation}
\par
Take a partition\footnote{\thanks{Throughout, we allow members of a partition of any subset of  $\Gamma$ to be empty.}} $\{\Gamma^{+}_{d},\Gamma^{-}_{d}\}$ of $\Gamma_{d}$ and let
\begin{equation*}
U_{j}=
\begin{cases}
(0,2),& j\in \Gamma^{-}_{d};\\[0.05 cm]
[2,\infty),& j\in \Gamma^{+}_{d}.
\end{cases}
\end{equation*}
%By this set up, each partition
%$\{\Gamma^{+}_{d},\Gamma^{-}_{d}\}$ of $\Gamma_{d}$ determines a distinct
%region
%\begin{equation}\label{cp1}
%U:=U_{1}\times\dots\times U_{d-1}
%\end{equation}
%in $\mathbb{R}_{+}^{d-1}$.
We shall estimate each of the integrals in (\ref{hyg0}) on $U \times \mathbb{R}_{+}$ where
\begin{equation*}
U:=U_{1}\times\dots\times U_{d-1}.
\end{equation*}
%The pointwise bound we obtain for each of these integrals depends explicitly on
%the arbitrarily chosen partition $\{\Gamma^{+}_{d},\Gamma^{-}_{d}\}$. We sum this pointwise bound over all possible partitions $\{\Gamma^{+}_{d},\Gamma^{-}_{d}\}$ of $\Gamma_{d}$. Since $\mathbb{R}_{+}^{d-1}$ is
%the union of the regions $U$ over these partitions,
%we eventually obtain an estimate for $G$ on $\mathbb{R}^{d}_{+}$.
The integral $W$ is obviously the most difficult to estimate.
We further decompose the sets $U_{j}$ with $j\in \Gamma^{-}_{d}$
in different ways, each corresponding to a certain partition of $\Gamma^{-}_{d}$, and then obtain local estimates
for $W$. This localization is discussed in more detail in Section \ref{adobh}. Observe that $\alpha_{j}\ell_{j}<1$ if and only if $j\in E_{1}$.
\subsection{Estimates for \texorpdfstring{{\mathversion{bold}$G_{1}$}}{}
and \texorpdfstring{{\mathversion{bold}$G_{2}$}}{}}\label{sc41}
Update Suppose that
$\tilde{d}\leq \psi(x_{d}^{\prime})\leq 1$
and $\psi(x_{d}^{\prime})\geq 2 \xi_{d}$. Then
\begin{equation*}
\psi(x_{d}^{\prime})-\xi_{d}\geq
\psi(x_{d}^{\prime})/2\geq \tilde{d}/2,\quad
0<\xi_{d}\leq 1/2.
\end{equation*}
If $\xi_{d}> 1/2$, the function $G_{1}$ does not appear in the sum (\ref{hyg0}). So we may as well assume that
$G_{1}$ is supported in the strip
$0\leq \xi_{d}\leq 1/2$. Moreover, if $\Gamma^{+}_{d}\neq\emptyset$ and $j\in \Gamma^{+}_{d}$, one clearly has $|x-\xi_{j}|\approx \xi_j$, for $0<x_{j}<1$. The restriction $1\geq \psi(x_{d}^{\prime})\geq 2\xi_{d}$ implies that
$x_{j}\leq \widetilde{\zeta}_{d}^{\,\ell_{j}}$,
for all $j\in \Gamma_{d}$, where $\widetilde{\zeta}_{d}=(1-
(2\xi_{d})^{2\overline{\ell}_{d}})^{1/2}$.
We therefore have
\begin{align}
\nonumber G_{1}(\xi)&\lesssim\,
\prod_{j\in \Gamma_{d}}
\int_{0<x_{j}\leq \widetilde{\zeta}_{d}^{\,\ell_{j}}}
x_{j}^{\overline{\ell}_{j}-1}
|x_{j}-\xi_{j}|^{
-\alpha_{j}}
{d x_{j}}\\
\nonumber &\lesssim\, \widetilde{\zeta}_{d}^{\,|\Gamma^{+}_{d}|+
\sum_{j\in \Gamma^{-}_{d}\cap E_{1}}(1-\alpha_{j})}\prod_{j\in \Gamma^{+}_{d}}
\xi_{j}^{-\alpha_{j}}
\prod_{k\in \Gamma^{-}_{d}\cap E_{2}}
\xi_{k}^{-\alpha_{k}}
\min\{\xi_{k}^{\overline{\ell}_{k}},
\widetilde{\zeta}_{d}\},\\
\label{hyg1} &\lesssim\,
\prod_{j\in \Gamma^{+}_{d}}
\xi_{j}^{-\alpha_{j}}
\prod_{k\in \Gamma^{-}_{d}\cap E_{2}}
\xi_{k}^{\overline{\ell}_{k}-\alpha_{k}},\quad
\xi \in U\times (0,1/2].
\end{align}
\smallskip\par
Let us turn to $G_{2}$. If $\xi_{d}\geq 2 \psi(x_{d}^{\prime})$ and $\tilde{d}\leq \psi(x_{d}^{\prime})\leq 1$, we have
\begin{equation*}
|\psi(x_{d}^{\prime})-\xi_{d}|\geq
\xi_{d}/2,\quad \xi_{d}\geq 2 \tilde{d}.
\end{equation*}
The function $G_{2}$ does not contribute to the sum (\ref{hyg0}) wherever $\xi_{d}<2 \tilde{d}$. We can thus assume
that $G_{2}$ is supported in the region
$\xi_{d}\geq 2 \tilde{d}$. By a slight modification of the proof of (\ref{hyg1}), it follows that
\begin{equation}\label{hyg2}
 G_{2}(\xi)\lesssim\,\xi_{d}^{-\alpha_{d}}
\prod_{j\in \Gamma^{+}_{d}}
\xi_{j}^{-\alpha_{j}}\,
\prod_{k\in \Gamma^{-}_{d}\cap E_{2}}\xi_{k}^{\overline{\ell}_{k}-\alpha_{k}},\quad
\xi \in U\times[2 \tilde{d},\infty).
\end{equation}
\subsection{Preliminary estimates for
\texorpdfstring{{\mathversion{bold}$W$}}{}}\label{prhstr}
In the regions $\xi_{d}\leq \tilde{d}/2$ and $\xi_{d}\geq 2$, the function $W$
does not appear in the sum (\ref{hyg0}). We may therefore assume that
$W$ is supported in the strip
$\tilde{d}/2<\xi_{d}< 2$. It is convenient here to change
variables $x_{j}\rightarrow x_{j}^{\ell_{j}}$
to see that
\begin{equation}\label{cnsi}
W(\xi)\lesssim\int_{R}
\Pi(x_d^{\prime},\xi)
{d x^{\prime}_{d}},
\end{equation}
where
\begin{equation*}
\Pi(x_d^{\prime},\xi):=|(1-|x^{\prime}_{d}|^{2})^{
\ell_{d}/2}-\xi_{d}|^{-\alpha_{d}}\prod_{j\in \Gamma_{d}}
|x^{\ell_{j}}_{j}-\xi_{j}|^{
-\alpha_{j}},
\end{equation*}
and
\begin{equation*}
R:=\left\{x^{\prime}_{d}\in \mathbb{R}_{+}^{d-1}:
0\leq
|x^{\prime}_{d}|\leq \sqrt{1-1/d^2},\,
\xi_{d}^{2{\overline{\ell}}_{d}}/4<1-|x^{\prime}_{d}|^2<4\xi_{d}^{2{\overline{\ell}}_{d}}\right\}.
\end{equation*}
We must estimate $W(\xi)$ on $\mathbb{R}_{+}^{d-1}\times V$, where $V$ denotes the interval
$(\tilde{d}/2,2)$.
Note here that $(1/2,2)\subseteq V\subsetneq (0,2)$ for any $d\geq 2$, and that
\begin{equation}\label{1mxdz}
\zeta_{d}^{2}:=|1-\xi_{d}|\approx|1-\xi^{2{\overline{\ell}}_{d}}_{d}|,
\quad\text{whenever}\;\;\tilde{d}/2<\xi_{d}<2.
\end{equation}
For all $x^{\prime}_{d}\in R$
and $\xi_{d}\in V$, the factor
\begin{equation*}
|(1-|x^{\prime}_{d}|^{2})^{
\ell_{d}/2}-\xi_{d}| \approx
||x^{\prime}_{d}|^{2}-(1-\xi^{
2{\overline{\ell}}_{d}}_{d})|.
\end{equation*}
By this estimate, it follows from  (\ref{cnsi}) that
\begin{equation}\label{cnsin}
W(\xi)\lesssim \prod_{j\in \Gamma^{+}_{d}}\xi_{j}^{
-\alpha_{j}}\int_{R}
\widetilde{W}(x_d^{\prime},\xi)
{d x^{\prime}_{d}},
\end{equation}
where
\begin{equation*}
\widetilde{W}(x_d^{\prime},\xi):=
||x^{\prime}_{d}|^{2}-(1-\xi^{
2{\overline{\ell}}_{d}}_{d})|^{-\alpha_{d}}
\prod_{k\in \Gamma^{-}_{d}}
|x^{\ell_{k}}_{k}-\xi_{k}|^{
-\alpha_{k}}.
\end{equation*}
\subsubsection{An estimate for
\texorpdfstring{{\mathversion{bold}$W$}}{}
on
\texorpdfstring{
{\mathversion{bold}$U\times V$}}{}
 when \texorpdfstring{{\mathversion{bold}$\Gamma^{-}_{d}=\emptyset$}}{}}
If $\Gamma^{-}_{d}=\emptyset$, we are in the region where
$\xi_{j}> 2$, for all $j\in \Gamma_{d}$, whence, by (\ref{1mxdz}),
\begin{align}
\nonumber W(\xi)&\lesssim\,
\prod_{j\in \Gamma_{d}} \xi^{-\alpha_{j}}_{j}
\int_{R}\frac{d x^{\prime}_{d}}{
||x^{\prime}_{d}|^{2}-(1-\xi^{
2{\overline{\ell}}_{d}}_{d})|^{\alpha_{d}}}\\[0.05 cm]
\nonumber &\lesssim\,
\max{\{1,\zeta_{d}^{d-1-2\alpha_{d}}\}}\prod_{j\in \Gamma_{d}} \xi^{-\alpha_{j}}_{j}\\
\label{sh0} &\lesssim\,\left\{
                 \begin{array}{ll}
                  \xi^{-\alpha_{1}}_{1}
\zeta_{2}^{1-2\alpha_{2}},& \hbox{if $d=2$;}\vspace{0.2 cm} \\
                   \prod_{j\in \Gamma_{d}} \xi^{-\alpha_{j}}_{j},& \hbox{if $d\geq 3$.}
                 \end{array}
               \right.
\end{align}
\subsubsection{A decomposition of \texorpdfstring{
{\mathversion{bold}$W$}}{} when
\texorpdfstring{{\mathversion{bold}$\Gamma^{-}_{d}\neq \emptyset$}}{}}\label{442hstr}
Assume now that $\Gamma^{-}_{d}\neq \emptyset$.
Decompose $R=R_{1}\cup R_{2}\cup
R_{3}$ where
\begin{align*}
R_{1} &:=\,
\left\{x^{\prime}_{d}\in R:
2|x^{\prime}_{d}|^{2}< |1-\xi^{
2{\overline{\ell}}_{d}}_{d}|\right\}, \\
R_{2} &:=\,
\left\{x^{\prime}_{d}\in R:
|x^{\prime}_{d}|^{2}> 2|1-\xi^{
2{\overline{\ell}}_{d}}_{d}|\right\}, \\
R_{3} &:=\,
\left\{x^{\prime}_{d}\in R:
|1-\xi^{
2{\overline{\ell}}_{d}}_{d}|/2<|x^{\prime}_{d}|^{2}<2|1-\xi^{
2{\overline{\ell}}_{d}}_{d}|\right\}.
\end{align*}
By (\ref{1mxdz}), one has
\begin{equation}\label{xdr12}
||x^{\prime}_{d}|^{2}-(1-\xi^{
2{\overline{\ell}}_{d}}_{d})|\approx
\left\{
  \begin{array}{ll}
  \zeta_{d}^{2},& \hbox{$x^{\prime}_{d}\in R_{1}$;} \\[0.2 cm]
  |x^{\prime}_{d}|^{2},& \hbox{$x^{\prime}_{d}\in R_{2}$.}
  \end{array}
\right.
\end{equation}
Furthermore, for every $k\in \Gamma^{-}_{d}$, one has
\begin{equation}\label{xdr123p}
|x^{\ell_{k}}_{k}-\xi_{k}|\approx
\left\{
  \begin{array}{ll}
    x^{\ell_{k}}_{k},& \hbox{$x_{k}>2\xi_{k}^{\overline{\ell}_{k}}$;} \\[0.05 cm]
    \xi_{k},& \hbox{$x_{k}<\xi_{k}^{\overline{\ell}_{k}}/2$;}  \\[0.05 cm]
    \xi_{k}^{1-{\overline{\ell}}_{k}}
|x_{k}-\xi^{\overline{\ell}_{k}}_{k}|,& \hbox{$\xi_{k}^{\overline{\ell}_{k}}/2<x_{k}<2\xi_{k}^{
{\overline{\ell}}_{k}}$.}
  \end{array}
\right.
\end{equation}
The estimates (\ref{xdr12}) and  (\ref{xdr123p}) above suggest splitting
\begin{equation}\label{smsm}
\int_{ R}\widetilde{W}(x_d^{\prime},\xi)
{d x^{\prime}_{d}}=\sum_{m=1,2,3}\,\sum_{\substack{
I_{1}\cup I_{2}\cup I_{3}=\Gamma^{-}_{d},\\
I_{j}\cap I_{k}=\emptyset,\; j\neq k}}
\int_{ R_{m}\cap{\Theta}(I_{1},I_{2},I_{3})}\widetilde{W}(x_d^{\prime},\xi)
{d x^{\prime}_{d}},
\end{equation}
where $ {\Theta}(I_{1},I_{2},I_{3})\subseteq \mathbb{R}^{d-1}_{+}$ is the region
\begin{equation*}
\left\{x_{k}>2\xi_{k}^{\overline{\ell}_{k}},\,k\in I_{1} \right\}\cap\left\{x_{k}<\xi_{k}^{\overline{\ell}_{k}}/2,\,k\in I_{2} \right\} \cap\left\{\xi_{k}^{\overline{\ell}_{k}}/2<
x_{k}<2\xi_{k}^{\overline{\ell}_{k}},\,k\in I_{3} \right\}.
\end{equation*}
\par Let $\{I_{1},I_{2},I_{3}\}$ be a partition
of $\Gamma^{-}_{d}$. Since we are working under the assumption that
 $\Gamma^{-}_{d}\neq \emptyset$, we allow at most two of the sets $I_{1}$, $I_{2}$ and $I_{3}$ to be empty. From (\ref{xdr12}) and (\ref{xdr123p}) together,
we have
\begin{equation*}
\widetilde{W}(x_d^{\prime},\xi)\approx
\prod_{k\in I_{2}}\xi^{-\alpha_{k}}_{k}
\prod_{n\in I_{3}}
\,\xi_{n}^{-\alpha_{n}(1-{\overline{\ell}}_{n})}
\mathcal{Q}(x_d^{\prime},\xi), \;\;\text{whenever}\; x^{\prime}_{d}\in \Theta(I_{1},I_{2},I_{3}),
\end{equation*}
where
\begin{equation}\label{xdr123}
\mathcal{Q}(x_d^{\prime},\xi):=\prod_{j\in I_{1}} x^{-\alpha_{j}\ell_{j}}_{j}
\prod_{n\in I_{3}}
|x_{n}-\xi^{\overline{\ell}_{n}}_{n}|^{-\alpha_{n}}
\widetilde{\mathcal{Q}}(x_d^{\prime},\xi),
\end{equation}
and
\begin{equation*}
\widetilde{\mathcal{Q}}(x_d^{\prime},\xi):=\left\{
  \begin{array}{ll}
    \zeta_{d}^{-2\alpha_{d}},& \hbox{$x^{\prime}_{d}\in R_{1}\cap \Theta(I_{1},I_{2},I_{3})$;} \\[0.2 cm]
  (\sum_{j\in \Gamma_{d}}x_{j})^{-2\alpha_{d}},& \hbox{$x^{\prime}_{d}\in R_{2}\cap \Theta(I_{1},I_{2},I_{3})$;}    \\[0.2 cm]
 |\sum_{j\in \Gamma_{d}}x^{2}_{j}-(1-\xi^{
2{\overline{\ell}}_{d}}_{d})|^{-\alpha_{d}},& \hbox{$x^{\prime}_{d}\in R_{3}\cap \Theta(I_{1},I_{2},I_{3})$.}
  \end{array}
\right.
\end{equation*}
In view of (\ref{cnsin}), (\ref{smsm}) and the latter
estimate for $\widetilde{W}(\cdot,\xi)$ on
$\Theta(I_{1},I_{2},I_{3})$, we have
\begin{equation}\label{smsm1}
W(\xi)\lesssim Q_{\Gamma_{d}^{+}}(\xi) \sum_{m=1,2,3} W_{m}(\xi),
\end{equation}
where
\begin{align*}
Q_{\Gamma_{d}^{+}}(\xi)&:=\prod_{j\in \Gamma_{d}^{+}} \xi^{-\alpha_{j}}_{j},\\
W_{m}(\xi)&:=\sum_{\substack{
I_{1}\cup I_{2}\cup I_{3}=\Gamma^{-}_{d},\\
I_{j}\cap I_{k}=\emptyset,\; j\neq k}}
w_{I_{1}, I_{2}, I_{3}}(\xi)
\int_{R_{m}\cap{\Theta}(I_{1},I_{2},I_{3})}
\mathcal{Q}(x_d^{\prime},\xi)
d x^{\prime}_{d},
\end{align*}
and
\begin{equation*}
w_{I_{1}, I_{2}, I_{3}}(\xi):=
\prod_{k\in I_{2}}\xi^{-\alpha_{k}}_{k}
\prod_{n\in I_{3}}
\xi_{n}^{-\alpha_{n}(1-{\overline{\ell}}_{n})}.
\end{equation*}
Our main task is to estimate $W_{m}$, $m=1,2,3$, on $U\times V$.
This is the purpose of Sections \ref{scw1}, \ref{w2sc} and \ref{w3sc},
respectively.
Before we proceed, we introduce some more notations.
\subsection{Some notations and remarks}\label{adobh}
In addition to the functions introduced in Section \ref{pstm1}, we gather some auxiliary functions used in the proof. Then we introduce some more conventions and notations that will come handy. Recall that
$\zeta_{d}=\sqrt{|1-\xi_{d}|}$. For each $I\subseteq \Gamma_{d}^{-}$, we define
\begin{align*}
\psi_{I}(x,\xi)&:=\,\sum_{j\in I}
x_{j}+M_{\Gamma_{d}^{-}\setminus I}(\xi),\\
\eta_{I}(x,\xi_{d})&:=
|\sum_{j\in I} x_{j}-(1-\xi^{
2{\overline{\ell}}_{d}}_{d})|,\\
\Psi_{I}(x,\xi)&:=\,|1-\sum_{j\in I_{3}\setminus I}{
x_{j}}-\sum_{n\in I}{
\xi^{2\overline{\ell}_{n}}_{n}}-\xi^{
2{\overline{\ell}}_{d}}_{d}|,\quad I\subseteq I_{3}.
\end{align*}
 For each partition $\{I_{1},I_{2},I_{3}\}$ of $\Gamma^{-}_{d}$, let
$X_{I_{1}, I_{2}, I_{3}}$, $Y_{I_{1}, I_{2}, I_{3}}$, and
$Z_{I_{1}, I_{2}, I_{3}}$ be the functions given by the integral
of $x^{\prime}_{d}\mapsto \mathcal{Q}(x_d^{\prime},\xi)$ over
$R_{m}\cap{\Theta}(I_{1},I_{2},I_{3})$, for
$m=1,2,3$, respectively. Roughly speaking, we estimate each of these integrals first on the subregions of $U\times V$ generated by decomposing the sets $U_{j}$ with $j\in \Gamma^{-}_{d}$ as follows. Take a partition $\{J^{+},J^{-}\}$ of $\Gamma^{-}_{d}$ and, for each $j\in \Gamma_{d}$, set
\begin{equation*}
\widetilde{U}_{j}=
\begin{cases}
U_{j},& j\in \Gamma^{+}_{d};\\[0.05 cm]
U_{j}\cap  \{\xi_{j}^{\overline{\ell}_{j}}\gg  \zeta_{d}\},& j\in J^{+};\\[0.05 cm]
               U_{j}\cap \{\xi_{j}^{\overline{\ell}_{j}}\lesssim \zeta_{d}\},& j\in J^{-}.
\end{cases}
\end{equation*}
For $r_{*}\in \Gamma_{d}^{-}$, let
\begin{equation*}
R_{r_{*}}:=\{(\xi_{1},\dots,\xi_{d-1})\in \widetilde{U}:\xi_{r_{*}}^{{
\overline{\ell}}_{r_{*}}}\geq\xi_{j}^{{
\overline{\ell}}_{j}},\;\forall j\in \Gamma_{d}^{-}\},
\end{equation*}
where $\widetilde{U}:=\widetilde{U}_{1}\times\dots\times \widetilde{U}_{d-1}$. In Section \ref{scw1}, we obtain a pointwise bound for $X_{I_{1}, I_{2}, I_{3}}$ on $\widetilde{U}\times V$. We find it more convenient however to
estimate $Y_{I_{1}, I_{2}, I_{3}}$ first on the subregion
$R_{r_{*}}\times V$ (see Section \ref{w2sc}).\smallskip
\par Given a partition $\{I_{1},I_{2},I_{3}\}$ of $\Gamma^{-}_{d}$,
and assuming that $\xi\in R_{r_{*}}$, we let $j_{*}$, $k_{*}$ and $n_{*}$ denote the integers in
$I_{1}$, $I_{2}$ and $I_{3}$, respectively, such that
\begin{align*}
\xi^{\overline{\ell}_{j_{*}}}_{j_{*}}&\geq\,
\xi^{\overline{\ell}_{j}}_{j},\quad
\forall\, j\in I_{1},\\[0.05 cm]
\xi^{\overline{\ell}_{k_{*}}}_{k_{*}}&\geq\,
\xi^{\overline{\ell}_{k}}_{k},\quad
\forall\, k\in I_{2},\\[0.05 cm]
\xi^{\overline{\ell}_{n_{*}}}_{n_{*}}&\geq\,
\xi^{\overline{\ell}_{n}}_{n},\quad
\forall\, n\in I_{3}.
\end{align*}
In addition, we let $j_{1}^{*}$, $j_{2}^{*}$, and $k^{\prime}_{*}$ denote the integers in $I_{1}\cap E_{1}$, $I_{1}\cap E_{2}$, and $I_{2}$, respectively, such that
\begin{align*}
\xi_{j_{1}^{*}}&\geq\,
\xi_{j},\quad
\forall\, j\in I_{1}\cap E_{1},\\[0.05 cm]
\xi^{\overline{\ell}_{j_{2}^{*}}}_{j_{2}^{*}}&\geq\,
\xi^{\overline{\ell}_{j}}_{j},\quad
\forall\, j\in I_{1}\cap E_{2},\\[0.05 cm]
\xi^{\overline{\ell}_{k^{\prime}_{*}}}_{k^{\prime}_{*}}&\geq\,
\xi^{\overline{\ell}_{k}}_{k},\quad
\forall\, k\in I_{2}\setminus\{k_{*}\}.
\end{align*}
In particular, we have
\begin{equation*}
\lambda_{I_{3}}(\xi)\,\approx\,
\max{\{\xi_{k_{*}}^{2{\overline{\ell}}_{k_{*}}},\,
\xi^{2\overline{\ell}_{j_{*}}}_{j_{*}}\}},\quad
\lambda_{I_{3}\cup \{k_{*}\}}\,\approx\,
\max{\{
    \xi_{k_{*}^{\prime}}^{2{\overline{\ell}}_{k_{*}^{\prime}}},\,
\xi^{2\overline{\ell}_{j_{*}}}_{j_{*}}\}}.
\end{equation*}
We shall use $\Gamma^{-}_{k, d}$
to denote $\Gamma^{-}_{d}\setminus \{k\}$, for $k\in \Gamma^{-}_{d}$.
\smallskip\par
To obtain a pointwise bound for $Z_{I_{1}, I_{2}, I_{3}}$, we further localize to a subregion of $R_{r_{*}}\times V$. Let $I_{*}$ denote either $I_{3}$ or $I_{3}\cup \{k_{*}\}$.
We consider subregions of $R_{r_{*}}$ where, for each $j\in \Gamma^{-}_{d}$,
either
\begin{equation*}
\xi_{j}^{2\overline{\ell}_{j}}\gg
\max{\{\Phi_{I_{3}}(\xi),\lambda_{I_{*}}(\xi)\}},\;\;\text{or}\;\;
\xi_{j}^{2\overline{\ell}_{j}}\lesssim
\max{\{\Phi_{I_{3}}(\xi),\lambda_{I_{*}}(\xi)\}}.
\end{equation*}
More precisely, let
\begin{equation}\label{ukk}
\mathcal{U}_{j}(I_{*}):=
\begin{cases}
\widetilde{U}_{j}\cap \{\xi_{j}^{2\overline{\ell}_{j}}\gtrsim
\max{\{\Phi_{I_{3}}(\xi),\lambda_{I_{*}}(\xi)\}}\},&
j\in \Gamma^{-}_{d}\cap
K^{+}(I_{*});\\[0.05 cm]
\widetilde{U}_{j}\cap \{\xi_{j}^{2\overline{\ell}_{j}}\ll
\max{\{\Phi_{I_{3}}(\xi),\lambda_{I_{*}}(\xi)\}}\},&
j\in \Gamma^{-}_{d}\cap K^{-}(I_{*});\\[0.05 cm]
U_{j},&j\in \Gamma^{+}_{d}.
\end{cases}
\end{equation}
where $\{K^{+}(I_{*}), K^{-}(I_{*})\}$ is a partition of $\Gamma^{-}_{d}$. Obviously, the region $U$ is the union of the regions
\begin{equation*}
\mathcal{U}(I_{*}):=\mathcal{U}_{1}(I_{*})\times\dots\times \mathcal{U}_{d-1}(I_{*})
\end{equation*}
over all partitions $\{K^{+}(I_{*}), K^{-}(I_{*})\}$
of $\Gamma^{-}_{d}$.
The region $U$ is decomposed using (\ref{ukk}) in two different ways using two generally independent sets of partitions, according to whether $I_{*}=I_{3}$ or $I_{*}=I_{3}\cup\{k_{*}\}$.
We estimate $Z_{I_{1}, I_{2}, I_{3}}$ first on $S_{r_{*}}(I_{*})\times V$,
where
\begin{equation*}
S_{r_{*}}(I_{*}):=
\{(\xi_{1},\dots,\xi_{d-1})\in \mathcal{U}(I_{*}):\xi_{r_{*}}^{{
\overline{\ell}}_{r_{*}}}\geq\xi_{j}^{{
\overline{\ell}}_{j}},\;\forall j\in \Gamma_{d}^{-}\}.
\end{equation*}
\par
Finally we introduce some notations that will be used heavily in
Sections \ref{str4} and \ref{lstp}. Recall that, for $\sigma \in (-2,0]$ and $ I\subseteq \Gamma_{d}$, we define
$\theta_{I}(\sigma)=|I|-
\sum_{j\in I}\alpha_{j}+\sigma/2$.
Suppose that
$l=\{n_{1},\dots,n_{N}\}\subseteq \Gamma_{d}$ and that
\begin{equation*}
\xi_{n_{1}}^{{
\overline{\ell}}_{n_{1}}}\geq\dots\geq
\xi_{n_{N}}^{{
\overline{\ell}}_{n_{N}}}.
\end{equation*}
Since $\alpha_{j}<1$ for all $j\in \Gamma$, $\theta_{I}$ has the properties that
if $I\subsetneq I^{\prime}$, then $\theta_{I}(\sigma)<
\theta_{I^{\prime}}(\sigma)$, and if
$I\cap I^{\prime}=\emptyset$, then
$\theta_{I}(2\theta_{I^{\prime}}(\sigma))=
\theta_{I\cup I^{\prime}}(\sigma)$.
By the former property, either $1-\alpha_{n_{1}}+\sigma/2\geq 0$, or there exists a
maximal nonempty subset
$l^{-}(\sigma)=\{n_{1},n_{2},\dots,n_{N^{-}(\sigma)}\}$ of $l$ such that
\begin{equation*}
\theta_{l^{-}(\sigma)}(s)<0,\quad
\forall -2<s\leq \sigma<0.
\end{equation*}
If $1-\alpha_{n_{1}}+\sigma/2\geq 0$,
we set $l^{-}(\sigma)=\emptyset$, $l^{+}(\sigma)=\{n_{1}\}$, and $N_{k}(\sigma)=k$,
$k=1,\dots,N$.
If $\theta_{l}(\sigma)< 0$, then $l^{-}(\sigma)=l$
(i.e., $N^{-}(\sigma)=N$) and we set $l^{+}(\sigma)=\emptyset$.
In the case that $1-\alpha_{n_{1}}+\sigma/2<0$
and $\theta_{l}(\sigma)\geq 0$, we clearly have
$l^{-}(\sigma)
\subsetneq l$ and $l^{-}(\sigma)\neq \emptyset$,
and we denote $n_{N^{-}(\sigma)+k}$
by $n_{N_{k}(\sigma)}$ and the set
$\{n_{1},\dots,n_{N^{-}(\sigma)},n_{N_{1}(\sigma)}\}$
by $l^{+}(\sigma)$. In this notation, we have
\begin{equation*}
l\setminus l^{-}(\sigma)\,=\,
\{n_{N_{1}(\sigma)},n_{N_{2}(\sigma)},\dots,n_{N}\},\quad
l\setminus l^{+}(\sigma)\,=\,
\{n_{N_{2}(\sigma)},\dots,n_{N}\},
\end{equation*}
and we are in the region
\begin{equation*}
\xi_{n_{1}}^{{
\overline{\ell}}_{n_{1}}}\geq\dots\geq\xi_{n_{N^{-}(\sigma)}}^{{
\overline{\ell}}_{n_{N^{-}(\sigma)}}}\geq
\xi_{n_{N_{1}(\sigma)}}^{{
\overline{\ell}}_{n_{N_{1}(\sigma)}}}\geq
\xi_{n_{N_{2}(\sigma)}}^{{
\overline{\ell}}_{n_{N_{2}(\sigma)}}}\geq
\dots\geq
\xi_{n_{N}}^{{
\overline{\ell}}_{n_{N}}},
\end{equation*}
with the obvious convention if $N^{-}(\sigma)=N-1$.\smallskip\par
\smallskip\par
Let us start with $W_{1}(\xi)$.
\subsection{An estimate for \texorpdfstring{
{\mathversion{bold}$W_{1}(\xi)$}}{} on
\texorpdfstring{
{\mathversion{bold}$U\times V$}}{}}\label{scw1}
We are going to show that
\begin{equation}\label{hr1}
W_{1}(\xi) \lesssim \sum_{I\subseteq \Gamma^{-}_{d}}
\sum_{J\subseteq \Gamma^{-}_{d}\setminus I}
M^{\beta_{|\Gamma^{+}_{d}|}(I)+|J|}_{\Gamma^{-}_{d}\setminus J}(\xi)
\,P_{\Gamma^{-}_{d}\setminus ((I\cap E_{1})\cup J)}(\xi)
\,Q^{+}_{J}(\xi).
\end{equation}
\smallskip\par
Let $x_{d}^{\prime}\in R_{1}$. Then
$x_{1}+\dots+x_{d-1}\lesssim \zeta_{d}$. Thus, by (\ref{xdr123}), one has
\begin{multline*}
\begin{aligned}
X_{I_{1}, I_{2}, I_{3}}(\xi)&=\int_{R_{1}\cap{\Theta}(I_{1},I_{2},I_{3})}
\mathcal{Q}(x_d^{\prime},\xi)
d x^{\prime}_{d}\\
&\lesssim
\zeta_{d}^{-2\alpha_{d}}
\prod_{\ell\in \Gamma^{+}_{d}}\int_{0<x_{\ell}< \zeta_{d}}dx_{\ell}
 \prod_{j\in I_{1}}
\int_{\xi_{j}^{\overline{\ell}_{j}}<x_{j}<\zeta_{d}}
x^{-\alpha_{j}\ell_{j}}_{j} dx_{j}
\end{aligned}\\
\prod_{k\in I_{2}}\int_{0<x_{k}<\min{\{\xi_{k}^{\overline{\ell}_{k}},
\zeta_{d}\}}}dx_{k}
\prod_{n\in I_{3}}\int_{\xi_{n}^{\overline{\ell}_{n}}/2<
x_{n}<2\xi_{n}^{\overline{\ell}_{n}}}
|x_{n}-\xi^{\overline{\ell}_{n}}_{n}|^{-\alpha_{n}}dx_{n},
\end{multline*}
and, if $I_{1}\cup I_{3}\neq \emptyset$, we are necessarily in the region
\begin{equation}\label{rgnstr}
\zeta_{d}\gtrsim \sum_{j\in I_{1}\cup I_{3}}\xi_{j}^{\overline{\ell}_{j}}.
\end{equation}
Straightforward computation shows that
\begin{equation}\label{xi123}
X_{I_{1}, I_{2}, I_{3}}(\xi)\lesssim
\zeta_{d}^{|\Gamma^{+}_{d}|-2\alpha_{d}}
\prod_{j\in I_{1}}\varsigma_{j}(\xi)
\prod_{k\in I_{2}}
\min{\{\xi_{k}^{\overline{\ell}_{k}},\zeta_{d}\}}
\prod_{n\in I_{3}}
\xi^{\overline{\ell}_{n}(1-\alpha_{n})}_{n},
\end{equation}
where
\begin{equation*}
\varsigma_{j}(\xi):=
\begin{cases}
    \zeta_{d}^{1-\alpha_{j}},& j\in  E_{1}; \\[0.05 cm]
    \xi_{j}^{\overline{\ell}_{j}-\alpha_{j}},& j\in E_{2}.
\end{cases}
\end{equation*}
We have taken into account that $1-\alpha_{j}\ell_{j}>0$, if $j\in E_{1}$ (i.e., if $\ell_{j}=1$), and $1-\alpha_{j}\ell_{j}<0$, if $j\in E_{2}$ (i.e., if $\ell_{j}>1$), when we computed the $x_{j}$-integrals with $j\in I_{1}$. Recall that $\xi_{j}$ lives in $(0,2)\cap \{\xi_{j}^{\overline{\ell}_{j}}\lesssim \zeta_{d}\}$ if $j\in J^{-}$, and in $(0,2)\cap \{\xi_{j}^{\overline{\ell}_{j}}\gg  \zeta_{d}\}$ if $j\in J^{+}$. The restriction (\ref{rgnstr}) implies  $I_{1}\cup I_{3}\subseteq J^{-}$ (equivalently
$J^{+}\subseteq I_{2}$). Thus
\begin{equation*}
\sum_{j\in \Gamma_{d}^{-}\setminus (I_{2}\cap J^{+})}
\xi_{j}^{\overline{\ell}_{j}}
\lesssim
\zeta_{d},\qquad
\min{\{\xi_{k}^{\overline{\ell}_{k}},\zeta_{d}\}}\lesssim \begin{cases}
\xi_{k}^{\overline{\ell}_{k}},&k\in I_{2}\cap J^{-};
\\[0.2 cm]
\zeta_{d},&k\in I_{2}\cap J^{+}.
\end{cases}
\end{equation*}
With these remarks in mind, we have
\begin{equation}\label{miz}
\zeta_{d}\approx M_{\Gamma_{d}^{-}\setminus (I_{2}\cap J^{+})}\ll \xi_{k}^{\overline{\ell}_{k}},\quad \forall\, k\in I_{2}\cap J^{+},
\end{equation}
and (\ref{xi123}) implies that
\begin{multline*}
W_{1}(\xi)=\,\sum_{\substack{
I_{1}\cup I_{2}\cup I_{3}=\Gamma^{-}_{d},\\
I_{j}\cap I_{k}=\emptyset,\; j\neq k}}
w_{I_{1}, I_{2}, I_{3}}(\xi)X_{I_{1}, I_{2}, I_{3}}(\xi)\\
 \lesssim\,
\sum_{\substack{
I_{1}\cup I_{2}\cup I_{3}=\Gamma^{-}_{d},\\
I_{j}\cap I_{k}=\emptyset,\; j\neq k,\\
I_{1}\cup I_{3}\subseteq J^{-}}}
M_{\Gamma_{d}^{-}\setminus (I_{2}\cap J^{+})}^{\beta_{|\Gamma^{+}_{d}|}(I_{1})+|J^{+}|}(\xi)
 \prod_{j\in (I_{1}\cap E_{2})\cup (I_{2}\cap J^{-})\cup I_{3}}\xi_{j}^{\overline{\ell}_{j}-\alpha_{j}}
\prod_{k\in I_{2}\cap J^{+}}\xi_{k}^{-\alpha_{k}}.
\end{multline*}
Together with (\ref{miz}), the latter estimate proves (\ref{hr1}). We turn to $W_{2}$.
\subsection{An estimate for \texorpdfstring{
{\mathversion{bold}$W_{2}(\xi)$}}{} on
\texorpdfstring{
{\mathversion{bold}$U\times V$}}{}}\label{w2sc}
For $x_{d}^{\prime}\in R_{2}$, we have
$\sum_{\kappa\in \Gamma_{d}}x_{\kappa}\approx \sum_{\kappa\in \Gamma_{d}}x_{\kappa}+\zeta_{d}$.
Using (\ref{xdr123}) again, we see that
\begin{equation}\label{r20}
Y_{I_{1}, I_{2}, I_{3}}(\xi)\lesssim
\int_{\substack{0<x_{\ell}<1,\:\ell\in \Gamma^{+}_{d},\\
2\xi_{j}^{\overline{\ell}_{j}}<x_{j}<1,\:j\in I_{1},\\ 0<x_{k}<\xi_{k}^{\overline{\ell}_{k}}/2,\:k\in I_{2},\\ \xi_{n}^{\overline{\ell}_{n}}/2<
x_{n}<2\xi_{n}^{\overline{\ell}_{n}}\!,\:n\in I_{3},\\
x_{1}+\dots+x_{d-1}\gtrsim \zeta_{d}}}
\frac{\prod_{j\in I_{1}} x^{-\alpha_{j}\ell_{j}}_{j}
\prod_{n\in I_{3}} |x_{n}-\xi^{\overline{\ell}_{n}}_{n}|^{-\alpha_{n}}}
{(\sum_{\kappa\in \Gamma_{d}}x_{\kappa}+\zeta_{d})^{2\alpha_{d}}}
d x^{\prime}_{d}.
\end{equation}
Suppose that $\Gamma^{+}_{d}\neq \emptyset$. Integrating
in $x_{\ell}$, for each $\ell\in \Gamma^{+}_{d}$, we obtain
\begin{subequations}\label{w2g111}
\begin{numcases}
{\int_{0<x_{\ell}<1,\, \ell \in \Gamma^{+}_{d}}\,
(\sum_{\kappa\in \Gamma_{d}}x_{\kappa}+\zeta_{d})^{-2\alpha_{d}}
\prod_{\ell \in \Gamma^{+}_{d}} d x_{\ell}\approx}\label{w2g11}1,&$|\Gamma^{+}_{d}|\geq 2$;\\[0.1 cm]
\label{w2g1}\psi_{\Gamma^{-}_{d}}^{1-2\alpha_{d}}
(x,\xi),&
$|\Gamma^{+}_{d}|=1$.
\end{numcases}
\end{subequations}
\par
The implicit constant in (\ref{w2g111}) is independent of $x$ and $\xi$. Using (\ref{w2g11}) in (\ref{r20}),
it follows that, if $|\Gamma^{+}_{d}|\geq 2$,
\begin{equation}\label{yips}
Y_{I_{1}, I_{2}, I_{3}}(\xi) \lesssim
\prod_{j\in I_{1}\cap E_{2}}
\xi_{j}^{\overline{\ell}_{j}-\alpha_{j}}
\prod_{k\in I_{2}}
\xi_{k}^{\overline{\ell}_{k}}
\prod_{n\in I_{3}}
\xi_{n}^{\overline{\ell}_{n}(1-\alpha_{n})},
\end{equation}
from which we deduce that
\begin{equation}\label{w2pls2}
W_{2}(\xi) \lesssim \sum_{I\subseteq \Gamma_{d}^{-}}
P_{\Gamma^{-}_{d}\setminus (I\cap E_{1})}(\xi),\quad \text{if}\;\;
|\Gamma^{+}_{d}|\geq 2.
\end{equation}
\par Suppose now that $\Gamma^{+}_{d}$
is either empty or a singleton. In order to estimate $Y_{I_{1}, I_{2}, I_{3}}$, we estimate the iterated integrals in $x_{j}$, $j\in \Gamma^{-}_{d}$. If $|\Gamma^{+}_{d}|=1$,
we plug the estimate (\ref{w2g1})
into (\ref{r20}) first.
We start with the integrals in $x_{j}$ with $j\in I_{1}\cap E_{1}$. Let $\eta\geq 0$ and let $s$ be a real number. Then
\begin{subequations}\label{st47}
\begin{numcases}
{\int_{2\xi_{j}<t<1,\,t+\eta\gtrsim
\zeta_{d}}
\frac{t^{-\alpha_{j}} d t}{
(t+\eta+\zeta_{d})^{s}}
\lesssim}
|\log{(\xi_{j}+\eta+\zeta_{d})}|,&$s=1-\alpha_{j}$;\label{w2i12}\\[0.2 cm]
(\xi_{j}+\eta+\zeta_{d})^{1-\alpha_{j}-s},
&$s>1-\alpha_{j}$.\label{w2i13}
\end{numcases}
\end{subequations}
For brevity, we slightly abuse notations and use $\beta_{m}$ to refer to the number $\beta_{m}(I_{1})$, for $m=0,1$. \medskip
\par
\textbullet\ If $\beta_{m}>0$, we have the integral
\begin{equation}\label{tldps}
\int_{2\xi_{j}<x_{j}<1,\,
j\in I_{1}\cap E_{1}}
\psi_{\Gamma^{-}_{d}}^{m-2\alpha_{d}}
\prod_{j\in I_{1}\cap E_{1}} x^{-\alpha_{j}}_{j} d x_{j}
\end{equation}
bounded by a constant independent of $x$ and $\xi$.
Computing the remaining integrals in (\ref{r20}), we find that, if $|\Gamma^{+}_{d}|\leq 1$ and $\beta_{|\Gamma^{+}_{d}|}(I_{1})>0$, the estimate (\ref{yips}) holds.\smallskip\par
\textbullet\ Assume $\beta_{m}=0$.
Using the estimates (\ref{st47}), one has the integral (\ref{tldps})
bounded\footnote{\thanks{Throughout, we say $X$ is bounded or majorized by $Y$ if $X\lesssim Y$.}} by (see the proof of the slightly more involved estimate (\ref{ti1e12}) in Section \ref{471lgs})
\begin{equation*}
|\log{\psi_{\Gamma^{-}_{d}\setminus (I_{1}\cap E_{1})}}| \lesssim |\log{\psi_{I_{2}}}|,
\end{equation*}
on $\{2\xi_{j}^{\overline{\ell}_{j}}<x_{j}<1,\;j\in I_{1}\cap E_{2},\;\;\xi_{n}^{\overline{\ell}_{n}}/2<
x_{n}<2\xi_{n}^{\overline{\ell}_{n}},\:n\in I_{3}\}$.
By a direct computation, it follows from (\ref{r20}) that
\begin{equation}\label{w2lg}
Y_{I_{1}, I_{2}, I_{3}}(\xi) \lesssim
|\log{M_{\Gamma^{-}_{d}}(\xi)}|
\prod_{j\in I_{1}\cap E_{2}}
\xi_{j}^{\overline{\ell}_{j}-\alpha_{j}}
\prod_{k\in I_{2}}
\xi_{k}^{\overline{\ell}_{k}}
\prod_{n\in I_{3}}
\xi_{n}^{\overline{\ell}_{n}(1-\alpha_{n})},
\end{equation}
whenever $|\Gamma^{+}_{d}|\leq 1$ and $\beta_{|\Gamma^{+}_{d}|}(I_{1})=0$.
\smallskip \par
\textbullet\ Suppose that $\beta_{m}< 0$. Then, applying the estimate (\ref{w2i13}) repeatedly to the integrals in $x_{j}$, $j\in I_{1}\cap E_{1}$, we see that
the integral (\ref{tldps}) is bounded by
$\psi_{\Gamma^{-}_{d}\setminus (I_{1}\cap E_{1})}^{\beta_{m}}$. We integrate next in $x_{j}$, for all $j\in I_{1}\cap E_{2}$, employing the estimate
\begin{equation*}
\int_{2\xi_{j}^{\overline{\ell}_{j}}< t<1}t^{-\alpha_{j}\ell_{j}}
(t+\eta+\zeta_{d})^{\beta_{m}}d t
\lesssim \xi_{j}^{\overline{\ell}_{j}-\alpha_{j}}
(\xi_{j}^{\overline{\ell}_{j}}+
\eta+\zeta_{d})^{\beta_{m}},
\end{equation*}
to see that
\begin{equation}\label{r2i1}
\int_{2\xi_{j}^{\overline{\ell}_{j}}<x_{j}<1,\,
j\in I_{1}\cap E_{2}}\,  \psi_{\Gamma^{-}_{d}\setminus (I_{1}\cap E_{1})}^{\beta_{m}}
\prod_{j\in I_{1}\cap E_{2}}x^{-\alpha_{j}\ell_{j}}_{j}d x_{j} \lesssim
\psi_{\Gamma^{-}_{d}\setminus I_{1}}^{\beta_{m}}\prod_{j\in I_{1}\cap E_{2}}\xi_{j}^{\overline{\ell}_{j}-\alpha_{j}}.
\end{equation}
Postpone the integrals in $x_{k}$ with $k\in I_{2}$,
and turn to the integrals in $x_{n}$ with $n\in I_{3}$. If $|I_{3}|\geq 1$, we apply the estimate
\begin{equation*}
\int_{\xi_{n}^{\overline{\ell}_{n}}/2<
t<2\xi_{n}^{\overline{\ell}_{n}}}
|t-\xi^{\overline{\ell}_{n}}_{n}|^{-\alpha_{n}}
(t+\eta+\zeta_{d})^{\beta_{m}}d t\approx
\xi_{n}^{\overline{\ell}_{n}(1-\alpha_{n})}
(\xi_{n}^{\overline{\ell}_{n}}+\eta+\zeta_{d})^{\beta_{m}},
\end{equation*}
successively for each $n\in I_{3}$, to obtain
\begin{equation}\label{r2i3}
\int_{\xi_{n}^{\overline{\ell}_{n}}/2<
x_{n}<2\xi_{n}^{\overline{\ell}_{n}}\!,\,n\in I_{3}}
\psi_{\Gamma^{-}_{d}\setminus I_{1}}^{\beta_{m}}
\prod_{n\in I_{3}}
|x_{n}-\xi^{\overline{\ell}_{n}}_{n}|^{-\alpha_{n}} d x_{n}\\
\lesssim
\psi_{I_{2}}^{\beta_{m}}\prod_{n\in I_{3}}
\xi_{n}^{\overline{\ell}_{n}(1-\alpha_{n})}.
\end{equation}
It follows from the inequalities (\ref{r2i1}) and (\ref{r2i3})
that
\begin{equation}\label{qqw2}
Y_{I_{1}, I_{2}, I_{3}}(\xi)\lesssim
\widetilde{Y}_{I_{1}, I_{2}, I_{3}}(\xi)
\prod_{j\in I_{1}\cap E_{2}}
\xi_{j}^{\overline{\ell}_{j}-\alpha_{j}}
\prod_{n\in I_{3}}
\xi_{n}^{\overline{\ell}_{n}(1-\alpha_{n})},
\end{equation}
where
\begin{equation*}
\widetilde{Y}_{I_{1}, I_{2}, I_{3}}(\xi):=
\int_{0<x_{k}<
\xi_{k}^{\overline{\ell}_{k}},\;k\in I_{2}}
\psi_{I_{2}}^{\beta_{m}}(x,\xi)
\prod_{k\in I_{2}}dx_{k}.
\end{equation*}
Let us estimate $Y_{I_{1}, I_{2}, I_{3}}$ on $R_{r_{*}}\times V$ where we have $\xi_{r_{*}}^{\overline{\ell}_{r_{*}}}\geq
\xi_{j}^{\overline{\ell}_{j}}$, for all $j\in \Gamma^{-}_{d}$.
If $r_{*}\in I_{2}$ and $I_{2}\cap J^{+}\neq \emptyset$, then we necessarily have $r_{*}\in I_{2}\cap J^{+}$.
So we have two possibilities:
\begin{description}
  \item[\normalfont{(i)}] $r_{*}\in I_{1}\cup I_{3}$,
or $r_{*}\in I_{2}$, $I_{2}\cap J^{+}=\emptyset$,
  \item[\normalfont{(ii)}]
$r_{*}\in I_{2}\cap J^{+}$.
\end{description}
\smallskip\par
Consider the case \normalfont{(i)} and let $k\in I_{2}$. If
$r_{*}\in I_{1}\cup I_{3}$, then
$\xi_{k}^{\overline{\ell}_{k}}\leq \sum_{j\in \Gamma_{d}^{-}\setminus I_{2}} \xi_{j}^{\overline{\ell}_{j}}$, while if
$I_{2}\cap J^{+}=\emptyset$,
then $I_{2}\subseteq J^{-}$ and
we have $\xi_{k}^{\overline{\ell}_{k}}\lesssim
\zeta_{d}$. Hence, we have $\sum_{k\in I_{2}}
\xi_{k}^{\overline{\ell}_{k}}
\lesssim M_{\Gamma^{-}_{d}\setminus I_{2}}(\xi)$. Since  $0<x_{k}<
\xi_{k}^{\overline{\ell}_{k}}$ for all $k\in I_{2}$,
\begin{equation*}
\psi_{I_{2}}(x,\xi)\approx
M_{\Gamma^{-}_{d}}(\xi).
\end{equation*}
Thus, we have
\begin{equation}\label{r2i4}
\widetilde{Y}_{I_{1}, I_{2}, I_{3}}(\xi)\lesssim
M^{\beta_{m}}_{\Gamma^{-}_{d}}(\xi)
\prod_{k\in I_{2}}\xi_{k}^{\overline{\ell}_{k}}.
\end{equation}
\smallskip\par
In the case \normalfont{(ii)}, one has
$k_{*}=r_{*}\in I_{2}\cap J^{+}$. Hence $\xi_{r_{*}}^{\overline{\ell}_{r_{*}}}\gtrsim
M_{\Gamma^{-}_{d}\setminus \{r_{*}\}}(\xi)$.
Since $\beta_{1}\geq 1-2\alpha_{d}>-1$, it is impossible
to have $1+\beta_{m}\leq 0$ unless $m=0$.
A straightforward computation yields
\begin{subequations}\label{bmm1p}
\begin{numcases}
{\int_{0<x_{r_{*}}<
\xi_{r_{*}}^{\overline{\ell}_{r_{*}}}}
\psi_{I_{2}}^{\beta_{m}}
dx_{r_{*}}\lesssim}
\psi_{I_{2}\setminus \{r_{*}\}}^{\beta_{m}+1},&$\beta_{m}+1>0$; \label{bmm1p1}\\[0.05 cm]
|\log{\psi^{\!*}_{I_{2}}}|,&$m=0$, $\beta_{0}+1=0$;
\label{bmm1p2}\\[0.05 cm]
\psi^{\!*\,\beta_{0}+1}_{I_{2}},&$m=0$, $\beta_{0}+1<0$
\label{bmm1p3},
\end{numcases}
\end{subequations}
where
\begin{equation*}
\psi^{\!*}_{I_{2}}(x,\xi):=
\sum_{k\in I_{2}\setminus \{r_{*}\}} x_{k}+
M_{\Gamma_{d}^{-}\setminus I_{2}}(\xi).
\end{equation*}
\par\textbullet\ Let $\beta_{m}+1>0$. We have
\begin{equation*}
\sup_{0<x_{k}<\xi_{k}^{\overline{\ell}_{k}},\,
k\in I_{2}\setminus \{r_{*}\}}{\psi_{I_{2}\setminus \{r_{*}\}}^{\beta_{m}+1}(x,\xi)}=\,
\,M^{\beta_{m}+1}_{\Gamma^{-}_{d}}(\xi)\,\approx\,
M^{\beta_{m}}_{\Gamma^{-}_{d}}(\xi)\,\xi_{r_{*}}^{\overline{\ell}_{r_{*}}}.
\end{equation*}
It follows therefore from (\ref{bmm1p1}) that
the estimate (\ref{r2i4}) holds in the case \normalfont{(ii)}
as well, if $\beta_{m}+1>0$.
\smallskip\par
\textbullet\ Now let $\beta_{0}+1=0$. Recall that $\Gamma^{-}_{r_{*}, d}$ denotes $\Gamma^{-}_{d}\setminus \{r_{*}\}$. Using (\ref{bmm1p2})
and the fact that, for small $M_{\Gamma^{-}_{r_{*}, d}}(\xi)$,
\begin{equation*}
\int_{0<x_{k}<
\xi_{k}^{\overline{\ell}_{k}},\;k\in I_{2}}
|\log{\psi^{\!*}_{I_{2}}}|
\prod_{k\in I_{2}\setminus\{r_{*}\}}dx_{k}\\
\approx
|\log{M_{\Gamma^{-}_{r_{*}, d}}(\xi)}|\prod_{k\in I_{2}\setminus\{r_{*}\}}\xi_{k}^{\overline{\ell}_{k}},
\end{equation*}
we obtain
\begin{equation}\label{r2i4bm1}
\widetilde{Y}_{I_{1}, I_{2}, I_{3}}(\xi)\lesssim |\log{M_{\Gamma^{-}_{r_{*}, d}}(\xi)}|\prod_{k\in I_{2}\setminus\{r_{*}\}}\xi_{k}^{\overline{\ell}_{k}}.
\end{equation}
\par
\textbullet\ Finally assume that $\beta_{0}+1<0$. We are under the assumption that $r_{*}=k_{*}$. By definition, $\xi_{k_{*}}^{\overline{\ell}_{k_{*}}}\geq
\xi_{k_{*}^{\prime}}^{\overline{\ell}_{k_{*}^{\prime}}}\geq
\xi_{k}^{\overline{\ell}_{k}}$, for all $k\in I_{2}\setminus \{k_{*},k_{*}^{\prime}\}$. Since $\beta_{0}+2\geq 2-2\alpha_{d}>0$, integrating
both sides of (\ref{bmm1p3}) in $x_{k_{*}^{\prime}}$
and then taking the supremum over $R_{*}=\{0<x_{k}<
\xi_{k}^{\overline{\ell}_{k}},\, k\in I_{2}\setminus\{
r_{*},k_{*}^{\prime}\}\}$, we obtain
\begin{align*}
\sup_{R_{*}}\int_{0<x_{k_{*}^{\prime}}<
\xi_{k_{*}^{\prime}}^{\overline{\ell}_{k_{*}^{\prime}}}}
\psi^{\!*\,\beta_{0}+1}_{I_{2}}
dx_{k_{*}^{\prime}}
&\lesssim
\begin{cases}
M^{\beta_{0}+1}_{\Gamma^{-}_{d}\setminus I_{2}}(\xi)\,
\xi_{k_{*}^{\prime}}^{\overline{\ell}_{k_{*}^{\prime}}},&
\xi_{k_{*}^{\prime}}^{\overline{\ell}_{k_{*}^{\prime}}}\ll
M_{\Gamma^{-}_{d}\setminus I_{2}}(\xi);\\[0.05 cm]
M_{\Gamma^{-}_{r_{*}, d}}^{\beta_{0}+2}(\xi),&
\xi_{k_{*}^{\prime}}^{\overline{\ell}_{k_{*}^{\prime}}}\gtrsim
M_{\Gamma^{-}_{d}\setminus I_{2}}(\xi)
\end{cases}\\[0.05 cm]
&\approx
M_{\Gamma^{-}_{r_{*}, d}}^{\beta_{0}+1}(\xi)\,
\xi_{k_{*}^{\prime}}^{\overline{\ell}_{k_{*}^{\prime}}}.
\end{align*}
Integrating both sides of the latter estimate
in $x_{k}$, $k\in I_{2}\setminus\{r_{*},k_{*}^{\prime}\}$, implies
\begin{equation}\label{b0ng}
\widetilde{Y}_{I_{1}, I_{2}, I_{3}}(\xi)\lesssim
M_{\Gamma^{-}_{r_{*}, d}}^{\beta_{0}+1}(\xi)
\prod_{k\in I_{2}\setminus\{r_{*}\}}
\xi_{k}^{\overline{\ell}_{k}}.
\end{equation}
\par
We have so far shown that, when $|\Gamma^{+}_{d}|\leq 1$, $\widetilde{Y}_{I_{1}, I_{2}, I_{3}}$ satisfies the estimate (\ref{r2i4}) if $r_{*}\in \Gamma^{-}_{d}\setminus I_{2}$, or
$r_{*}\in I_{2}$, $I_{2}\cap J^{+}=\emptyset$,
or $r_{*}\in I_{2}\cap J^{+}$, $\beta_{m}+1>0$.
If $m=0$ and $r_{*}\in I_{2}\cap J^{+}$, then
$\widetilde{Y}_{I_{1}, I_{2}, I_{3}}$ satisfies the estimate (\ref{r2i4bm1}) when
$\beta_{0}+1=0$, and satisfies (\ref{b0ng}) when $\beta_{0}+1<0$.
Combining the estimates (\ref{yips}), (\ref{w2lg}), and the estimates that follow by plugging each of the inequalities (\ref{r2i4}), (\ref{r2i4bm1}) and (\ref{b0ng}) into (\ref{qqw2}), we conclude that, if $\Gamma^{+}_{d}=\emptyset$, we have
\begin{multline}\label{w2gm0}
W_{2}(\xi)\lesssim
\sum_{\substack{I\subseteq\Gamma^{-}_{d} \cap J^{-}\\
\beta_{0}(I)>-2}}
\mathcal{M}_{I,\Gamma_{d}^{-}}(\xi;\beta_{0}(I))\\+
\sum_{j\in \Gamma_{d}^{-}}
\sum_{\substack{I\subseteq\Gamma^{-}_{j, d} \cap J^{-}\\
-2<\beta_{0}(I)\leq -1}}
Q^{+}_{\{j\}}(\xi)
\mathcal{M}_{I,\Gamma^{-}_{j, d}}(\xi;\beta_{0}(I)+1),
\end{multline}
where, for each $I\subseteq J\subseteq \Gamma_{d}^{-}$,
\begin{equation*}
\mathcal{M}_{I,J}(\xi;\gamma):=
\left(\mathds{1}_{\gamma>0}+\mathds{1}_{\gamma=0}|\log{M_{J}(\xi)}|+
\mathds{1}_{-2<\gamma<0}
M_{J}^{\gamma}(\xi)\right)
P_{J\setminus (I\cap E_{1})}(\xi).
\end{equation*}
While if $|\Gamma^{+}_{d}|=1$, it follows from (\ref{yips}), (\ref{w2lg}), and (\ref{qqw2}) combined with (\ref{r2i4}), that
\begin{equation}\label{w2gm1}
W_{2}(\xi)\lesssim
\sum_{\substack{I\subseteq\Gamma^{-}_{d} \cap J^{-}\\
\beta_{1}(I)>-1}}\mathcal{M}_{I,\Gamma^{-}_{d}}(\xi;\beta_{1}(I)).
\end{equation}
Putting together (\ref{w2pls2}), (\ref{w2gm0}) and (\ref{w2gm1}), we deduce that
\begin{multline}\label{w2f}
W_{2}(\xi)\lesssim
\sum_{\substack{I\subseteq\Gamma^{-}_{d} \cap J^{-}\\
\beta_{|\Gamma_{d}^{+}|}(I)>-2}}
\mathcal{M}_{I,\Gamma_{d}^{-}}(\xi;\beta_{|\Gamma_{d}^{+}|}(I))\\+
\sum_{j\in \Gamma_{d}^{-}}
\sum_{\substack{I\subseteq\Gamma^{-}_{j, d} \cap J^{-}\\
-2<\beta_{|\Gamma_{d}^{+}|}(I)\leq -1}}
Q^{+}_{\{j\}}(\xi)\,
\mathcal{M}_{I,\Gamma^{-}_{j, d}}(\xi;\beta_{|\Gamma_{d}^{+}|}(I)+1).
\end{multline}
\par The next task is to estimate $W_{3}$ on $U\times V$. This is the last term to estimate in the sum (\ref{smsm1}).
\subsection{An estimate for
\texorpdfstring{{\mathversion{bold}$W_{3}(\xi)$}}{} on
\texorpdfstring{
{\mathversion{bold}$U\times V$}}{}}\label{w3sc}
First we introduce some notations, obtain some
preliminary estimates and get the case
$|\Gamma^{+}_{d}|\geq 2$ out of the way.
\subsubsection{Preliminary estimates}\label{prmlst}
Recall that
\begin{equation}\label{r3s}
W_{3}(\xi)=
\sum_{\substack{
I_{1}\cup I_{2}\cup I_{3}=\Gamma^{-}_{d},\\
I_{j}\cap I_{k}=\emptyset,\;j\neq k}}
\prod_{k\in I_{2}} \xi^{-\alpha_{k}}_{k}
\prod_{n\in I_{3}}\xi_{n}^{-\alpha_{n}(1-{\overline{\ell}}_{n})}
\,Z_{I_{1}, I_{2}, I_{3}}(\xi),
\end{equation}
and
\begin{align*}
Z_{I_{1}, I_{2}, I_{3}}(\xi)&=
\int_{R_{3}\cap{\Theta}(I_{1},I_{2},I_{3})}
\mathcal{Q}(x_d^{\prime},\xi)
d x^{\prime}_{d}\\
&\lesssim \int_{\substack{
0<x_{\ell}<1,\:\ell\in \Gamma^{+}_{d},\\
2\xi_{j}^{\overline{\ell}_{j}}<x_{j}<1,\:j\in I_{1},\\ 0<x_{k}<\xi_{k}^{\overline{\ell}_{k}}/2,\:k\in I_{2},\\ \xi_{n}^{\overline{\ell}_{n}}/2<
x_{n}<2\xi_{n}^{\overline{\ell}_{n}}\!,\:n\in I_{3},\\
x^2_{1}+\dots+x^2_{d-1}\approx \zeta_{d}^{2}}}
\frac{
\prod_{j\in I_{1}} x^{-\alpha_{j}\ell_{j}}_{j} \prod_{n\in I_{3}}
|x_{n}-\xi^{\overline{\ell}_{n}}_{n}|^{-\alpha_{n}}}
{|x_{1}^{2}+\dots+x_{d-1}^{2}-(1-\xi^{
2{\overline{\ell}}_{d}}_{d})|^{\alpha_{d}}}
dx^{\prime}_{d}.
\end{align*}
It is useful to introduce here the functions
\begin{eqnarray*}
A_{j}(t,\xi)&:=&{
t^{-(1+\alpha_{j}\ell_{j})/2}}
 {\chi_{\{4 {\xi}^{2\overline{\ell}_{j}}_{j}<t\lesssim\zeta_{d}^{2}\}}(t)},\quad
j\in I_{1},\\
B_{k}(t,\xi)&:=& {t^{-1/2}}{\chi_{\{0< 4 t< \xi_{k}^{2\overline{\ell}_{k}}\!,
\;t\lesssim\zeta_{d}^{2}\}}(t)}, \quad k\in I_{2},\\
C_{n}(t,\xi)&:=&
|t-\xi^{2\overline{\ell}_{n}}_{n}|^{-\alpha_{n}}
{\chi_{\{\frac{1}{4} \xi_{n}^{2\overline{\ell}_{n}}<
t< 4 \xi_{n}^{2\overline{\ell}_{n}}\!,
\;t\lesssim\zeta_{d}^{2}\}}(t)}, \quad n\in I_{3}.
\end{eqnarray*}
Recall also that
\begin{equation*}
\eta_{I}(x,\xi_{d}):=
|\sum_{j\in I} x_{j}-(1-\xi^{
2{\overline{\ell}}_{d}}_{d})|,\quad I\subseteq \Gamma_{d},
\end{equation*}
Changing variables $x_j^{2}\rightarrow x_j$, $j\in \Gamma_{d}$, then using the fact that
$|\sqrt{t}-s|\approx s^{-1}|t-s^{2}|$, when $s,t>0$ are such that $t\approx s^2$, we get
\begin{multline}\label{kximn}
Z_{I_{1}, I_{2}, I_{3}}(\xi)\approx
\prod_{n\in I_{3}} \xi^{(\alpha_{n}-1){\overline{\ell}}_{n}}_{n}
\int_{x_{1}+\dots+x_{d-1}\approx \zeta_{d}^{2}}
\eta_{\Gamma_{d}}^{-\alpha_{d}}(x,\xi_{d})
\prod_{\ell\in \Gamma^{+}_{d}} x^{-1/2}_{\ell}\, d x_{\ell}\\
\prod_{j\in I_{1}} A_{j}(x_{j},\xi) dx_{j}
\prod_{k\in I_{2}} B_{k}(x_{k},\xi) dx_{k}
 \prod_{n\in I_{3}}C_{n}(x_{n},\xi) dx_{n}.
\end{multline}
Observe that the restriction
\begin{equation}\label{dots}
x_{1}+\dots+x_{d-1}\approx\zeta_{d}^{2}
\end{equation}
implies that $x_j \lesssim \zeta_{d}^{2}$, for each $j\in \Gamma_{d}$, and $x_j \approx \zeta_{d}^{2}$, for at least one $j\in \Gamma_{d}$.
So it puts us in the region (\ref{rgnstr}) where
$\xi^{\overline{\ell}_{j}}_{j}\lesssim \zeta_{d}$, for all
$j\in I_{1}\cup I_{3}$. In other words, $I_{1}\cap J^{+}=I_{3}\cap J^{+}=\emptyset$, which is equivalent to
 $I_{2}\cap J^{+}=\Gamma^{-}_{d}\cap J^{+}=J^{+}$.
Another noteworthy consequence of
(\ref{dots}) is that
\begin{equation}\label{unetzd}
\eta_{I}(x,\xi_{d})\lesssim \zeta_{d}^{2},
\end{equation}
uniformly in $x$, for every $I\subseteq \Gamma_{d}$. Note also that, by the triangle inequality
\begin{equation*}
\eta_{I}(x,\xi_{d})\geq
|x_{j}-\eta_{I\setminus \{j\}}(x,\xi_{d})|,\quad j\in I,\;I\subseteq \Gamma_{d}.
\end{equation*}
For brevity, we shall write $\eta_{I}$,
$A_{j}(\cdot)$, $B_{k}(\cdot)$, and $C_{n}(\cdot)$,
for $\eta_{I}(x,\xi_{d})$, $A_{j}(\cdot,\xi)$, $B_{k}(\cdot,\xi)$, and $C_{n}(\cdot,\xi)$, respectively.
 \smallskip \par If $\Gamma^{+}_{d}\neq \emptyset$, we may start with integration in $x_{\ell}$ for all $\ell \in \Gamma^{+}_{d}$. One has
\begin{equation}\label{ltr1}
\int_{0<t\lesssim \zeta_{d}^{2}}
t^{-1/2} |t-\eta|^{-s} d t\lesssim
\begin{cases}
\eta^{1/2-s},& s>1/2;\\[0.2 cm]
\zeta_{d}^{1-2s},& s<1/2,
\end{cases}
\end{equation}
for any number $s<1$ and any $\eta>0$ such that $\eta\lesssim \zeta_{d}^2$. Since $\alpha_{d}\in (1/2,1)$, and since
$|1-\xi^{2{\overline{\ell}}_{d}}_{d}|\approx \zeta_{d}^{2}$,
successive application of the estimate (\ref{ltr1}) yields
\begin{equation}\label{lplsk}
\int_{\substack{ x_{\ell}\lesssim \zeta_{d}^{2},\;\ell\in \Gamma^{+}_{d},\\
\eta_{I}\lesssim \zeta_{d}^2,\;I\subseteq \Gamma_{d}}}
\eta_{\Gamma_{d}}^{-\alpha_{d}}\prod_{\ell\in \Gamma^{+}_{d}} x^{-1/2}_{\ell}\,dx_{\ell}
\lesssim
\begin{cases}
\eta_{\Gamma^{-}_{d}}^{1/2-\alpha_{d}},&
|\Gamma^{+}_{d}|=1;\\[0.2 cm]
\zeta_{d}^{|\Gamma^{+}_{d}|-2\alpha_{d}},&|\Gamma^{+}_{d}|\geq 2.
\end{cases}
\end{equation}
Hence, if $|\Gamma^{+}_{d}|\geq 2$, computing the remaining integrals in (\ref{kximn}) (the integrals in $x_{j}$ with $j\in \Gamma^{-}_{d}$) becomes straightforward. We see then that
$Z_{I_{1}, I_{2}, I_{3}}(\xi)$ has, modulo a multiplicative constant, the same pointwise bound as that of
$X_{I_{1}, I_{2}, I_{3}}(\xi)$ given by
(\ref{xi123}).
Inserting that into (\ref{r3s}) and taking into account that we are in the region (\ref{rgnstr}), we see that
$W_{3}$ has the bound (\ref{hr1}), when
$|\Gamma^{+}_{d}|\geq 2$. \smallskip\par
We proceed then to estimate the integral $Z_{I_{1}, I_{2}, I_{3}}$ under the assumption that $\Gamma^{+}_{d}$ is either empty or a singleton. By (\ref{kximn}) and (\ref{lplsk}), we have
\begin{equation}\label{kximnm}
Z_{I_{1}, I_{2}, I_{3}}(\xi)\lesssim
\mathcal{K}_{I_{1}, I_{2}, I_{3}}(\xi)\prod_{n\in I_{3}} \xi^{(\alpha_{n}-1){\overline{\ell}}_{n}}_{n},
\end{equation}
when $|\Gamma^{+}_{d}|=0,1$, where
\begin{equation*}
\mathcal{K}_{I_{1}, I_{2}, I_{3}}(\xi):=
\int \eta_{\Gamma^{-}_{d}}^{|\Gamma^{+}_{d}|/2-\alpha_{d}}
\prod_{j\in I_{1}} A_{j}(x_{j}) dx_{j}
\prod_{k\in I_{2}} B_{k}(x_{k}) dx_{k}
 \prod_{n\in I_{3}}C_{n}(x_{n}) dx_{n}.
\end{equation*}
Using (\ref{kximnm}) in (\ref{r3s}), we see that
\begin{equation}\label{r3s1}
W_{3}(\xi)\lesssim
\sum_{\substack{
I_{1}\cup I_{2}\cup I_{3}=\Gamma^{-}_{d},\\
I_{j}\cap I_{k}=\emptyset,\;j\neq k}}\,
\prod_{k\in I_{2}} \xi^{-\alpha_{k}}_{k}
\prod_{n\in I_{3}}\xi_{n}^{2\alpha_{n}{\overline{\ell}}_{n}-
\alpha_{n}-{\overline{\ell}}_{n}}
\,{\mathcal{K}}_{I_{1}, I_{2}, I_{3}}(\xi).
\end{equation}
The rest of Section \ref{w3sc} is for the most part devoted to estimating the integral $\mathcal{K}_{I_{1}, I_{2}, I_{3}}$, when
$\Gamma^{+}_{d}$ is either empty or a singleton.
\subsubsection{The integrals \texorpdfstring{
{\mathversion{bold}$\mathcal{K}_{I_{1}, I_{2}, I_{3}}(\xi)$}}{}}\label{intmxi}
Let $\Gamma^{+}_{d}$ be either empty or a singleton. We assume initially that the sets $I_{1}\cap E_{1}$, $I_{1}\cap E_{2}$, $I_{2}$ and $I_{3}$ are all nonempty. It will be clear what estimates $\mathcal{K}_{I_{1}, I_{2}, I_{3}}$ has if one or more of these sets is empty.
Recall that  $\beta_{m}(I)=|I\cap E_{1}|-\sum_{j\in I\cap E_{1}}\alpha_{j}+m-2\alpha_{d}$, $m=0,1$, and note that
$\beta_{1}(I)=\beta_{0}(I)+1$ and $\beta_{0}(I)\geq -2\alpha_{d}>-2$. As we did in Section \ref{w2sc},
we slightly abuse notations and let $\beta_{m}$ denote the particular value $\beta_{m}(I_{1})$, for brevity.
Let
 \begin{align*}
F_{I_{1}, I_{2}, I_{3}}(x,\xi)&:=
\int \eta_{\Gamma^{-}_{d}}^{|\Gamma^{+}_{d}|/2-\alpha_{d}}
\prod_{j\in I_{1}} A_{j}(x_{j}) dx_{j},\\
G_{I_{1}, I_{2}, I_{3}}(x,\xi)&:=\int
F_{I_{1}, I_{2}, I_{3}}(x,\xi)
\prod_{k\in I_{2}} B_{k}(x_{k}) dx_{k},
\end{align*}
so that
\begin{equation*}
\mathcal{K}_{I_{1}, I_{2}, I_{3}}(\xi)=
\int G_{I_{1}, I_{2}, I_{3}}(x,\xi) \prod_{n\in I_{3}}C_{n}(x_{n}) dx_{n}.
\end{equation*}
We estimate $F_{I_{1}, I_{2}, I_{3}}$ in Section \ref{471lgs}.
As shall be seen in Section \ref{471lgs}, estimating
$G_{I_{1}, I_{2}, I_{3}}$ and $\mathcal{K}_{I_{1}, I_{2}, I_{3}}$
is easier to complete in the cases:
\begin{description}
  \item[(a1)] $\beta_{|\Gamma^{+}_{d}|}\geq 0$ (i.e., if $\Gamma^{+}_{d}=\emptyset$ and $\beta_{0}\geq 0$,
or if $|\Gamma^{+}_{d}|=1$ and $\beta_{1}\geq 0$),
  \item[(a2)] $J^{+}\neq \emptyset$, and $-1<\beta_{|\Gamma^{+}_{d}|}<0$,
  \item[(a3)] $\Gamma^{+}_{d}=\emptyset$, $|J^{+}|\geq 2$, and $-2<\beta_{0}<-1$,
\item[(a4)] $\Gamma^{+}_{d}=\emptyset$, $J^{+}\neq \emptyset$, and $\beta_{0}=-1$.
\end{description}
\par \textbullet\ The case \textbf{(a1)} is only possible if $I_{1}\cap E_{1}\neq \emptyset$.
In this case, we show that
\begin{subequations}\label{lge2p}
\begin{numcases}
{\mathcal{K}_{I_{1}, I_{2}, I_{3}}(\xi)\lesssim}\label{lge2p1}
\zeta_{d}^{\beta_{|\Gamma^{+}_{d}|}(I_{1})}
\widetilde{\mathcal{K}}_{I_{1}, I_{2}, I_{3}}(\xi),&
$\beta_{|\Gamma^{+}_{d}|}(I_{1})>0$,\\
\label{lge2p2}
|\log{(\Phi_{I_{3}}(\xi)+\sum_{j\in \Gamma^{-}_{d}}\xi^{
2{\overline{\ell}}_{j}}_{j})}|\,\widetilde{\mathcal{K}}_{I_{1}, I_{2}, I_{3}}(\xi),&
$\beta_{|\Gamma^{+}_{d}|}(I_{1})=0$,
\end{numcases}
\end{subequations}
where
\begin{equation*}
\widetilde{\mathcal{K}}_{I_{1}, I_{2}, I_{3}}(\xi):=\prod_{j\in I_{1}\cap E_{2}}
\xi^{\overline{\ell}_{j}-\alpha_{j}}_{j}
\prod_{k\in I_{2}}
\xi^{\overline{\ell}_{k}}_{k}
\prod_{n\in I_{3}} \xi^{2{\overline{\ell}}_{n}(1-\alpha_{n})}_{n}.
\end{equation*}
Let us denote the sum
(\ref{r3s1}) restricted to $a<\beta_{|\Gamma_{d}^{+}|}(I_{1})<b$ by
$\mathcal{R}_{a<\beta_{|\Gamma_{d}^{+}|}(I_{1})<b}(\xi)$, for any numbers $b>a\geq -2$.
We similarly denote the sum (\ref{hr1}) restricted to $a<\beta_{|\Gamma_{d}^{+}|}(I)<b$ by
$\mathcal{S}_{a<\beta_{|\Gamma_{d}^{+}|}(I)<b}(\xi)$.
The estimates (\ref{lge2p}) show that, wherever
$|\Gamma_{d}^{+}|\leq 1$,
\begin{subequations}\label{lge}
\begin{align}
\label{lge121}\mathcal{R}_{\beta_{|\Gamma^{+}_{d}|}(I_{1})>0}(\xi)&\lesssim
\sum_{\substack{I\subseteq \Gamma^{-}_{d}\cap J^{-}\\
\beta_{|\Gamma^{+}_{d}|}(I)>0}}
M_{\Gamma^{-}_{d}}^{\beta_{|\Gamma^{+}_{d}|}(I)}(\xi)\,
P_{\Gamma^{-}_{d}\setminus (I\cap E_{1})}(\xi):=\widetilde{\mathcal{S}}(\xi),\\
\label{lge2}\mathcal{R}_{\beta_{|\Gamma^{+}_{d}|}(I_{1})=0}(\xi)
&\lesssim
\sum_{\substack{I\subseteq \Gamma^{-}_{d}\cap J^{-}\\
\beta_{|\Gamma^{+}_{d}|}(I)= 0}}
\sum_{J\subseteq (\Gamma^{-}_{d}\setminus I)\cap J^{-}}{
|\log{(\Phi_{J}(\xi)+\sum_{j\in \Gamma^{-}_{d}}\xi^{
2{\overline{\ell}}_{j}}_{j})}|\,
P_{\Gamma^{-}_{d}\setminus (I\cap E_{1})}(\xi)}.
\end{align}
\end{subequations}
\par \textbullet\ In the cases
\textbf{(a2)} and \textbf{(a3)}, we prove that
\begin{equation}\label{ma24}
\mathcal{K}_{I_{1}, I_{2}, I_{3}}(\xi)\lesssim
\widetilde{\mathcal{K}}_{I_{1}, I_{2}\cap J^{-}, I_{3}}(\xi)
\zeta_d^{\beta_{|\Gamma^{+}_{d}|}(I_{1})+|J^{+}|}.
\end{equation}
The fact that the estimate (\ref{ma24}) holds in the case
\textbf{(a2)} means that, in the region that corresponds to $|\Gamma^{+}_{d}|\leq 1$ and $J^{+}\neq\emptyset$, we have
\begin{equation}\label{lq1jnp1}
\mathcal{R}_{-1<\beta_{|\Gamma^{+}_{d}|}(I_{1})<0}(\xi)\lesssim
\mathcal{S}_{-1<\beta_{|\Gamma^{+}_{d}|}(I)<0}(\xi).
\end{equation}
Since (\ref{ma24}) holds
in the case \textbf{(a3)} as well, one has
\begin{equation}\label{lq1jnp2}
\mathcal{R}_{-2<\beta_{0}(I_{1})<-1}(\xi)\lesssim
\mathcal{S}_{-2<\beta_{0}(I)<-1}(\xi),
\end{equation}
in the region that corresponds to $\Gamma^{+}_{d}=\emptyset$ and $|J^{+}|\geq 2$.\medskip \par
\textbullet\
We shall also prove that, in the case \textbf{(a4)}, we have
\begin{equation}\label{lge31}
\mathcal{K}_{I_{1}, I_{2}, I_{3}}(\xi)\lesssim
\widetilde{\mathcal{K}}_{I_{1}, I_{2}\setminus \{k_{*}\}, I_{3}}(\xi)
|\log{(\Phi_{I_{3}}(\xi)+\sum_{j\in \Gamma^{-}_{k_{*}, d}}\xi^{2\overline{\ell}_{j}}_{j})}|,
\vspace{-0.3 cm}
\end{equation}
and we necessarily have $k_{*}\in J^{+}$, $\xi_{k_{*}}^{\overline{\ell}_{k_{*}}}\geq \xi_{j}^{\overline{\ell}_{j}}$, for all $j\in \Gamma_{d}$,
and since $I_{3}\subseteq J^{-}$,
\begin{equation*}
\Phi_{I_{3}}(\xi)\leq
\sum_{n\in I_{3}}
\xi^{2\overline{\ell}_{n}}_{n}+|1-
\xi^{2\overline{\ell}_{d}}_{d}|
\approx \sum_{n\in I_{3}}
\xi^{2\overline{\ell}_{n}}_{n}+\zeta_{d}^{2}\approx \zeta_{d}^{2}\ll
\xi^{2{\overline{\ell}}_{k_{*}}}_{k_{*}}.
\end{equation*}
Using the estimate (\ref{lge31}) in (\ref{r3s1}) implies that, in the region that corresponds to
$\Gamma^{+}_{d}=\emptyset$ and $J^{+}\neq \emptyset$,
\begin{multline}\label{lge3}
\mathcal{R}_{\beta_{0}(I_{1})=-1}(\xi)\lesssim
\sum_{k\in \Gamma_{d}}
\sum_{\substack{I\subseteq \Gamma^{-}_{k, d}\cap J^{-},\\
\beta_{0}(I)=-1}}
\sum_{J\subseteq (\Gamma^{-}_{k, d}\setminus I)\cap J^{-}}
\xi_{k}^{-\alpha_{k}}
\chi_{\{\xi_{k}^{2\overline{\ell}_{k}}\gtrsim
\Phi_{J}(\xi)+\sum_{j\in \Gamma^{-}_{k, d}}\xi^{2\overline{\ell}_{j}}_{j}\}}(\xi_{k})\\
|\log{(\Phi_{J}(\xi)+\sum_{j\in \Gamma^{-}_{k, d}}\xi^{2\overline{\ell}_{j}}_{j})}|\,
P_{\Gamma^{-}_{k, d}\setminus (I\cap E_{1})}(\xi).
\end{multline}
Once we get the easier cases
\textbf{(a1)} - \textbf{(a4)} out of the way, we arrive at the more interesting cases:
\begin{description}
  \item[(b1)] $\Gamma^{+}_{d}=J^{+}=\emptyset$, and $-2<\beta_{0}<0$
  \item[(b2)] $\Gamma^{+}_{d}=\emptyset$, $|J^{+}|=1$, and $-2<\beta_{0}<-1$,
  \item[(b3)] $|\Gamma^{+}_{d}|=1$, $J^{+}=\emptyset$,
and $-1<\beta_{1}<0$.
\end{description}
We estimate $\,G_{I_{1}, I_{2}, I_{3}}$
in the cases \textbf{(b1)} through \textbf{(b3)} in Section \ref{mi2xk}. In section \ref{rdz}, we
estimate $\mathcal{K}_{I_{1}, I_{2}, I_{3}}$ in the region
$\{\sum_{n\in I_{3}}\xi_{n}^{\overline{\ell}_{n}}
\ll \zeta_{d}\}$ in each of theses cases, and show that
\begin{equation}\label{xinsmzt}
\left.
\begin{aligned}
\mathcal{R}_{-2<\beta_{0}(I_{1})<0}(\xi)&\lesssim
\mathcal{S}_{-2<\beta_{0}(I)<0}(\xi), \qquad \Gamma^{+}_{d}=J^{+}=\emptyset,\\
\mathcal{R}_{-2<\beta_{0}(I_{1})<-1}(\xi)&\lesssim
\mathcal{S}_{-2<\beta_{0}(I)<-1}(\xi), \qquad \Gamma^{+}_{d}=\emptyset,\;\;|J^{+}|=1,\\
\mathcal{R}_{-1<\beta_{1}(I_{1})<0}(\xi)&\lesssim
\mathcal{S}_{-1<\beta_{1}(I)<0}(\xi), \qquad |\Gamma^{+}_{d}|=1,\;\;J^{+}=\emptyset,
\end{aligned}\right\}
\end{equation}
in that region.
It will remain to estimate $\mathcal{K}_{I_{1}, I_{2}, I_{3}}$
in the region
$\{\sum_{n\in I_{3}}\xi_{n}^{\overline{\ell}_{n}}
\approx \zeta_{d}\}$. Assuming that
$\sum_{n\in I_{3}}\xi_{n}^{\overline{\ell}_{n}}
\approx \zeta_{d}$, we manage to remove portions of the integral of $G_{I_{1}, I_{2}, I_{3}}$ that lead to the same estimate
$\mathcal{K}_{I_{1}, I_{2}, I_{3}}$ has in the region
$\{\sum_{n\in I_{3}}\xi_{n}^{\overline{\ell}_{n}}
\ll \zeta_{d}\}$. We prove some technical results
in Section \ref{str4} and then we complete the task
of estimating $\mathcal{K}_{I_{1}, I_{2}, I_{3}}$ in the cases \textbf{(b1)} - \textbf{(b3)} in Section \ref{lstp}.
\subsubsection{Estimates for
\texorpdfstring{{\mathversion{bold}$\mathcal{K}_{I_{1}, I_{2}, I_{3}}$}}{}
in the easier cases \textbf{(a1)} - \textbf{(a4)}}\label{471lgs}
We begin with the integral of $\eta_{\Gamma^{-}_{d}}^{|\Gamma^{+}_{d}|/2-\alpha_{d}}
\prod_{j\in I_{1}\cap E_{1}} A_{j}(x_{j})$ in $x_{j}$, $j\in I_{1}\cap E_{1}$. For any positive $\eta$ such that $\eta\lesssim \zeta_{d}^{2}$, one has
\begin{subequations}\label{dmi12s}
\begin{numcases}
{\int_{\xi^{2}_{j}<t\lesssim
\zeta_{d}^{2}}
\frac{t^{-(1+\alpha_{j})/2}}{ |t-\eta|^{s/2}}d t
\lesssim}
\label{dmi2} |\log{\max\{\xi^{2}_{j},\eta\}}|,& $s=1-\alpha_{j}$;
\\[0.2 cm]
\label{dmi3}   \min{\{
\eta^{(1-\alpha_{j}-s)/2},
\xi_{j}^{1-\alpha_{j}-s}\}}, &$1-\alpha_{j}<s<2$.
\end{numcases}
\end{subequations}
As shown below, the following estimates follow from the coarea formula for $\beta_{m}>0$, and from (\ref{dmi12s})
when $\beta_{m}\leq 0$.
\begin{subequations}
\begin{numcases}
{\int \eta_{\Gamma^{-}_{d}}^{m/2-\alpha_{d}}
\prod_{j\in I_{1}\cap E_{1}} A_{j}(x_{j})\,dx_{j}\lesssim}
\label{ti1e13} \zeta_{d}^{\beta_{m}},&
$\beta_{m}>0$;
\\[0.2 cm]
\label{ti1e12}{|\log{\max{\{\eta_{\Gamma^{-}_{d}\setminus (I_{1}\cap E_{1})},
\xi^{2}_{j_{1}^{*}}\}}}|},&
$\beta_{m}=0$;
\\[0.2 cm]
\label{ti1e11}\min{\{
\eta_{\Gamma^{-}_{d}\setminus (I_{1}\cap E_{1})}^{{\beta_{m}}/{2}},\,
\xi_{j_{1}^{*}}^{\beta_{m}}\}},&
$\beta_{m}<0$,
\end{numcases}
\end{subequations}
for $m=0,1$. \par
\textbullet\ First, let $\beta_{m}>0$. Let $d\sigma_{t}$ denote the $(|I_{1}\cap E_{1}|-1)$-dimensional coarea measure on the simplex $\Sigma_{t}:=\{x_{j}>0,\, j\in I_{1}\cap E_{1},\; \sum_{j\in I_{1}\cap E_{1}}x_{j}=t\}$. We have
\begin{multline*}
\hspace{-0.25 cm}
\int \eta_{\Gamma^{-}_{d}}^{m/2-\alpha_{d}}
\prod_{j\in I_{1}\cap E_{1}} A_{j}(x_{j})\,dx_{j}\lesssim
\int_{\substack{x_{j}>0,\, j\in I_{1}\cap E_{1}\\
\sum_{j\in I_{1}\cap E_{1}}x_{j}\lesssim \zeta^{2}_{d}}} \eta_{\Gamma^{-}_{d}}^{m/2-\alpha_{d}}
\prod_{j\in I_{1}\cap E_{1}} x_{j}^{-(1+\alpha_{j})/2}\,dx_{j}\\
\begin{aligned}
&\lesssim \int_{0}^{\zeta^{2}_{d}}
|t-\eta_{\Gamma_{d}^{-}\setminus (I_{1}\cap E_{1})}|^{m/2-\alpha_{d}}
\int_{\Sigma_{t}}
\prod_{j\in I_{1}\cap E_{1}} x_{j}^{-(1+\alpha_{j})/2} d\sigma_{t}(x) dt\\
&\approx \int_{0}^{\zeta^{2}_{d}}
|t-\eta_{\Gamma_{d}^{-}\setminus (I_{1}\cap E_{1})}|^{m/2-\alpha_{d}}\;t^{(|I_{1}\cap E_{1}|-\sum_{j\in I_{1}\cap E_{1}}\alpha_{j})/2-1}dt\approx \zeta_{d}^{\beta_{m}},
\end{aligned}
\end{multline*}
by the uniform estimate (\ref{unetzd}). This shows (\ref{ti1e13}).
 \smallskip \par
\textbullet\ Suppose that $\beta_{m}=0$. Recall that $j_{1}^{*}$ denotes the integer in $I_{1}\cap E_{1}$ for which  $\,\xi_{j_{1}^{*}}\geq \xi_{j}$, for all $j\in I_{1}\cap E_{1}$. Let $I\subsetneq (I_{1}\cap E_{1})\setminus \{j_{1}^{*}\}$, and let $j\in (I_{1}\cap E_{1})\setminus (I\cup \{j_{1}^{*}\})$. Clearly,
$-2<\beta_{m}(I)<\beta_{m}(I\cup\{j\})<\beta_{m}$. Consider the following weaker version of (\ref{dmi3}):
\begin{equation*}
\int_{\xi^{2}_{j}<t\lesssim
\zeta_{d}^{2}}
|t-\eta|^{-s/2}\, t^{-(1+\alpha_{j})/2} d t
\lesssim
\eta^{(1-\alpha_{j}-s)/2}, \quad 1-\alpha_{j}<s<2.
\end{equation*}
Using the triangle inequality and then applying the latter auxiliary estimate with $s=-\beta_{m}(I)$,
and $\eta=\eta_{\Gamma^{-}_{d}\setminus(I\cup \{j\})}$, we readily get
\begin{equation*}
\int \eta\mathstrut^{\beta_{m}(I)/2}_{\Gamma^{-}_{d}\setminus I} A_{j}(x_{j})\,dx_{j}\,\leq\,
\int |x_{j}-\eta_{\Gamma^{-}_{d}\setminus(I\cup \{j\})}|^{\beta_{m}(I)/2} A_{j}(x_{j})\,dx_{j}\,\lesssim\,
\eta\mathstrut^{\beta_{m}(I\cup \{j\})/2}_{\Gamma^{-}_{d}\setminus(I\cup \{j\})}.
\end{equation*}
Iterating this estimate, we find that
\begin{equation}\label{dmi3w}
\int
\eta_{\Gamma^{-}_{d}}^{m/2-\alpha_{d}}
\prod_{j\in (I_{1}\cap E_{1})\setminus \{j_{1}^{*}\}} A_{j}(x_{j})\,dx_{j}\lesssim
\eta_{(\Gamma^{-}_{d}\setminus (I_{1}\cap E_{1}))\cup \{j_{1}^{*}\}}^{\beta_{m}(I_{1}\setminus \{j_{1}^{*}\})/2}.
\end{equation}
We are left then with the integral in $x_{j_{1}^{*}}$ that we estimate using (\ref{dmi2}).
Note that $\beta_{m}(I_{1}\setminus \{j_{1}^{*}\})=-(1-\alpha_{j_{1}^{*}})$,
when $\beta_{m}=0$.
According to (\ref{dmi2}) with $s=1-\alpha_{j_{1}^{*}}$
and $\eta=\eta_{\Gamma^{-}_{d}\setminus (I_{1}\cap E_{1})}$,
\begin{equation*}
\int \eta_{(\Gamma^{-}_{d}\setminus (I_{1}\cap E_{1}))\cup \{j_{1}^{*}\}}^{
-(1-\alpha_{j_{1}^{*}})/2}\,
A_{j_{1}^{*}}(x_{j_{1}^{*}})\,dx_{j_{1}^{*}}
\lesssim
|\log{\max{\{\eta_{\Gamma^{-}_{d}\setminus (I_{1}\cap E_{1})},
\xi^{2}_{j_{1}^{*}}\}}}|.
\end{equation*}
This shows (\ref{ti1e12}). \smallskip \par
\textbullet\
Finally, let $\beta_{m}<0$. Then (\ref{dmi3w}) holds true.
Using (\ref{dmi3}), with $s=-\beta_{m}(I_{1}\setminus \{j_{1}^{*}\})$ and $
\eta=\eta_{\Gamma^{-}_{d}\setminus (I_{1}\cap E_{1})}$,
to estimate the integral of (\ref{dmi3w}) in $x_{j_{1}^{*}}$ gives (\ref{ti1e11}).
\medskip \par Assume the case \textbf{(a1)}. If $\beta_{|\Gamma^{+}_{d}|}(I_{1})>0$, the estimate
(\ref{lge2p1}) follows from (\ref{ti1e13})
and direct computation of the remaining iterated integrals that define $\mathcal{K}_{I_{1}, I_{2}, I_{3}}$. If $\beta_{|\Gamma^{+}_{d}|}(I_{1})=0$, we
start from (\ref{ti1e12}) and then apply the following auxiliary estimates that hold
for $\lambda>0$, small enough $\zeta_{d}$, and $|\eta|\lesssim \zeta^{2}_{d}$:
\begin{subequations}
\begin{align}
\label{lgt1}&\int_{\xi_{j}^{2\overline{\ell}_{j}}<t\lesssim
\zeta_{d}^{2}}
\frac{ |\log{\max{\{|t-\eta|,\lambda\}}}|}{t^{(1+\alpha_{j}\ell_{j})/2}} d t\lesssim
\xi_{j}^{\overline{\ell}_{j}-\alpha_{j}}|\log{
\max{\{\eta,\lambda,\xi_{j}^{2\overline{\ell}_{j}}\}}}|,\\[0.1 cm]
\label{lgt2}&\int_{0<t\lesssim\min{\{\xi_{k}^{2\overline{\ell}_{k}},
\zeta_{d}^2\}}}
|\log{\max{\{|t-\eta|,\lambda\}}}| \frac{d t}{\sqrt{t}}\lesssim
\xi_{k}^{\overline{\ell}_{k}}|\log{\max{\{\eta,\lambda,\xi_{k}^{2{
\overline{\ell}}_{k}}\}}}|,\\[0.1 cm]
\nonumber &\int_{t\approx \xi_{n}^{2\overline{\ell}_{n}}}
\frac{|\log{\max{\{|t-\eta|,\lambda\}}}|}{
|t-\xi^{2\overline{\ell}_{n}}_{n}|^{\alpha_{n}}} dt\\
\label{lgt3}&\hspace{4.5 cm}\lesssim
\xi^{2{\overline{\ell}}_{n}(1-\alpha_{n})}_{n}
|\log{\max{\{|\xi^{2\overline{\ell}_{n}}_{n}-\eta|,\,
\lambda,
\xi^{2\overline{\ell}_{n}}_{n}\}}}|.
\end{align}
\end{subequations}
Use (\ref{lgt1}) repeatedly for the integrals in $x_{j}$, $j\in I_{1}\cap E_{2}$, and then apply (\ref{lgt2}) repeatedly for the integrals in $x_{k}$, $k\in I_{2}$, and finally use
(\ref{lgt3}) repeatedly to estimate the integrals in $x_{n}$, $n\in I_{3}$. This gives (\ref{lge2p2}).
\bigskip \par Suppose now that $\beta_{|\Gamma^{+}_{d}|}(I_{1})< 0$ ( i.e.,
$-2<\beta_{0}<0$, if $\Gamma^{+}_{d}=\emptyset$, and
$-1<\beta_{1}<0$, if $|\Gamma^{+}_{d}|=1$). We proceed to estimate the integral of the inequality (\ref{ti1e11}) in $x_{j}$, $j\in I_{1}\cap E_{2}$, against
$\prod_{j\in I_{1}\cap E_{2}} A_{j}(x_{j})$. For these integrals, we will use the estimate
\begin{equation}\label{stmi1e1}
\int_{\xi^{2 {\overline{\ell}}_{j}}_{j}<t\lesssim
\zeta_{d}^{2},\;\ell_{j}\geq 2}
t^{-(1+\alpha_{j}\ell_{j})/2} |t-\eta|^{-s}d t
\lesssim
\xi_{j}^{\overline{\ell}_{j}-\alpha_{j}}
    \min{\{
\xi_{j}^{-2{\overline{\ell}}_{j}s},\eta^{-s}\}},
\end{equation}
that holds true for all $s\in (0,1)$ and $\eta\lesssim \zeta^{2}_{d}$. Since $0<-{\beta_{m}}/{2}<1$, for $m=0,1$, successive application of the estimate (\ref{stmi1e1}) with $s=-{\beta_{m}}/{2}$, for all integers $j\in I_{1}\cap E_{2}$, yields
\begin{equation*}
\hspace{-0.05 cm}
\int \eta_{\Gamma^{-}_{d}\setminus (I_{1}\cap E_{1})}^{\beta_{m}/2}
\prod_{j\in I_{1}\cap E_{2}} A_{j}(x_{j})\,dx_{j}
\lesssim \min{\{\eta_{\Gamma^{-}_{d}\setminus I_{1}}^{\beta_{m}/2},\,
\xi_{j_{2}^{*}}^{\beta_{m}{{\overline{\ell}}_{j_{2}^{*}}}} \}}\prod_{j\in I_{1}\cap E_{2}}
\xi^{\overline{\ell}_{j}-\alpha_{j}}_{j},
\end{equation*}
for $m=0,1$. Recall that $j_{2}^{*}$ denotes the integer in $I_{1}\cap E_{2}$ for which  $\,
\xi^{\overline{\ell}_{j_{2}^{*}}}_{j_{2}^{*}}
\geq \xi^{\overline{\ell}_{j}}_{j}$, for all $j\in I_{1}\cap E_{2}$. With $m=|\Gamma^{+}_{d}|$, combining this last estimate with (\ref{ti1e11}), we see that
\begin{equation}\label{i1eta}
F_{I_{1}, I_{2}, I_{3}}(x,\xi)\lesssim
\min{\{\eta_{\Gamma^{-}_{d}\setminus I_{1}}^{\beta_{|\Gamma^{+}_{d}|}(I_{1})/2},
\,\xi_{j_{*}}^{\beta_{|\Gamma^{+}_{d}|}(I_{1}){{\overline{\ell}}_{j_{*}}}}\}}
\prod_{j\in I_{1}\cap E_{2}}
\xi^{\overline{\ell}_{j}-\alpha_{j}}_{j},
\end{equation}
whenever $\beta_{|\Gamma^{+}_{d}|}(I_{1})<0$.
 \smallskip \par We estimate next $G_{I_{1}, I_{2}, I_{3}}$ under the assumption that $\beta_{|\Gamma^{+}_{d}|}(I_{1})<0$. Let $k\in I_{2}$. Then
\begin{equation}\label{auxi21}
\int_{t\lesssim\min{\{\xi^{2\overline{\ell}_{k}}_{k},\zeta_{d}^{2}\}}}
|t-\eta|^{-s}{t^{-1/2}}{dt} \lesssim
\dfrac{\min{\{\xi_{k}^{\overline{\ell}_{k}},\zeta_{d}
\}}}{\max{\{
\min{\{\xi_{k}^{2{\overline{\ell}}_{k}s},\zeta^{2s}_{d}
\}},\eta^{s}\}}},
\end{equation}
for any number $ 0<s<1/2$ and $\eta\lesssim \zeta^{2}_{d}$. Either
$k\in I_{2}\cap J^{-}$ (i.e., $\xi_{k}^{\overline{\ell}_{k}}\lesssim \zeta_{d}$) or $k\in J^{+}$ (i.e., $\xi_{k}^{\overline{\ell}_{k}}\gg  \zeta_{d}$).
As pointed out in Section \ref{prmlst}, we
are estimating $\mathcal{K}_{I_{1}, I_{2}, I_{3}}$ under the assumption
that $I_{1}\cup I_{3}\subseteq J^{-}$ and
$J^{+}\subseteq I_{2}$. With this and the uniform estimate (\ref{unetzd}) in mind, apply (\ref{auxi21}) to obtain
\begin{numcases}
{\int\eta_{I}^{-s}
B_{k}(x_{k}) \,dx_{k}\lesssim}
\label{i2slhf1}\zeta_{d}^{1-2s},& \hspace{-0.1 cm}$k\in  J^{+}$;
\\[0.2 cm]
\label{i2slhf20}\xi_{k}^{\overline{\ell}_{k}}
\min{\{
\xi_{k}^{-2{\overline{\ell}}_{k}s},
\eta_{I\setminus \{k\}}^{-s}\}},&
\hspace{-0.1 cm} $k\in I_{2}\cap J^{-}$,
\end{numcases}
for any number $0<s<1/2$, and every $I\subseteq \Gamma^{-}_{d}$
that contains $k$. An immediate consequence of (\ref{i2slhf20})
that will be useful later is that
\begin{equation}\label{i2slhf2}
\int
\min{\{\eta_{I}^{-s},\,
\xi_{j_{*}}^{-2{\overline{\ell}}_{j_{*}}s}\}}
B_{k}(x_{k}) \,dx_{k}\lesssim
\xi_{k}^{\overline{\ell}_{k}}
\min{\{
\eta_{I\setminus \{k\}}^{-s},\,
\xi_{j_{*}}^{-2{\overline{\ell}}_{j_{*}}s},\,
\xi_{k}^{-2{\overline{\ell}}_{k}s}\}},
\end{equation}
for any $0<s<1/2$, and for every
$I\subseteq \Gamma_{d}$ such that $k\in I_{2}\cap J^{-}\cap I$.
\smallskip \par
At this point, we are ready to demonstrate (\ref{ma24}).
Observe first that whenever $J^{+}\neq \emptyset$, one has
$k_{*}\in J^{+}$. Using (\ref{i2slhf1}) with $s=-\beta_{|\Gamma_{d}^{+}|}/2$ in the case \textbf{(a2)}
one obtains
\begin{equation*}
\int \eta_{\Gamma^{-}_{d}\setminus I_{1}}^{\beta_{|\Gamma^{+}_{d}|}(I_{1})/2}
B_{k_{*}}(x_{k_{*}}) dx_{k_{*}}
\lesssim \zeta_{d}^{1+\beta_{|\Gamma^{+}_{d}|}(I_{1})}.
\end{equation*}
We also have
\begin{equation*}
\int \min{\{\eta_{\Gamma^{-}_{d}\setminus I_{1}}^{\beta_{0}/2},\,
\xi_{j_{*}}^{\overline{\ell}_{j_{*}}\beta_{0}}\}}
B_{k_{*}}(x_{k_{*}}) dx_{k_{*}}
\lesssim
\eta^{\beta_{1}/2}_{\Gamma^{-}_{d}\setminus (I_{1}\cup \{k_{*}\})},\quad
-1<\beta_{1}<0.
\end{equation*}
This is a consequence of the estimates
(\ref{i2slhf3}) and (\ref{i2slhf5}) that will be proved later.
In the case \textbf{(a3)}, let
$k^{\prime}\in J^{+}\setminus \{k_{*}\}$, apply the latter estimate
 to the $x_{k_{*}}$-integral,
and then apply (\ref{i2slhf1}) with
$s=-\beta_{1}/2$ to the $x_{k^{\prime}}$-integral, to obtain
\begin{multline*}
\int \min{\{\eta_{\Gamma^{-}_{d}\setminus I_{1}}^{\beta_{0}/2},\,
\xi_{j_{*}}^{\overline{\ell}_{j_{*}}\beta_{0}}\}}
B_{k_{*}}(x_{k_{*}}) B_{k^{\prime}}(x_{k^{\prime}}) dx_{k_{*}}dx_{k^{\prime}}
\\ \lesssim
\int \eta_{\Gamma^{-}_{d}\setminus (I_{1}\cup\{k_{*}\})}^{\beta_{1}/2}
B_{k^{\prime}}(x_{k^{\prime}}) dx_{k^{\prime}}\lesssim \zeta_{d}^{2+\beta_{0}}.
\end{multline*}
Doing the remaining integrals in $x_{k}$,
$k\in I_{2}$ (i.e., the integrals in
$x_{k}$, $k\in I_{2}\setminus\{k_{*}\}$ in the case
\textbf{(a2)}, and
the integrals in $x_{k}$, $k\in I_{2}\setminus\{k_{*},\,k^{\prime}\}$ in the case \textbf{(a3)}), we find that
\begin{multline*}
\int
\min{\{\eta_{\Gamma^{-}_{d}\setminus  I_{1}}^{\beta_{|\Gamma^{+}_{d}|}(I_{1})/2},\,
\xi_{j_{*}}^{\overline{\ell}_{j_{*}}\beta_{|\Gamma^{+}_{d}|}(I_{1})}\}}
\prod_{k\in I_{2}} B_{k}(x_{k}) dx_{k}\\
\begin{aligned}
&\lesssim
\begin{cases}
\zeta_{d}^{1+\beta_{|\Gamma^{+}_{d}|}(I_{1})}\prod_{k\in I_{2}\setminus\{k_{*}\}}
\min{\{\xi^{\overline{\ell}_{k}}_{k},\zeta_{d}\}},&
\textbf{(a2)};
\\[0.2 cm]
\zeta_{d}^{2+\beta_{|\Gamma^{+}_{d}|}(I_{1})}\prod_{k\in I_{2}\setminus\{k_{*},\,k^{\prime}\}}
\min{\{\xi^{\overline{\ell}_{k}}_{k},\zeta_{d}\}},&
\textbf{(a3)}
\end{cases}\\[0.05 cm]
&\lesssim
\zeta_d^{|J^{+}|+\beta_{|\Gamma^{+}_{d}|}(I_{1})}\prod_{k\in I_{2}\cap J^{-}}\xi^{\overline{\ell}_{k}}_{k}.
\end{aligned}
\end{multline*}
Combined with (\ref{i1eta}), this implies
(\ref{ma24}).
\smallskip \par
Assume now the case \textbf{(a4)}.
Again, we must have $k_{*}\in J^{+}$.
Since $I_{1}\cup I_{3}\subseteq J^{-}$, we
have $\xi^{\overline{\ell}_{k_{*}}}_{k_{*}}\geq
\xi^{\overline{\ell}_{j}}_{j}$, for all $j\in \Gamma^{-}_{d}$.
Hence $r_{*}=k_{*}$.
Also, since $\xi^{\overline{\ell}_{k_{*}}}_{k_{*}}\gg \zeta_{d}$,
it follows from (\ref{unetzd}) that $
\xi^{2\overline{\ell}_{k_{*}}}_{k_{*}}\gg
\eta_{\Gamma^{-}_{d}\setminus \left
(I_{1}\cup \{k_{*}\}\right)}$,
uniformly. Verifiably, one has
\begin{align*}
\int
\min{\{\eta_{\Gamma^{-}_{d}\setminus I_{1}}^{-1/2},
\,\xi_{j_{*}}^{-{\overline{\ell}}_{j_{*}}}\}}
B_{k_{*}}(x_{k_{*}}) dx_{k_{*}}
&\lesssim\,
\int_{t\lesssim \zeta_{d}^{2}}
\frac{1}{\sqrt{|t-\eta_{\Gamma^{-}_{d}\setminus (I_{1}\cup \{k_{*}\})}|}+\xi_{j_{*}}^{{\overline{\ell}}_{j_{*}}}}\frac{dt}{\sqrt{t}}\\
&\lesssim\, 1+
\int_{\max{\{\eta_{\Gamma^{-}_{d}\setminus (I_{1}\cup \{k_{*}\})},\, \xi_{j_{*}}^{2{\overline{\ell}}_{j_{*}}}\}}\ll
t\lesssim \zeta_{d}^{2}}
\frac{dt}{t}\\
&\approx\, |\log{\max{\{\eta_{\Gamma^{-}_{d}\setminus (I_{1}\cup \{k_{*}\})},\, \xi_{j_{*}}^{2{\overline{\ell}}_{j_{*}}}\}}}|.
\end{align*}
Combining this with (\ref{i1eta}) gives
\begin{multline*}
\int F_{I_{1}, I_{2}, I_{3}}(x,\xi)B_{k_{*}}(x_{k_{*}}) dx_{k_{*}}
=\int
\eta_{\Gamma^{-}_{d}}^{-\alpha_{d}}
\prod_{j\in I_{1}}A_{j}(x_{j})\,dx_{j} B_{k_{*}}(x_{k_{*}}) dx_{k_{*}}\\ \lesssim
|\log{\max{\{\eta_{\Gamma^{-}_{d}\setminus (I_{1}\cup \{k_{*}\})},\, \xi_{j_{*}}^{2{\overline{\ell}}_{j_{*}}}\}}}|
\prod_{j\in I_{1}\cap E_{2}}
\xi^{\overline{\ell}_{j}-\alpha_{j}}_{j},
\end{multline*}
when $\Gamma^{+}_{d}=\emptyset$ and $\beta_{0}=-1$.
Analogously to the proof of (\ref{lge2}), we can estimate next the integrals
in $x_{k}$, for each $k\in I_{2}\setminus \{k_{*}\}$, using
(\ref{lgt2}), and then estimate the integrals in $x_{n}$,
for each $n\in I_{3}$, using (\ref{lgt3}). This shows (\ref{lge31}).
Our next step is to estimate
$G_{I_{1}, I_{2}, I_{3}}$ in each of the cases \textbf{(b1)} - \textbf{(b3)}.
\subsubsection{Estimates for \texorpdfstring{{\mathversion{bold}$G_{I_{1}, I_{2}, I_{3}}$}}{} in the cases (b1) - (b3)}\label{mi2xk}
In the case \textbf{(b1)}, we prove that
$G_{I_{1}, I_{2}, I_{3}}(x,\xi)$ is bounded pointwise by
\begin{equation}\label{fc1}
\begin{cases}
A(x,\xi)\tau(\xi),&\;\text{if}\;\;
\xi^{2{\overline{\ell}}_{k_{*}}}_{k_{*}}\lesssim
\max{\{\eta_{I_{3}}(x,\xi_{d}),\,\xi_{j_{*}}^{2{\overline{\ell}}_{j_{*}}}\}},\;
-2<\beta_{0}\leq -1;
\\[0.05 cm]
B(x,\xi)\tau_{k_{*}}(\xi),&\;\text{if}\;\;
\xi^{2{\overline{\ell}}_{k_{*}}}_{k_{*}}\gg
\max{\{\eta_{I_{3}}(x,\xi_{d}),\,
\xi_{j_{*}}^{2{\overline{\ell}}_{j_{*}}}\}},\;
-2<\beta_{0}<-1,\\[0.05 cm]
C(x,\xi)\tau_{k_{*}}(\xi),&\;\text{if}\;\;
\xi^{2{\overline{\ell}}_{k_{*}}}_{k_{*}}\gg
\max{\{\eta_{I_{3}}(x,\xi_{d}),\,
\xi_{j_{*}}^{2{\overline{\ell}}_{j_{*}}}\}},\;
\beta_{0}=-1,\\[0.05 cm]
D_{0}(x,\xi)\tau(\xi),&\;\text{if}\;\;-1<\beta_{0}<0,
\end{cases}
\end{equation}
where\footnote{\thanks{See also Section \ref{adobh} for the notations $j_{*}$, $k_{*}$, $k^{\prime}_{*}$, $n_{*}$.}}
\begin{align*}
A(\cdot,\xi)&:=\,
\min{\{\eta_{I_{3}}^{\beta_{0}/2},\,
\xi_{j_{*}}^{{
\overline{\ell}}_{j_{*}}\beta_{0}}\}},\\
B(\cdot,\xi)&:=\,\min{\{
\eta^{\beta_{1}/2}_{I_{3}},\,
\xi_{j_{*}}^{\overline{\ell}_{j_{*}}\beta_{1}}
,\,\xi_{k_{*}^{\prime}}^{\overline{\ell}_{k_{*}^{\prime}}\beta_{1}}
\}},\\
C(\cdot,\xi)&:=\,|\log{\max{\{
\eta_{I_{3}},\,
\xi_{j_{*}}^{2{\overline{\ell}}_{j_{*}}}
,\,\xi_{k_{*}^{\prime}}^{2
{\overline{\ell}}_{k_{*}^{\prime}}}
\}}}|,\\
D_{m}(\cdot,\xi)&:=\,
\min{\{\eta_{I_{3}}^{\beta_{m}/2},\,
\xi_{j_{*}}^{{
\overline{\ell}}_{j_{*}}\beta_{m}},\,
\xi_{k_{*}}^{{
\overline{\ell}}_{k_{*}}\beta_{m}}\}}, \quad m=0,1,\\
\tau(\xi)&:=\,\prod_{j\in I_{1}\cap E_{2}}
\xi^{\overline{\ell}_{j}-\alpha_{j}}_{j}\prod_{k\in I_{2}}
\xi^{\overline{\ell}_{k}}_{k},\\
\tau_{k_{*}}(\xi)&:=\,\prod_{j\in I_{1}\cap E_{2}}\xi^{\overline{\ell}_{j}-\alpha_{j}}_{j}
\prod_{k\in I_{2}\setminus\{k_{*}\}}
\xi^{\overline{\ell}_{k}}_{k}.
\end{align*}
Whereas in the case \textbf{(b2)}, we have
\begin{equation}\label{fc2}
G_{I_{1}, I_{2}, I_{3}}(x,\xi)\lesssim
{B}(x,\xi)\tau_{k_{*}}(\xi),
\end{equation}
and we necessarily have
\begin{equation*}
J^{+}=\{k_{*}\}=\{r_{*}\},\quad
\xi^{2{\overline{\ell}}_{k_{*}}}_{k_{*}}\gg
\max{\{\eta_{I_{3}}(x,\xi_{d}),\,
\xi_{j_{*}}^{2{\overline{\ell}}_{j_{*}}}\}}.
\end{equation*}
Finally, in the case \textbf{(b3)}, we show that
\begin{equation}\label{fc3}
G_{I_{1}, I_{2}, I_{3}}(x,\xi)\lesssim
D_{1}(x,\xi)\tau(\xi).
\end{equation}
We keep in mind that the functions $A$, $B$,
$C$, $D_{0}$, $D_{1}$, $\tau$ and $\tau_{k_{*}}$ are all dependent on the partition $\{I_{1},I_{2},I_{3}\}$. We drop this dependence on the partition
from the notation for simplicity.
We shall make use of the following auxiliary estimates that hold for any
$\lambda>0$ and real $\eta$ such that
$\max{\{|\eta|,\lambda^{2}\}}\lesssim \zeta_{d}^{2}$:\smallskip\par
\begin{subequations}\label{i2c3}
$\displaystyle
\int_{t\lesssim\min{\{\xi^{2\overline{\ell}_{k}}_{k},\zeta_{d}^{2}\}}}
\min{\{|t-\eta|^{-s},\lambda^{-2s}\}}{t^{-1/2}}{dt}$
\begin{numcases}
{\lesssim}\label{i2c31}\xi^{\overline{\ell}_{k}}_{k}
{\min{\{|\eta|^{-s},\lambda^{-2s}\}}},&
\text{if}\;\;$\xi^{\overline{\ell}_{k}}_{k}\lesssim\max{\{\sqrt{|\eta|},
\lambda\}}$,\;\; $1/2\leq s<1$;
\\[0.05 cm]
\label{i2c32}{\min{\{|\eta|^{1/2-s},\lambda^{1-2s}\}}},&
\text{if}\;\;$\xi^{\overline{\ell}_{k}}_{k}\gg \max{\{\sqrt{|\eta|},
\lambda\}}$,\;\; $1/2<s<1$;
\\[0.05 cm]
\label{i2c33}|\log{\max{\{|\eta|,\lambda^{2}\}}}|,&
\text{if}\;\;$\xi^{\overline{\ell}_{k}}_{k}\gg \max{\{\sqrt{|\eta|},
\lambda\}}$,\;\; $s=1/2$,
\end{numcases}
\end{subequations}
along with
\begin{equation}\label{i2c34}
\int_{t<\xi^{2\overline{\ell}_{k}}_{k}/2}
|\log{\max{\{|t-\eta|,\,\lambda^{2}\}}}| t^{-1/2}dt
\lesssim
\xi^{\overline{\ell}_{k}}_{k}
|\log{\max{\{|\eta|,\,\lambda^{2},\,
\xi^{2\overline{\ell}_{k}}_{k}\}}}|.
\end{equation}
\smallskip
\begin{description}
  \item[(b1)] $\Gamma^{+}_{d}=J^{+}=\emptyset$, $-2<\beta_{0}<0$
\end{description}
\smallskip\par  Assume for the moment that
$-2<\beta_{0}\leq -1$.
We estimate the integral of (\ref{i1eta}) against
$B(x_{k_{*}})$ in $x_{k_{*}}$ first, using (\ref{i2c3}). The estimates (\ref{i2slhf3}), (\ref{i2slhf5}) and
(\ref{i2slog}) below follow from
(\ref{i2c31}), (\ref{i2c32}) and (\ref{i2c33}), respectively,
with $s=-\beta_{0}/2$, $\lambda=\xi_{j_{*}}^{\overline{\ell}_{j_{*}}}$, and  $\eta=1-\xi^{
2{\overline{\ell}}_{d}}_{d}-
\sum_{j\in \Gamma^{-}_{d}\setminus (I_{1}\cup \{k_{*}\})}x_{j}$.
We have\medskip
\begin{subequations}\label{hf35lg}
$\displaystyle
\int \min{\{\eta_{\Gamma^{-}_{d}\setminus I_{1}}^{\beta_{0}/2},\,
\xi_{j_{*}}^{\overline{\ell}_{j_{*}}\beta_{0}}\}}
B_{k_{*}}(x_{k_{*}}) dx_{k_{*}}$
\begin{numcases}
{\lesssim}
\xi^{\overline{\ell}_{k_{*}}}_{k_{*}}
\min{\{
{\eta^{\beta_{0}/2}_{\Gamma^{-}_{d}\setminus (I_{1}\cup\{k_{*}\})}},\,
\xi_{j_{*}}^{\overline{\ell}_{j_{*}}\beta_{0}}\}},&
\nonumber \\
&\hspace{-3.5 cm}\text{if}\;
$\xi^{2{\overline{\ell}}_{k_{*}}}_{k_{*}}\lesssim
\max{\{
\eta_{\Gamma^{-}_{d}\setminus (I_{1}\cup\{k_{*}\})},\,
\xi_{j_{*}}^{2{\overline{\ell}}_{j_{*}}}\}}$, $-2<\beta_{0}\leq-1$;
\label{i2slhf3}\\[0.05 cm]
\min{\{
\eta^{\beta_{1}/2}_{\Gamma^{-}_{d}\setminus (I_{1}\cup\{k_{*}\})},\,
\xi_{j_{*}}^{\overline{\ell}_{j_{*}}\beta_{1}}\}},&
\nonumber \\
&\hspace{-3.5 cm}\text{if}\;
$\xi^{2{\overline{\ell}}_{k_{*}}}_{k_{*}}\gg
\max{\{
\eta_{\Gamma^{-}_{d}\setminus (I_{1}\cup\{k_{*}\})},\,
\xi_{j_{*}}^{2{\overline{\ell}}_{j_{*}}}\}}$,
$-2<\beta_{0}<-1$; \label{i2slhf5}\\[0.05 cm]
|\log{\max{\{
\eta_{\Gamma^{-}_{d}\setminus (I_{1}\cup\{k_{*}\})},\,
\xi_{j_{*}}^{2{\overline{\ell}}_{j_{*}}}\}}}|,&
\nonumber \\
&\hspace{-3.5 cm}\text{if}\;
$\xi^{2{\overline{\ell}}_{k_{*}}}_{k_{*}}\gg
\max{\{
\eta_{\Gamma^{-}_{d}\setminus (I_{1}\cup\{k_{*}\})},\,
\xi_{j_{*}}^{2{\overline{\ell}}_{j_{*}}}\}}$, $\beta_{0}=-1$.
\label{i2slog}
\end{numcases}
\end{subequations}
The following remarks are useful when $|I_{2}|\geq 1$.
\begin{remark}\label{rst11}
Assume that $x_{k}\ll \xi^{2\overline{\ell}_{k}}_{k}$,
$\forall k\in I_{2}$. Then,
for each $I\subseteq I_{2}\setminus \{k_{*}\}$, the restriction
$\xi^{2{\overline{\ell}}_{k_{*}}}_{k_{*}}\gg
\max{\{
\eta_{\Gamma^{-}_{d}\setminus (I_{1}\cup\{k_{*}\})},\,
\xi_{j_{*}}^{2{\overline{\ell}}_{j_{*}}}\}}$
is equivalent to
$\xi^{2{\overline{\ell}}_{k_{*}}}_{k_{*}}\gg
\max{\{
\eta_{I\cup I_{3}},\,
\xi_{j_{*}}^{2{\overline{\ell}}_{j_{*}}}\}}$.
To see this, let $I\subseteq I_{2}\setminus \{k_{*}\}$.
Since $\xi^{\overline{\ell}_{k}}_{k}\leq
\xi^{\overline{\ell}_{k_{*}}}_{k_{*}}$, $\forall k\in I_{2}$,
we have $
|\eta_{\Gamma^{-}_{d}\setminus (I_{1}\cup\{k_{*}\})}-
\eta_{I\cup I_{3}}|\leq
\sum_{k\in I_{2}\setminus (I\cup \{k_{*}\})} x_{k}\ll
\xi^{2{\overline{\ell}}_{k_{*}}}_{k_{*}}$.
\end{remark}
We note an immediate consequence of
Remark \ref{rst11}.
\begin{remark}\label{rst12}
Let $I\subseteq I_{2}\setminus \{k_{*}\}$,
 and  $x_{j}\ll \xi^{2\overline{\ell}_{j}}_{j}$,
$\forall j\in I_{2}$. Then, whenever
$\xi^{2{\overline{\ell}}_{k_{*}}}_{k_{*}}\lesssim
\max{\{
\eta_{\Gamma^{-}_{d}\setminus (I_{1}\cup\{k_{*}\})},
\xi_{j_{*}}^{2{\overline{\ell}}_{j_{*}}}\}}$, we necessarily have
$\,\xi^{2\overline{\ell}_{k}}_{k}\lesssim
\max{\{
\eta_{I\cup I_{3}},
\xi_{j_{*}}^{2{\overline{\ell}}_{j_{*}}}\}}$, for every $k\in I_{2}$.
\end{remark}
\smallskip\par
The inequality (\ref{i2slhf3}) holds if
$\xi^{2{\overline{\ell}}_{k_{*}}}_{k_{*}}\lesssim
\max{\{
\eta_{\Gamma_{d}\setminus (I_{1}\cup\{k_{*}\})},\,
\xi_{j_{*}}^{2{\overline{\ell}}_{j_{*}}}\}}$.
By Remark \ref{rst12}, we necessarily have $
\xi^{2{\overline{\ell}}_{k_{*}^{\prime}}}_{k_{*}^{\prime}}\lesssim
\max{\{
{\eta_{\Gamma_{d}\setminus (I_{1}\cup\{k_{*},
k_{*}^{\prime}\})}},\,
\xi_{j_{*}}^{2{\overline{\ell}}_{j_{*}}}\}}$ in this case.
Thus, we can use (\ref{i2c31}) with $s=-\beta_{0}/2$,
$\lambda=\xi_{j_{*}}^{\overline{\ell}_{j_{*}}}$, and $\eta=1-\xi^{
2{\overline{\ell}}_{d}}_{d}-\sum_{j\in \Gamma^{-}_{d}\setminus
\left(I_{1}\cup\{k_{*},\,k^{\prime}_{*}\}\right)}x_{j}$,
to estimate the
integral of the inequality
(\ref{i2slhf3}) in $x_{k^{\prime}_{*}}$ against
$B_{k^{\prime}_{*}}(x_{k^{\prime}_{*}})$. This gives
(\ref{tmlxi1}) below.
The inequality (\ref{tmlxi3})
follows using (\ref{i2c34}), with the same $\lambda$ and $\eta$, to
estimate the integral of (\ref{i2slog})
in $x_{k_{*}}^{\prime}$. To obtain (\ref{tmlxi2}), apply
(\ref{i2slhf2}) with $s=-\beta_{1}/2$, $I=\Gamma^{-}_{d}\setminus (I_{1}\cup\{k_{*}\})$, and $k=k_{*}^{\prime}$ to the integral of
(\ref{i2slhf5}). So we have \medskip\\
\begin{subequations}
$\displaystyle \iint \min{\{\eta_{\Gamma^{-}_{d}\setminus I_{1}}^{\beta_{0}/2},
\xi_{j_{*}}^{\beta_{0}{\overline{\ell}}_{j_{*}}}\}}
B_{k_{*}}(x_{k_{*}})B_{k^{\prime}_{*}}(x_{k^{\prime}_{*}})
dx_{k_{*}}dx_{k^{\prime}_{*}}$\vspace{-0.2 cm}
\begin{numcases}
{\lesssim}
\nonumber \xi^{\overline{\ell}_{k_{*}}}_{k_{*}}
\xi^{\overline{\ell}_{k^{\prime}_{*}}}_{k^{\prime}_{*}}
\min{\{\eta^{\beta_{0}/2}_{\Gamma^{-}_{d}
\setminus\left(I_{1}\cup\{k_{*},\,k^{\prime}_{*}\}\right)},\,
\xi_{j_{*}}^{\overline{\ell}_{j_{*}}\beta_{0}}\}},&\\
\label{tmlxi1} &
\hspace{-6 cm} if
$\xi^{2\overline{\ell}_{k^{\prime}_{*}}}_{k^{\prime}_{*}}\lesssim
\max{\{\eta_{\Gamma^{-}_{d}
\setminus\left(I_{1}\cup\{k_{*},\,k^{\prime}_{*}\}\right)},\,
\xi_{j_{*}}^{2\overline{\ell}_{j_{*}}}\}}$,
$-2<\beta_{0}\leq -1$;\\
\nonumber
\xi_{k^{\prime}_{*}}^{\overline{\ell}_{k^{\prime}_{*}}}
\min{\{
\eta^{\beta_{1}/2}_{\Gamma^{-}_{d}
\setminus\left(I_{1}\cup\{k_{*},\,k^{\prime}_{*}\}\right)},\,
\xi_{j_{*}}^{\overline{\ell}_{j_{*}}\beta_{1}}
,\,\xi_{k^{\prime}_{*}}^{\overline{\ell}_{k^{\prime}_{*}}\beta_{1}}
\}},&\\
\label{tmlxi2} &
\hspace{-6 cm} if
$\xi^{2\overline{\ell}_{k^{\prime}_{*}}}_{k^{\prime}_{*}}\gg
\max{\{\eta_{\Gamma^{-}_{d}
\setminus\left(I_{1}\cup\{k_{*},\,k^{\prime}_{*}\}\right)},\,
\xi_{j_{*}}^{2\overline{\ell}_{j_{*}}}\}}$,
$-2<\beta_{0}<-1$;\\
\nonumber\xi_{k^{\prime}_{*}}^{\overline{\ell}_{k^{\prime}_{*}}}
|\log{\max{\{
\eta_{\Gamma^{-}_{d}
\setminus\left(I_{1}\cup\{k_{*},\,k^{\prime}_{*}\}\right)},\,
\xi_{j_{*}}^{2{\overline{\ell}}_{j_{*}}}
,\,\xi_{k^{\prime}_{*}}^{2{\overline{\ell}}_{k^{\prime}_{*}}}
\}}}|,&\\
\label{tmlxi3} &
\hspace{-6 cm} if
$\xi^{2\overline{\ell}_{k^{\prime}_{*}}}_{k^{\prime}_{*}}\gg
\max{\{\eta_{\Gamma^{-}_{d}
\setminus\left(I_{1}\cup\{k_{*},\,k^{\prime}_{*}\}\right)},\,
\xi_{j_{*}}^{2\overline{\ell}_{j_{*}}}\}}$,
$\beta_{0}=-1$.
\end{numcases}
\end{subequations}
The estimates (\ref{tmlxi2}) and (\ref{tmlxi3}) hold only if $\xi^{2\overline{\ell}_{k^{\prime}_{*}}}_{k^{\prime}_{*}}\gg
\max{\{\eta_{\Gamma^{-}_{d}
\setminus\left(I_{1}\cup\{k_{*},\,k^{\prime}_{*}\}\right)},\,
\xi_{j_{*}}^{2\overline{\ell}_{j_{*}}}\}}$. This follows from Remark \ref{rst11}
and the fact that
 (\ref{i2slhf5}) and (\ref{i2slog}) hold only if $\xi^{2\overline{\ell}_{k^{\prime}_{*}}}_{k^{\prime}_{*}}\gg
\max{\{\eta_{\Gamma^{-}_{d}
\setminus\left(I_{1}\cup\{k_{*}\}\right)},\,
\xi_{j_{*}}^{2\overline{\ell}_{j_{*}}}\}}$.
With the help of Remark \ref{rst12}, we can
apply (\ref{i2c31}) repeatedly to estimate the integral
of the inequality (\ref{tmlxi1}) in $x_{k}$,
$k\in I_{2}\setminus\{k_{*},k_{*}^{\prime}\}$. This gives (\ref{i2f1}) below.
We can also use (\ref{i2slhf2}) and
(\ref{i2c34}) repeatedly to estimate the integral of
(\ref{tmlxi2}) and (\ref{tmlxi3}), respectively, in $x_{k}$,
$k\in I_{2}\setminus\{k_{*},k_{*}^{\prime}\}$.
This yields (\ref{i2f2}) and (\ref{i2f3}), respectively.
We find that\medskip\\
\begin{subequations}\label{i2f123}
$\displaystyle\int\,
\min{\{\eta_{\Gamma^{-}_{d}\setminus I_{1}}^{\beta_{0}/2},
\xi_{j_{*}}^{\beta_{0}{\overline{\ell}}_{j_{*}}}\}}
\prod_{k\in I_{2}} B_{k}(x_{k}) dx_{k}$
\begin{numcases}
{\lesssim}
\nonumber \min{\{\eta_{I_{3}}^{\beta_{0}/2},\,
\xi_{j_{*}}^{{
\overline{\ell}}_{j_{*}}\beta_{0}}\}}
\prod_{k\in I_{2}}
\xi^{\overline{\ell}_{k}}_{k},&\\
&\hspace{-4 cm}\label{i2f1}if
$\xi^{2{\overline{\ell}}_{k_{*}}}_{k_{*}}\lesssim
\max{\{\eta_{I_{3}},\,\xi_{j_{*}}^{2{
\overline{\ell}}_{j_{*}}}\}}$,
$-2<\beta_{0}\leq -1$;
\\
\nonumber
\min{\{
\eta^{\beta_{1}/2}_{I_{3}},\,
\xi_{j_{*}}^{\overline{\ell}_{j_{*}}\beta_{1}}
,\,\xi_{k_{*}^{\prime}}^{\overline{\ell}_{k_{*}^{\prime}}\beta_{1}}
\}}\prod_{k\in I_{2}\setminus\{k_{*}\}}
\xi^{\overline{\ell}_{k}}_{k},&\\
&\hspace{-4 cm}\label{i2f2}if
$\xi^{2{\overline{\ell}}_{k_{*}}}_{k_{*}}\gg
\max{\{\eta_{I_{3}},\,
\xi_{j_{*}}^{2{\overline{\ell}}_{j_{*}}}\}}$,
$-2<\beta_{0}<-1$,\vspace{0.12 cm}\\
\nonumber |\log{\max{\{
\eta_{I_{3}},\,
\xi_{j_{*}}^{2{\overline{\ell}}_{j_{*}}}
,\,\xi_{k_{*}^{\prime}}^{2
{\overline{\ell}}_{k_{*}^{\prime}}}
\}}}|\prod_{k\in I_{2}\setminus\{k_{*}\}}
\xi^{\overline{\ell}_{k}}_{k},&\\
&\hspace{-4 cm}\label{i2f3}if
$\xi^{2{\overline{\ell}}_{k_{*}}}_{k_{*}}\gg
\max{\{\eta_{I_{3}},\,
\xi_{j_{*}}^{2{\overline{\ell}}_{j_{*}}}\}}$,
$\beta_{0}=-1$.
\end{numcases}
\end{subequations}
The inequalities (\ref{fc1}) follow,
for $-2<\beta_{0}\leq -1$, by combining (\ref{i1eta})
with the estimates (\ref{i2f123}).
The restriction associated with each of the estimates (\ref{i2f123}) comes from Remark \ref{rst11}.
\smallskip\par If $-1<\beta_{0}<0$, one can apply (\ref{i2slhf2}), with
$s=-\beta_{0}/2$, successively to the $x_{k}$-integrals, $k\in I_{2}$. We see then that
\begin{equation*}
\int\,
\min{\{\eta_{\Gamma^{-}_{d}\setminus I_{1}}^{\beta_{0}/2},
\xi_{j_{*}}^{\beta_{0}{\overline{\ell}}_{j_{*}}}\}}
\prod_{k\in I_{2}} B_{k}(x_{k}) dx_{k}\lesssim
\min{\{\eta_{I_{3}}^{\beta_{0}/2},\,
{\xi_{j_{*}}^{{{\overline{\ell}}_{j_{*}}}\beta_{0}}},\
{\xi_{k_{*}}^{{{\overline{\ell}}_{k_{*}}}\beta_{0}}}\}}
\prod_{k\in I_{2}}
\xi^{\overline{\ell}_{k}}_{k},
\end{equation*}
when $-1<\beta_{0}<0$. Combined with (\ref{i1eta}), this completes the proof of  (\ref{fc1}).
\medskip
\begin{description}
  \item[(b2)] $\Gamma^{+}_{d}=\emptyset$, $|J^{+}|=1$, and $-2<\beta_{0}<-1$
\end{description}
\par We are estimating $G_{I_{1}, I_{2}, I_{3}}$ under the assumption that
$I_{1}\cup I_{3} \subseteq J^{-}$ and
$J^{+}\subseteq I_{2}$. So if $J^{+}$ is a singleton,
then $J^{+}=\{k_{*}\}=\{r_{*}\}$. Thus,
\begin{equation*}
\xi^{\overline{\ell}_{k_{*}}}_{k_{*}}\gg \zeta_{d}\gtrsim
\xi^{\overline{\ell}_{j_{*}}}_{j_{*}},
\end{equation*}
and by (\ref{unetzd}),
$\xi^{2\overline{\ell}_{k_{*}}}_{k_{*}}\gg
\eta_{\Gamma^{-}_{d}\setminus I}$, uniformly, for every $I\subseteq \Gamma_{d}$.
In particular, we have
$
\xi^{2\overline{\ell}_{k_{*}}}_{k_{*}}\gg \max{\{
\eta_{\Gamma^{-}_{d}\setminus\left(I_{1}\cup\{k_{*}\}\right)},\,
\xi^{2\overline{\ell}_{j_{*}}}_{j_{*}}\}}$. It follows then from
(\ref{i2slhf5}) that
\begin{equation*}
\int \min{\{\eta_{\Gamma^{-}_{d}\setminus I_{1}}^{\beta_{0}/2},
{\xi_{j_{*}}^{\overline{\ell}_{j_{*}}\beta_{0}}}
\}}
B_{k_{*}}(x_{k_{*}}) dx_{k_{*}}
\lesssim \min{\{
\eta_{\Gamma^{-}_{d}\setminus \left(I_{1}\cup\{k_{*}\}\right)}^{\beta_{1}/2},
{\xi_{j_{*}}^{\overline{\ell}_{j_{*}}\beta_{1}}}
\}}.
\end{equation*}
Start with the latter estimate, and then use  (\ref{i2slhf2})
repeatedly with $s=-\beta_{1}/2$ to
estimate the integrals in $x_{k}$, for all
$k\in I_{2}\setminus \{k_{*}\}$. We deduce that
\begin{equation*}
\int
\min{\{\eta_{\Gamma^{-}_{d}\setminus I_{1}}^{\beta_{0}/2},
{\xi_{j_{*}}^{\overline{\ell}_{j_{*}}\beta_{0}}}
\}}
\prod_{k\in I_{2}} B_{k}(x_{k}) dx_{k}
\lesssim
\min{\{\eta^{\beta_{1}/2}_{I_{3}},\,
\xi_{k_{*}^{\prime}}^{\overline{\ell}_{k_{*}^{\prime}}\beta_{1}},\xi_{j_{*}}^{\overline{\ell}_{j_{*}}\beta_{1}}\}}\prod_{k\in I_{2}\setminus\{k_{*}\}}
\xi^{\overline{\ell}_{k}}_{k}.
\end{equation*}
The estimate (\ref{fc2}) follows from
(\ref{i1eta}) and the latter estimate combined.
\medskip
\begin{description}
  \item[(b3)] $|\Gamma^{+}_{d}|=1$, $J^{+}=\emptyset$,
and $-1<\beta_{1}<0$.
\end{description}
\par
In this case, we use (\ref{i2slhf2}), with
$s=-\beta_{1}/2$, repeatedly
for all $k\in I_{2}$. We see then that
\begin{equation*}
\int
\min{\{\eta_{\Gamma^{-}_{d}\setminus I_{1}}^{\beta_{1}/2},
\xi_{j_{*}}^{\overline{\ell}_{j_{*}}\beta_{1}}\}}
\prod_{k\in I_{2}} B_{k}(x_{k}) dx_{k}
\lesssim
\min{\{\eta_{I_{3}}^{\beta_{1}/2},{\xi_{j_{*}}^{{
\overline{\ell}}_{j_{*}}\beta_{1}}},
{\xi_{k_{*}}^{{{\overline{\ell}}_{k_{*}}}\beta_{1}}}\}}\prod_{k\in I_{2}}
\xi^{\overline{\ell}_{k}}_{k}.
\end{equation*}
Together with (\ref{i1eta}), this last estimate proves
(\ref{fc3}).
\subsubsection{A reduction}\label{rdz}
The last step is to estimate
\begin{equation}\label{mirst}
\mathcal{K}_{I_{1}, I_{2}, I_{3}}(\xi)=
\int G_{I_{1}, I_{2}, I_{3}}(x,\xi) \prod_{n\in I_{3}}C_{n}(x_{n}) dx_{n},
\end{equation}
in each of the cases \textbf{(b1)} -  \textbf{(b3)} starting from  the inequalities
(\ref{fc1}) - (\ref{fc3}), respectively. Recall that
$\Gamma^{-}_{d}=J^{+}\cup J^{-}=I_{1}\cup I_{2}\cup I_{3}$, $\lambda_{I}(\xi)=\sum_{j\in \Gamma^{-}_{d}\setminus I}\xi_{j}^{2{
\overline{\ell}}_{j}}$, for each $I\subseteq \Gamma^{-}_{d}$, and that, by definition,
$\zeta_{d}\gtrsim \xi_{j}^{{
\overline{\ell}}_{j}}$, $\forall j\in J^{-}$,
$\zeta_{d}\ll  \xi_{j}^{{
\overline{\ell}}_{j}}$, $\forall j\in J^{+}$.
We are working under the assumptions that
\begin{equation*}
\begin{aligned}
&\xi_{j_{*}}^{{
\overline{\ell}}_{j_{*}}}\geq
\xi_{j}^{{
\overline{\ell}}_{j}},\;\forall j\in I_{1},\quad
\xi_{k_{*}}^{{
\overline{\ell}}_{k_{*}}}\geq
\xi_{k}^{{
\overline{\ell}}_{k}},\;\forall k\in I_{2},\quad
\xi_{n_{*}}^{{
\overline{\ell}}_{n_{*}}}\geq
\xi_{n}^{{
\overline{\ell}}_{n}},\;\forall n\in I_{3},\\[0.2 cm]
&J^{-}=
\begin{cases}
\Gamma^{-}_{d},& \textbf{(b1)},\;\textbf{(b3)};\\[0.05 cm]
\Gamma^{-}_{d}\setminus\{k_{*}\},& \textbf{(b2)}.
\end{cases}
\end{aligned}
\end{equation*}
If $\zeta_{d}\approx \xi_{n_{*}}^{{
\overline{\ell}}_{n_{*}}}$, these assumptions aimply
\begin{equation}\label{b123}
\begin{cases}
\xi_{n_{*}}^{2{
\overline{\ell}}_{n_{*}}}\gtrsim
\max{\{\xi_{j_{*}}^{2{
\overline{\ell}}_{j_{*}}},\,\xi_{k_{*}}^{2{
\overline{\ell}}_{k_{*}}}\}}\,\approx\,\lambda_{I_{3}}(\xi),&
\textbf{(b1)},\;\textbf{(b3)};\\[0.05 cm]
\xi_{k_{*}}^{2{
\overline{\ell}}_{k_{*}}}\gg
\xi_{n_{*}}^{2{
\overline{\ell}}_{n_{*}}}\gtrsim
\max{\{\xi_{j_{*}}^{2{
\overline{\ell}}_{j_{*}}},\,\xi_{k^{\prime}_{*}}^{2{
\overline{\ell}}_{k^{\prime}_{*}}}\}}\,\approx\,
\lambda_{I_{3}\cup\{k_{*}\}}(\xi),& \textbf{(b2)}.
\end{cases}
\end{equation}
Let us prove the estimates (\ref{xinsmzt}).
Suppose that $\sum_{n\in I_{3}}
\xi_{n}^{\overline{\ell}_{n}}$ $\ll$ $\zeta_{d}$.
Since $\xi_{n}^{2\overline{\ell}_{n}}/4<
x_{n}< 4 \xi_{n}^{2\overline{\ell}_{n}}$, for all $n\in I_{3}$, one has
$\sum_{n\in I_{3}}x_{n}\ll  \zeta_{d}^{2}$, and consequently
$\eta_{I_{3}}\approx \zeta_{d}^{2}$. We then have
\begin{equation}\label{abdo}
A\lesssim \zeta_{d}^{\beta_{0}},\quad
B\lesssim \zeta_{d}^{\beta_{1}},\quad
C\lesssim |\log{\zeta_{d}}|,\quad
D_{m}\lesssim \zeta_{d}^{\beta_{m}},\;\;m=0,1.
\end{equation}
Furthermore, since $J^{+}=\emptyset$ in the case \textbf{(b1)}, the scenario $\xi^{2{\overline{\ell}}_{k_{*}}}_{k_{*}}\gg \eta_{I_{3}}$
is impossible in that case in the region $\sum_{n\in I_{3}}
\xi_{n}^{\overline{\ell}_{n}}\ll \zeta_{d}$, and the estimates (\ref{fc1}) reduce to
\begin{equation*}
G_{I_{1}, I_{2}, I_{3}}(\cdot,\xi)
\lesssim
(A(\cdot,\xi)+D_{0}(\cdot,\xi))\tau(\xi)\lesssim%%
\zeta_{d}^{\beta_{0}}\tau(\xi),\quad
-2<\beta_{0}<0,
\end{equation*}
in the case \textbf{(b1)}. Using
(\ref{abdo}) in the estimates
(\ref{fc2}) and (\ref{fc3}), we obtain
\begin{equation*}
G_{I_{1}, I_{2}, I_{3}}(\cdot,\xi)
\lesssim
\begin{cases}
\zeta_{d}^{\beta_{1}}\tau_{k_{*}}(\xi),& \textbf{(b2)};\\[0.05 cm]
\zeta_{d}^{\beta_{1}}\tau(\xi),& \textbf{(b3)}.
\end{cases}
\end{equation*}
Therefore, in each of the cases \textbf{(b1)} through \textbf{(b3)}, we have the uniform estimate
\begin{equation*}
G_{I_{1}, I_{2}, I_{3}}(\cdot,\xi)\chi_{\{\sum_{n\in I_{3}}\xi_{n}^{\overline{\ell}_{n}}
(\xi)\ll  \zeta_{d}\}}(\xi)\lesssim
\zeta_{d}^{\beta_{|\Gamma^{+}_{d}|}+|J^{+}|}
\prod_{j\in I_{1}\cap E_{2}}
\xi^{\overline{\ell}_{j}-\alpha_{j}}_{j}
\prod_{k\in I_{2}\cap J^{-}}\xi^{\overline{\ell}_{k}}_{k}.
\end{equation*}
Consequently, $\mathcal{K}_{I_{1}, I_{2}, I_{3}}(\xi)\chi_{\{\sum_{n\in I_{3}}\xi_{n}^{\overline{\ell}_{n}}\ll
\zeta_{d}\}}(\xi)$ is bounded by
\begin{equation*}
\zeta_{d}^{\beta_{|\Gamma^{+}_{d}|}+|J^{+}|}
\prod_{j\in I_{1}\cap E_{2}}
\xi^{\overline{\ell}_{j}-\alpha_{j}}_{j}
\prod_{k\in I_{2}\cap J^{-}}\xi^{\overline{\ell}_{k}}_{k}
\prod_{n\in I_{3}}
\xi^{2{\overline{\ell}}_{n}(1-\alpha_{n})}_{n}.
\end{equation*}
Plugging this into (\ref{r3s1}) gives (\ref{xinsmzt}). \smallskip\par
By (\ref{xinsmzt}), it suffices to
 estimate $\mathcal{K}_{I_{1}, I_{2}, I_{3}}$ on the region $\{\sum_{n\in I_{3}}\xi_{n}^{\overline{\ell}_{n}}
\approx \zeta_{d}\}$, or equivalently the region
$\{\xi_{n_{*}}^{\overline{\ell}_{n_{*}}}
\approx \zeta_{d}\}$. Note that, once the estimates (\ref{fc1}) - (\ref{fc3}) are
plugged into (\ref{mirst}), and assuming $\xi_{n_{*}}^{\overline{\ell}_{n_{*}}}
\approx \zeta_{d}$, one can restrict the domain of integration of the integral of $A$ to the region $\Omega^{-}$, the integrals of
$D_{0}$ and $D_{1}$ to $\Omega_{I_{3}}$, and
the integrals of $B$ and $C$ to $\Omega^{+}$, where
\begin{align*}
\Omega_{I_{*}}&:=\,\{\xi_{n_{*}}^{2\overline{\ell}_{n_{*}}}\gg
\max{\{\eta_{I_{3}}(x,\xi_{d}),\lambda_{I_{*}}(\xi)\}}\},\quad
I_{*}=I_{3},\,I_{3}\cup\{k_{*}\},\\
\Omega^{-}&:=\,\{\xi_{n_{*}}^{2\overline{\ell}_{n_{*}}}\gg
\max{\{\eta_{I_{3}}(x,\xi_{d}),\,
\xi_{j_{*}}^{2{\overline{\ell}}_{j_{*}}}\}}\gtrsim \xi^{2{\overline{\ell}}_{k_{*}}}_{k_{*}}\},
\\[0.05 cm]
\Omega^{+}&:=\,\{
\xi^{2{\overline{\ell}}_{k_{*}}}_{k_{*}}\gg
\max{\{\eta_{I_{3}}(x,\xi_{d}),
\xi_{j_{*}}^{2{\overline{\ell}}_{j_{*}}}\}}\}\cap
\Omega_{I_{3}\cup \{k_{*}\}}.
\end{align*}
To see this, suppose that $\xi_{n_{*}}^{\overline{\ell}_{n_{*}}}
\approx \zeta_{d}$. Then we immediately get
\begin{equation*}
\begin{aligned}
&A(\cdot,\xi)\chi_{\{\xi_{n_{*}}^{
\overline{\ell}_{n_{*}}}\approx\,
\max{\{\sqrt{\eta_{I_{3}}},\,\xi_{j_{*}}^{
\overline{\ell}_{j_{*}}}\}}\}} &\approx&\; \zeta_{d}^{\beta_{0}},\\
&B(\cdot,\xi)\chi_{\{\xi_{n_{*}}^{
\overline{\ell}_{n_{*}}}\approx\,
\max{\{\sqrt{\eta_{I_{3}}},\,\xi_{j_{*}}^{
\overline{\ell}_{j_{*}}},\,\xi_{k^{\prime}_{*}}^{
\overline{\ell}_{k^{\prime}_{*}}}\}}\}} &\approx&\; \zeta_{d}^{\beta_{1}},\\
&C(\cdot,\xi)\chi_{\{\xi_{n_{*}}^{
\overline{\ell}_{n_{*}}}\approx\,
\max{\{\sqrt{\eta_{I_{3}}},\,\xi_{j_{*}}^{
\overline{\ell}_{j_{*}}},\,\xi_{k^{\prime}_{*}}^{
\overline{\ell}_{k^{\prime}_{*}}}\}}\}} &\approx&\;  |\log{\zeta_{d}}|,\\
&D_{m}(\cdot,\xi)\chi_{\{\xi_{n_{*}}^{
\overline{\ell}_{n_{*}}}\approx\,
\max{\{\sqrt{\eta_{I_{3}}},\,\xi_{j_{*}}^{
\overline{\ell}_{j_{*}}},\,\xi_{k_{*}}^{
\overline{\ell}_{k_{*}}}\}}\}} &\approx&\; \zeta_{d}^{\beta_{m}},\quad m=0,1.
\end{aligned}
\end{equation*}
On the other hand, if $\xi_{n_{*}}^{\overline{\ell}_{n_{*}}}
\approx \zeta_{d}$, then (\ref{unetzd}) and (\ref{b123})  together imply
\begin{equation*}
\begin{cases}
\xi_{n_{*}}^{2{
\overline{\ell}}_{n_{*}}}\gtrsim
\max{\{\eta_{I_{3}},\,\lambda_{I_{3}}\}},&
\textbf{(b1)},\;\textbf{(b3)};\\[0.05 cm]
\xi_{k_{*}}^{2{
\overline{\ell}}_{k_{*}}}\gg
\xi_{n_{*}}^{2{
\overline{\ell}}_{n_{*}}}\gtrsim
\max{\{\eta_{I_{3}},\,\lambda_{I_{3}\cup\{k_{*}\}}\}},& \textbf{(b2)}.
\end{cases}
\end{equation*}
This justifies the aforementioned cutoffs.
Let us redefine\footnote{\thanks{These functions are introduced after
(\ref{fc1}) in
Section \ref{mi2xk}. We are simply exploiting (\ref{b123}) and cutting off $A$, $D_{m}$, and both $B$ and $C$ by
$\chi_{\Omega^{-}}$, $\chi_{\Omega_{I_{3}}}$
and $\chi_{\Omega^{+}}$, respectively.}}
\begin{eqnarray*}
A(\cdot,\xi)&:=&
\min{\{\eta_{I_{3}}^{\beta_{0}/2},\,
\xi_{j_{*}}^{{
\overline{\ell}}_{j_{*}}\beta_{0}}\}}\,\chi_{\Omega^{-}},\\
B(\cdot,\xi)&:=&\min{\{
\eta^{\beta_{1}/2}_{I_{3}},\,
\lambda^{\beta_{1}/2}_{I_{3}\cup\{k_{*}\}}\}}\,\chi_{\Omega^{+}},\\
C(\cdot,\xi)&:=&|\log{\max{\{
\eta_{I_{3}},\,
\lambda_{I_{3}\cup\{k_{*}\}}
\}}}|\,\chi_{\Omega^{+}},\\
D_{m}(\cdot,\xi)&:=&
\min{\{\eta_{I_{3}}^{\beta_{m}/2},\,
\lambda^{\beta_{m}/2}_{I_{3}}\}}\,\chi_{\Omega_{I_{3}}}, \quad m=0,1.
\end{eqnarray*}
Let $\mathcal{X}$ be the integral operator
\begin{equation*}
\mathcal{X}u(\xi):=
\int u(\cdot,\xi)\prod_{n\in I_{3}}C_{n}(x_{n})dx_{n}
\end{equation*}
on functions locally integrable in $x_{n}$ with respect to the measure $C_{n}(x_{n}) dx_{n}$, for every $n\in I_{3}$,
for almost every $\xi$. Also, let $\mathcal{A}$, $\mathcal{B}$, $\mathcal{C}$,
$\mathcal{D}_{0}$ and $\mathcal{D}_{1}$, denote
$\mathcal{X}A$, $\mathcal{X}B$, $\mathcal{X}C$,
$\mathcal{X}D_{0}$ and $\mathcal{X}D_{1}$, respectively. By the discussion above, we may assume that
\begin{subequations}
\begin{equation}\label{misb1}
\mathcal{K}_{I_{1}, I_{2}, I_{3}}(\xi)\lesssim
\begin{cases}
\mathcal{A}(\xi)\tau(\xi)+\mathcal{B}(\xi)\tau_{k_{*}}(\xi),&
-2<\beta_{0}<-1,\\[0.05 cm]
\mathcal{A}(\xi)\tau(\xi)+\mathcal{C}(\xi)\tau_{k_{*}}(\xi),&
\beta_{0}=-1,\\[0.05 cm]
\mathcal{D}_{0}(\xi)\tau(\xi),&-1<\beta_{0}<0,
\end{cases}
\end{equation}
in the case \textbf{(b1)}, and
\begin{numcases}
{\mathcal{K}_{I_{1}, I_{2}, I_{3}}(\xi)\lesssim}
\label{misb2}\mathcal{B}(\xi)\tau_{k_{*}}(\xi),&
$-2<\beta_{0}<-1$,\hspace{0.5 cm} in the case \textbf{(b2)},\\[0.05 cm]
\label{misb3}\mathcal{D}_{1}(\xi)\tau(\xi),&
$-1<\beta_{1}<0$,\hspace{0.5 cm} in the case \textbf{(b3)}.
\end{numcases}
\end{subequations}
The last task is reduced then to estimating
\begin{equation}\label{etc}
\left.
\begin{aligned}
\mathcal{A}&,\quad -2<\beta_{0}\leq -1,\;\;\text{in the case}\; \textbf{(b1)},\\[0.05 cm]
\mathcal{B}&,\quad -1<\beta_{1}<0,\;\;\text{in the cases}\;  \textbf{(b1)}\;
\text{and}\; \textbf{(b2)},\\[0.05 cm]
\mathcal{C}&,\quad \beta_{0}=-1,
\;\;\text{in the case}\;\textbf{(b1)},\\[0.05 cm]
\mathcal{D}_{0}&,\quad -1<\beta_{0}<0,\;\;\text{in the case}\;  \textbf{(b1)},\\[0.05 cm]
\mathcal{D}_{1}&,\quad -1<\beta_{1}<0,\;\;\text{in the case}\;\textbf{(b3)},
\end{aligned}
\right\}
\end{equation}
in the region $(S_{r_{*}}(I_{*})\cap
\{\xi_{n_{*}}^{\overline{\ell}_{n_{*}}}
\approx \zeta_{d}\})\times V$. Obviously,
each of the functions in (\ref{etc}) depends on the partition $\{I_{1},I_{2},I_{3}\}$.
We continue to drop this
dependence from the notation for simplicity.
\smallskip \par
We record some key technical estimates in the next section.
We refer to Section \ref{adobh} for the notations.
\subsubsection{Some technical estimates}\label{str4}
%The proofs of Lemmas \ref{lmt1}, \ref{lmt2} and Corollary \ref{cor43} below are collected in Appendix~\autoref{app:pnda}.
\begin{lemma}\label{lmt1}
Fix $n\in I_{3}$, and a number $-2<\sigma<0$, and let $0<\lambda<1$.
If $I\subseteq I_{3}\setminus \{n\}$ and $I^{\prime}\subseteq I$ is
such that $\theta_{I^{\prime}}(\sigma)<0$,
then we have\smallskip\par
$\displaystyle
\int \min{\{\Psi_{I}^{\theta_{I^{\prime}}(\sigma)},\,
\lambda^{\theta_{I^{\prime}}(\sigma)}\}}
C_{n}(x_{n}) d x_{n}$
\begin{subequations}\label{tf}
\begin{numcases}
{\lesssim}
\xi_{n}^{2{\overline{\ell}}_{n}
(1-\alpha_{n})}\min{\{\Psi_{I\cup \{n\}}^{\theta_{I^{\prime}}(\sigma)},\,
\lambda^{\theta_{I^{\prime}}(\sigma)},\,\xi_{n}^{2{\overline{\ell}}_{n}
\theta_{I^{\prime}}(\sigma)}\}}
,&$\theta_{I^{\prime}\cup\{n\}}(\sigma)>0$;\label{tf1}
\\[0.2 cm]
|\log{\max{\{\Psi_{I\cup \{n\}},\,
\lambda\}}}|,&
$\theta_{I^{\prime}\cup\{n\}}(\sigma)=0$;\label{tf2}\\[0.2 cm]
\dfrac{\min{\{\xi_{n}^{2{\overline{\ell}}_{n}},\,
\max{\{\Psi_{I\cup \{n\}},\,
\lambda\}}\}}^{1-\alpha_{n}}}{\max{\{\Psi_{I\cup \{n\}},\,
\lambda\}}^{-
\theta_{I^{\prime}}(\sigma)}},
&$\theta_{I^{\prime}\cup\{n\}}(\sigma)<0$.
\label{tf4}
\end{numcases}
And if $\,\theta_{I^{\prime}\cup\{n\}}(\sigma)= 0$, then
\begin{multline}\label{tf5}
\int_{\max{\{ \Psi_{I\cup
\{n\} },\, \lambda \}}\gtrsim \xi_{n}^{2\overline{\ell}_{n}}} \min{\{\Psi_{I}^{\theta_{I^{\prime}}(\sigma)},\,
\lambda^{\theta_{I^{\prime}}(\sigma)}\}}
C_{n}(x_{n}) d x_{n}\\
\lesssim \xi_{n}^{2{\overline{\ell}}_{n}
(1-\alpha_{n})}\min{\{\Psi_{I\cup \{n\}}^{\theta_{I^{\prime}}(\sigma)},\,
\lambda^{\theta_{I^{\prime}}(\sigma)}\}}.
\end{multline}
Moreover, we have
\begin{equation}\label{i3t4}
\int
{|\log{\max{\{\Psi_{I},\,
\lambda\}}}|}
C_{n}(x_{n}) dx_{n}
\lesssim
\xi_{n}^{2{
\overline{\ell}}_{n}(1-\alpha_{n})}
|\log{\max{\{\Psi_{I\cup \{n\}},\,
\lambda,\,
\xi_{n}^{2{
\overline{\ell}}_{n}}\}}}|.
\end{equation}
\end{subequations}
\end{lemma}
\begin{proof}
Let $0<\lambda,\,\nu\ll 1$ and $\alpha\in (1/2,1)$,
and let $\eta$ be a real number. Then\medskip\\
\begin{subequations}
$\displaystyle \int_{\nu/4<t<4\nu}
|t-\nu|^{-\alpha}
\min{\{|t-\eta|^{-s},\,\lambda^{-s}\}} d t$
\begin{numcases}
{\lesssim}\label{stf1}\nu^{1-\alpha}
\min{\{\nu^{-s},\,
|\nu-\eta|^{-s},\,
\lambda^{-s}\}},&$0<s<1-\alpha$;
\\[0.2 cm]
\label{stf2}|\log{\max{\{
|\nu-\eta|,\,\lambda\}}}|,&\hspace{- 1 cm}
$s=1-\alpha$, $\nu \gg  \max{\{
|\nu-\eta|,\,\lambda\}}$;
\\[0.2 cm]
\label{stf3}\nu^{1-\alpha}\min{\{
|\nu-\eta|^{-s},\,\lambda^{-s}\}},&\hspace{- 1 cm}
$s=1-\alpha$, $\nu \lesssim \max{\{
|\nu-\eta|,\,\lambda\}}$;
\\[0.2 cm]
\label{stf4}
\dfrac{\min{\{\nu,\,\max{\{
|\nu-\eta|,\,\lambda\}}\}}^{1-\alpha}}
{\max{\{
|\nu-\eta|,\,\lambda\}}^{s}},&
$1-\alpha<s<1$.
\end{numcases}
\end{subequations}
The  inequalities (\ref{tf1}), (\ref{tf4}) and (\ref{tf5})
follow directly from (\ref{stf1}), (\ref{stf4})
and (\ref{stf3}), respectively, with $s=-\theta_{I^{\prime}}(\sigma)$, $\alpha=\alpha_{n}$,
$\nu=\xi_{n}^{2{
\overline{\ell}}_{n}}$, and
$\eta=1-\sum_{j\in I_{3}\setminus (I\cup \{n\})}{
x_{j}}-\sum_{k\in I}{
\xi^{2\overline{\ell}_{k}}_{k}}-\xi^{
2{\overline{\ell}}_{d}}_{d}$.
Similarly, the inequality (\ref{tf2}) is a
direct consequence of
(\ref{stf2}) and (\ref{stf3}) together,
and (\ref{i3t4}) follows from the estimate
(\ref{lgt3}).
\end{proof}
\smallskip\par
We give further technical
estimates in the next lemma. See Section \ref{adobh} for
the notations.
\begin{lemma}\label{lmt2}
Let $\sigma\in (-2,0)$ and $l_{0} \subsetneq I_{3}$, and let
$l=I_{3}\setminus l_{0}$. Define
\begin{equation*}
\mathcal{Z}_{l_{0}}(\xi;\sigma):=\int \min{\{
\Psi_{l_{0}}^{\sigma/2},\,
\lambda_{I_{*}}^{\sigma/2}(\xi)\}}\,
\prod_{n\in l} C_{n}(x_{n}) d x_{n}.
\end{equation*}
Suppose that $l=\{n_{1},\dots,n_{N}\}$ and that
$\xi_{n_{1}}^{{
\overline{\ell}}_{n_{1}}}\geq
\xi_{n_{2}}^{{
\overline{\ell}}_{n_{2}}}\geq\dots\geq
\xi_{n_{N}}^{{
\overline{\ell}}_{n_{N}}}$.
Then $\mathcal{Z}_{l_{0}}(\xi;\sigma)$ is majorized
on $\,S_{r_{*}}(I_{*})\times V$ by
\begin{subequations}\label{abc}
\begin{numcases}
{}\min{\{
\Phi_{I_{3}}^{
\theta_{\!K^{+}(I_{*})\cap l}(\sigma)}(\xi),\,
\lambda_{I_{*}}^{
\theta_{\!K^{+}(I_{*})\cap l}(\sigma)}(\xi)\}}
\widetilde{P}_{K^{-}(I_{*})\cap l}(\xi),&\nonumber\\
&\hspace{-2 cm}\textnormal{\bf(c1)}\textnormal{($l$,$\,I_{*}$,$\,\sigma$)};
\label{istm31f}\\[0.05 cm]
|\log{\max{\{\Phi_{I_{3}}(\xi),\,
\lambda_{I_{*}}(\xi),\,
\xi_{n_{N_{2}(\sigma)}}^{2{\overline{\ell}}_{n_{N_{2}(\sigma)}}}
\}}}|\widetilde{P}_{l\setminus l^{+}(\sigma)}(\xi),&
\hspace{-2 cm}\textnormal{\bf(c2)}\textnormal{($l$,$\,I_{*}$,$\,\sigma$)};
\label{istm32f}
\\[0.05 cm]
\min{\{
\Phi_{I_{3}}^{\theta_{l^{-}(\sigma)}(\sigma)}(\xi),\,
\lambda_{I_{*}}^{\theta_{l^{-}(\sigma)}(\sigma)}(\xi),
\,\xi_{n_{N_{1}(\sigma)}}^{2{\overline{\ell}}_{n_{N_{1}(\sigma)}}
\theta_{l^{-}(\sigma)}(\sigma)}\}}
\widetilde{P}_{l\setminus l^{-}(\sigma)}(\xi),&\nonumber\\
&\hspace{-2 cm}\textnormal{\bf(c3)}\textnormal{($l$,$\,I_{*}$,$\,\sigma$)},
\label{istm33f}
\end{numcases}
\end{subequations}
where
\begin{equation*}
\widetilde{P}_{I}(\xi):=\prod_{n\in I}
\xi_{n}^{2{\overline{\ell}}_{n}(1-\alpha_{n})},\quad
I\subseteq \Gamma_{d},
\end{equation*}
and \textnormal{\bf(c1)}\textnormal{($l$,$\,I_{*}$,$\,\sigma$)} - \textnormal{\bf(c3)}\textnormal{($l$,$\,I_{*}$,$\,\sigma$)} refer to the cases
\begin{description}
  \item[(c1)\textnormal{($l$,$\,I_{*}$,$\,\sigma$)}.]
$l^{-}(\sigma)\neq\emptyset$, $K^{+}(I_{*})\cap l\subseteq l^{-}(\sigma)$,
  \item[(c2)\textnormal{($l$,$\,I_{*}$,$\,\sigma$)}.] \textnormal{($l^{-}(\sigma)= \emptyset$, or
$\,\emptyset\neq l^{-}(\sigma)\subsetneq K^{+}(I_{*})\cap l$),
$\,\theta_{l^{+}(\sigma)}(\sigma)=0$},
  \item[(c3)\textnormal{($l$,$\,I_{*}$,$\,\sigma$)}.] \textnormal{($l^{-}(\sigma)= \emptyset$, or
$\,\emptyset\neq l^{-}(\sigma)\subsetneq K^{+}(I_{*})\cap l$),
$\,\theta_{l^{+}(\sigma)}(\sigma)>0$}.
\end{description}
\end{lemma}
\begin{proof}
By abuse of notation,
we identify throughout the proof of this lemma
each integer $n_{j}$ in $l$ with its subscript $j$.
In this notation, $l=\{1,2,\dots,N\}$, and we are in the region
\begin{equation*}
\xi_{1}^{{
\overline{\ell}}_{1}}\geq
\xi_{2}^{{
\overline{\ell}}_{2}}\geq\dots\geq
\xi_{N}^{{
\overline{\ell}}_{N}}.
\end{equation*}
For the sake of brevity, we shall also denote $K^{+}(I_{*})\cap l$ by
$K^{+}$, $K^{-}(I_{*})\cap l$ by
$K^{-}$, $\lambda_{I_{*}}(\xi)$ by $\lambda$, and
$\mathcal{Z}_{l_{0}}(\xi;\sigma)$ by
$\mathcal{Z}$, throughout the proof. The dependence of the sets $l^{-}$ and $l^{+}$, and the subindices $N^{-}$ and $N_{k}$, on $\sigma$
will be dropped as well from the notation. The proof of Lemma \ref{lmt2} relies mainly
on the estimates (\ref{tf}) of Lemma \ref{lmt1}.
It is convenient for application to
rewrite the estimates (\ref{tf1}) - (\ref{tf4})
in the following form:
\smallskip\par
\begin{subequations} \label{ttfs}
$\displaystyle
\int \min{\{\Psi_{I}^{\theta_{I^{\prime}}(\sigma)},\,
\lambda^{\theta_{I^{\prime}}(\sigma)}\}}
C_{n}(x_{n}) d x_{n}$
\vspace{-0.3 cm}
\begin{numcases}
{\lesssim}\xi_{n}^{2{\overline{\ell}}_{n}
\theta_{I^{\prime}\cup\{n\}}(\sigma)},&
\hspace{-0.25 cm}$\theta_{I^{\prime}\cup\{n\}}(\sigma)>0$;\label{ttf1}
\\[0.2 cm]
|\log{\max{\{\Psi_{I\cup \{n\}},\,
\lambda\}}}|,&
\hspace{-0.25 cm}$\theta_{I^{\prime}\cup\{n\}}(\sigma)=0$;\label{ttf2}
\\[0.2 cm]
\min{\{\Psi_{I\cup \{n\}}^{
\theta_{I^{\prime}\cup\{n\}}(\sigma)},\,
\lambda^{\theta_{I^{\prime}\cup\{n\}}(\sigma)}\}},
&\hspace{-0.25 cm}$\theta_{I^{\prime}\cup\{n\}}(\sigma)<0$;
\label{ttf3}
\\[0.2 cm]
\xi_{n}^{2{\overline{\ell}}_{n}
(1-\alpha_{n})}\min{\{\Psi_{I\cup \{n\}}^{\theta_{I^{\prime}}(\sigma)},\,
\lambda^{\theta_{I^{\prime}}(\sigma)}\}},&\hspace{-0.25 cm}
$\theta_{I^{\prime}\cup\{n\}}(\sigma)\neq 0$.\label{ttf4}
\end{numcases}
\end{subequations}
\smallskip\par Assume first that $l^{-}=\emptyset$. Then
$1-\alpha_{1}+\sigma/2\geq 0$. If $1-\alpha_{1}+\sigma/2=0$, we use the estimate (\ref{ttf2}),
with $I=l_{0}$ and $I^{\prime}=\emptyset$, for the
integral in $x_{1}$. Then, we apply (\ref{i3t4}) repeatedly, with
$(n,I)=(2,l_{0}\cup\{1\})$, $(3,l_{0}\cup\{1,2\})$, $\dots$,
$(N,l_{0}\cup\{1,2,\dots,N-1\})$ to estimate the integrals
in $x_{2},\dots,x_{N}$, respectively. We see then that, if
$1-\alpha_{1}+\sigma/2=0$, we have
\begin{equation}\label{j0mty1}
\mathcal{Z}\lesssim
|\log{\max{\{\Phi_{I_{3}},\,
\lambda,\,
\xi_{2}^{2{\overline{\ell}}_{2}},\dots,
\xi_{N}^{2{\overline{\ell}}_{N}}
\}}}|\prod_{2\leq n\leq N}
\xi_{n}^{2{\overline{\ell}}_{n}(1-\alpha_{n})}.
\end{equation}
This proves (\ref{istm32f}) when $l^{-}=\emptyset$ and
$\theta_{l^{+}}(\sigma)=0$.
If $1-\alpha_{1}+\sigma/2>0$, we can use
(\ref{ttf1}) once with $n=1$, $I=l_{0}$,
and $I^{\prime}=\emptyset$, to estimate the integral in $x_{1}$.
We find that
\begin{equation*}
\int \min{\{
\Psi_{l_{0}}^{\sigma/2},\,
\lambda^{\sigma/2}\}}\,
C_{1}(x_{1}) d x_{1}\lesssim
{\xi}_{1}^{2{\overline{\ell}}_{1}
(1-\alpha_{1}+\sigma/2)}.
\end{equation*}
Then, using the simple fact that
\begin{equation*}
\int \prod_{2\leq n\leq N} C_{n}(x_{n})dx_{n}\approx
\prod_{2\leq n\leq N}
\xi_{n}^{2{\overline{\ell}}_{n}
(1-\alpha_{n})},
\end{equation*}
we immediately get
\begin{align*}
\mathcal{Z}&\lesssim
{\xi}_{1}^{2{\overline{\ell}}_{1}
(1-\alpha_{1}+\sigma/2)}
\prod_{2\leq n\leq N}
\xi_{n}^{2{\overline{\ell}}_{n}(1-\alpha_{n})}\\
&={\xi}_{1}^{\overline{\ell}_{1}\sigma}
\prod_{1\leq n\leq N}
\xi_{n}^{2{\overline{\ell}}_{n}(1-\alpha_{n})},\quad
1-\alpha_{1}+\sigma/2>0.
\end{align*}
One could alternatively use (\ref{ttf4}) repeatedly, with
$(n,I)=(1,l_{0})$, $(2,l_{0}\cup\{1\})$, $\dots$,
$(N,l_{0}\cup\{1,\dots,N-1\})$, to obtain
that
\begin{equation*}
\mathcal{Z}\lesssim \min{\{
\Phi_{I_{3}}^{\sigma/2},\,\lambda^{\sigma/2}\}}
\prod_{1\leq n\leq N}
\xi_{n}^{2{\overline{\ell}}_{n}(1-\alpha_{n})},\quad
1-\alpha_{1}+\sigma/2>0.
\end{equation*}
From the latter two estimates, we deduce that
\begin{equation}\label{j0mty2}
\mathcal{Z}\lesssim
\min{\{
\Phi_{I_{3}}^{\sigma/2},\,
\lambda^{\sigma/2},\,
{\xi}_{1}^{\overline{\ell}_{1}\sigma}\}}
\prod_{1\leq n\leq N}
\xi_{n}^{2{\overline{\ell}}_{n}(1-\alpha_{n})},\quad
1-\alpha_{1}+\sigma/2>0,
\end{equation}
which proves (\ref{istm33f}) when $l^{-}=\emptyset$ and
$\theta_{l^{+}}(\sigma)>0$.\smallskip\par
Suppose now that $l^{-}\neq\emptyset$.
We have three cases to consider:
\begin{description}
  \item[(i)] $\xi_{1}^{2{
\overline{\ell}}_{1}}\geq
\dots\geq
\xi_{N}^{2{
\overline{\ell}}_{N}}\gtrsim\max{\{\Phi_{I_{3}},\,\lambda\}}$,
 \item[(ii)] $\max{\{\Phi_{I_{3}},\,\lambda\}}\gg \xi_{1}^{2{
\overline{\ell}}_{1}}\geq
\dots\geq
\xi_{N}^{2{
\overline{\ell}}_{N}}$,
  \item[(iii)]$\xi_{1}^{2{
\overline{\ell}}_{1}}\geq\dots\geq
\xi_{k^{+}}^{2{
\overline{\ell}}_{k^{+}}}
\gtrsim\max{\{\Phi_{I_{3}},\,\lambda\}}\gg   \xi_{k^{-}}^{2{
\overline{\ell}}_{k^{-}}} \geq
\dots\geq
\xi_{N}^{2{
\overline{\ell}}_{N}}$.
\end{description}
\medskip
\par \textbf{(i)}. In this case, we prove that
\begin{subequations}\label{i1lm2}
\begin{numcases}
{\mathcal{Z}\lesssim}
\min{\{\Phi_{I_{3}}^{
\theta_{l}(\sigma)},\,
\lambda^{
\theta_{l}(\sigma)}\}},&\hspace{-1.7 cm} $l^{-}=l$;
\label{istm1}\\[0.15cm]
|\log{\max{\{\Phi_{I_{3}},\,
\lambda,\,
\xi_{N_{2}}^{2{\overline{\ell}}_{N_{2}}},\dots,
\xi_{N}^{2{\overline{\ell}}_{N}}
\}}}|\prod_{N_{2}\leq n\leq N }
\xi_{n}^{2{\overline{\ell}}_{n}(1-\alpha_{n})},&\nonumber\\[0.05 cm]
&\hspace{-4.5 cm}
$\emptyset\neq l^{-}\subsetneq l$, $\theta_{l^{+}}(\sigma)=0$;
\label{istm2}
\\[0.15cm]
\min{\{
\Phi_{I_{3}}^{\theta_{l^{-}}(\sigma)},\,
\lambda^{\theta_{l^{-}}(\sigma)},
\,{\xi}_{N_{1}}^{2{\overline{\ell}}_{N_{1}}
\theta_{l^{-}}(\sigma)}\}}\prod_{N_{1}\leq n\leq N }
\xi_{n}^{2{\overline{\ell}}_{n}(1-\alpha_{n})},&
\nonumber\\[0.05 cm]
&\hspace{-4.5 cm}
$\emptyset\neq l^{-}\subsetneq l$,
$\theta_{l^{+}}(\sigma)>0$\label{istm3}.
\end{numcases}
\end{subequations}
Observe that (\ref{istm1}), (\ref{istm2}) and (\ref{istm3})
imply (\ref{istm31f}), (\ref{istm32f}) and
(\ref{istm33f}), respectively, when
$K^{+}=l$ ($K^{-}=\emptyset$) and $l^{-}\neq \emptyset$.\smallskip\par
We start with (\ref{istm1}). Let
$l^{-}=l$. Since $\theta_{I}(\sigma)<0$, for every
$I\subseteq l^{-}$, we can apply the estimate
(\ref{ttf3}) repeatedly for
$(n,I,I^{\prime})=(1,l_{0},\emptyset)$, $(2,l_{0}\cup\{1\},\{1\})$,
$\dots$, $(N^{-},l_{0}\cup \{1,2,\dots,N^{-}-1\},\{1,2,\dots,N^{-}-1\})$.
This gives
\begin{equation}\label{luql}
\int \min{\{
\Psi_{l_{0}}^{\sigma/2},\,
\lambda^{\sigma/2}\}}\,
\prod_{n\in l^{-}} C_{n}(x_{n}) d x_{n}\lesssim
\min{\{\Psi_{l_{0}\cup l^{-}}^{
\theta_{l^{-}}(\sigma)},\,
\lambda^{
\theta_{l^{-}}(\sigma)}\}},
\end{equation}
which yields (\ref{istm1}), since
$l_{0}\cup l^{-}=I_{3}$ when $l^{-}=l$.\smallskip\par
Assume that $l^{-}\neq \emptyset$ and $l^{-}\subsetneq l$.
Then either
$\theta_{l^{+}}(\sigma)=0$, or $\theta_{l^{+}}(\sigma)>0$.
The integral $\mathcal{Z}$ can be estimated
in these cases similarly to (\ref{j0mty1}) and
 (\ref{j0mty2}), respectively. We sketch the details: \smallskip\par
If $\theta_{l^{+}}(\sigma)=0$, we use (\ref{ttf2}) with
$n=N_{1}$, $I=l_{0}\cup l^{-}$, and $I^{\prime}=l^{-}$, to obtain
\begin{equation*}
\int \min{\{\Psi_{l_{0}\cup l^{-}}^{
\theta_{l^{-}}(\sigma)},\,
\lambda^{
\theta_{l^{-}}(\sigma)}\}}C_{N_{1}}(x_{N_{1}}) d x_{N_{1}}
\lesssim |\log{\max{\{\Psi_{l_{0}\cup l^{+}},\,
\lambda\}}}|.
\end{equation*}
Next, using (\ref{i3t4}) repeatedly,
with $n=N_{2},\dots,N$, we see that
\begin{multline*}
\int |\log{\max{\{\Psi_{l_{0}\cup l^{+}},\,
\lambda\}}}| \prod_{N_{2}\leq n\leq N}
C_{n}(x_{n}) d x_{n}\\
\lesssim
|\log{\max{\{\Phi_{I_{3}},\,
\lambda,\,
\xi_{N_{2}}^{2{\overline{\ell}}_{N_{2}}},\dots,
\xi_{N}^{2{\overline{\ell}}_{N}}
\}}}|\prod_{N_{2}\leq n\leq N}
\xi_{n}^{2{\overline{\ell}}_{n}(1-\alpha_{n})}.
\end{multline*}
The estimate (\ref{istm2}) follows from (\ref{luql})
and the last two inequalities combined. \smallskip \par
Let $\theta_{l^{+}}(\sigma)>0$. Using
(\ref{ttf1}) once with $n=N_{1}$, $I=l_{0}\cup l^{-}$,
and $I^{\prime}=l^{-}$ to estimate the integral
of (\ref{luql}) in $x_{N_{1}}$, we see that
\begin{equation*}
\int \min{\{
\Psi_{l_{0}}^{\sigma/2},\,
\lambda^{\sigma/2}\}}\,
\prod_{n\in l^{+}} C_{n}(x_{n}) d x_{n}\lesssim
{\xi}_{N_{1}}^{2{\overline{\ell}}_{N_{1}}
\theta_{\!l^{+}}(\sigma)}.
\end{equation*}
Integrating this inequality against
$\prod_{N_{2}\leq n\leq N } C_{n}(x_{n})$ gives
\begin{equation*}
\mathcal{Z}\lesssim
{\xi}_{N_{1}}^{2{\overline{\ell}}_{N_{1}}
\theta_{l^{+}}(\sigma)}\prod_{N_{2}\leq n\leq N}
\xi_{n}^{2{\overline{\ell}}_{n}(1-\alpha_{n})}
\,=\,
{\xi}_{N_{1}}^{2{\overline{\ell}}_{N_{1}}
\theta_{l^{-}}(\sigma)}\prod_{N_{1}\leq n\leq N}
\xi_{n}^{2{\overline{\ell}}_{n}(1-\alpha_{n})}.
\end{equation*}
On the other hand, if we apply
(\ref{ttf4}), with $(n,I)=(N_{1},l_{0}\cup l^{-})$,
$(N_{2},l_{0}\cup l^{+})$, $\dots$,
$(N,l_{0}\cup \{n_{1},n_{2},\dots,n_{N-1}\})$, it follows from (\ref{luql}) that
\begin{equation*}
\mathcal{Z}\lesssim
\min{\{
\Phi_{I_{3}}^{\theta_{l^{-}}(\sigma)},\,
\lambda^{\theta_{l^{-}}(\sigma)}\}}
\prod_{1\leq n\leq N}
\xi_{n}^{2{\overline{\ell}}_{n}(1-\alpha_{n})}.
\end{equation*}
The last two estimates together show (\ref{istm3}).\smallskip\par
\textbf{(ii)}
If  $\Phi_{I_{3}}\geq \lambda$, then we are in the region
\begin{equation*}
\Phi_{I_{3}}\gg \xi_{1}^{2{
\overline{\ell}}_{1}}\geq
\dots\geq
\xi_{N}^{2{
\overline{\ell}}_{N}}.
\end{equation*}
In this region, since $
\xi_{n}^{2\overline{\ell}_{n}}/{4}<x_{n}<
4\xi_{n}^{2\overline{\ell}_{n}}$, for all $n\in l$, we have
\begin{align*}
\Psi_{l_{0}}(x,\xi)&=\,
|1-\sum_{j\in l}{
x_{j}}-\sum_{n\in l_{0}}{
\xi^{2\overline{\ell}_{n}}_{n}}-\xi^{
2{\overline{\ell}}_{d}}_{d}|\\
&=\,
|1-\sum_{n\in I_{3}}{
\xi^{2\overline{\ell}_{n}}_{n}}-\xi^{
2{\overline{\ell}}_{d}}_{d}-\sum_{j\in l}{
(x_{j}-\xi^{2\overline{\ell}_{j}}_{j})}|\\
&\approx\
\Phi_{I_{3}}(\xi).
\end{align*}
Therefore, in the case \textbf{(ii)},
if $\Phi_{I_{3}}\geq \lambda$, then
\begin{eqnarray}
\nonumber \mathcal{Z}
& \approx& \min{\{
\Phi_{I_{3}}^{\sigma/2},\,
\lambda^{\sigma/2}\}}\,
\int
\prod_{n\in l} C_{n}(x_{n}) d x_{n}\\
\label{ldaxi}& \approx &\min{\{
\Phi_{I_{3}}^{\sigma/2},\,
\lambda^{\sigma/2}\}} \prod_{1\leq n\leq N}
\xi_{n}^{2{\overline{\ell}}_{n}(1-\alpha_{n})}.
\end{eqnarray}
If  $\Phi_{I_{3}}\leq \lambda$, then we are in the region
\begin{equation}\label{lmdxi1}
\lambda\gg \xi_{1}^{2{
\overline{\ell}}_{1}}\geq
\dots\geq
\xi_{N}^{2{
\overline{\ell}}_{N}}.
\end{equation}
Assume that $1-\alpha_{n}+\sigma/2\neq 0$
for every $n=1,\dots,N$. In that case, we can employ the estimate
(\ref{ttf4}) with $I^{\prime}=\emptyset$ and $(n,I)=(1,l_{0})$, $(2,l_{0}\cup\{1\})$, $\dots$,
$(N,l_{0}\cup \{1,\dots,N-1\})$. If $1-\alpha_{n}+\sigma/2=0$,
for some $n\in l$, we are allowed, by the
restriction (\ref{lmdxi1}), to apply the estimate
(\ref{tf5}) in that step.
In fact, since we are under the assumption that
 $l^{-}\neq \emptyset$, we know that $1-\alpha_{1}+\sigma/2<0$.
So, by the estimate (\ref{ttf4}), we have
\begin{equation*}
\int \min{\{\Psi_{l_{0}}^{\sigma/2},
\lambda^{\sigma/2}\}}
C_{1}(x_{1}) d x_{1}\lesssim
{\xi}_{1}^{2{\overline{\ell}}_{1}
(1-\alpha_{1})}\min{\{\Psi_{l_{0}\cup \{1\}}^{
\sigma/2},
\lambda^{\sigma/2}\}}.
\end{equation*}
Then, for each
$n=2,\dots, N$, using
 (\ref{ttf4}), if $1-\alpha_{n}+\sigma/2\neq 0$,
 and (\ref{tf5}), if $1-\alpha_{n}+\sigma/2=0$, we obtain
the pointwise bound
\begin{multline*}
\int \min{\{\Psi_{l_{0}\cup\{1,\dots,n-1\}}^{\sigma/2},
\lambda^{\sigma/2}\}}
C_{n}(x_{n}) d x_{n} \lesssim
\xi_{n}^{2{\overline{\ell}}_{n}
(1-\alpha_{n})}\min{\{\Psi_{l_{0}\cup \{1,\dots,n\}}^{
\sigma/2},
\lambda^{\sigma/2}\}},
\end{multline*}
in the region (\ref{lmdxi1}).
Iterating this, we eventually obtain
\begin{equation}\label{iistm}
\mathcal{Z}\lesssim
\min{\{
\Phi_{I_{3}}^{\sigma/2},\,
\lambda^{\sigma/2}\}}
\prod_{1\leq n\leq N}
\xi_{n}^{2{\overline{\ell}}_{n}(1-\alpha_{n})}.
\end{equation}
As we have already seen in (\ref{ldaxi}),
the estimate (\ref{iistm}) holds in the region
$\Phi_{I_{3}}\geq \lambda$ as well.
The estimate (\ref{iistm}) proves (\ref{istm31f}) when
$K^{+}=\emptyset$, i.e., when $K^{-}=l$.
\smallskip\par
\textbf{(iii)}
By definition,
$K^{+}=\{1,\dots,k^{+}\}$ and  $K^{-}=\{k^{-},\dots,N\}$.
Denote $l_{0}\cup K^{+}$ by $K_{0}^{+}$.
By a straightforward adaptation of the proof of
the estimates (\ref{i1lm2}) obtained for the case \textbf{(i)}, we infer that
\medskip\\
\indent $\displaystyle
\int \min{\{
\Psi_{l_{0}}^{\sigma/2},\,
\lambda^{\sigma/2}\}}\,
\prod_{n\in K^{+}} C_{n}(x_{n}) d x_{n}$
\vspace{-0.2 cm}
\begin{subequations}
\begin{numcases}
{\lesssim}
\min{\{\Psi_{K_{0}^{+}}^{
\theta_{\!K^{+}}(\sigma)},\,
\lambda^{
\theta_{\!K^{+}}(\sigma)}\}},&\hspace{-1.5 cm}
$l^{-}\supseteq K^{+}$;
\label{istm31}\\[0.15cm]
|\log{\max{\{\Psi_{K_{0}^{+}},\,
\lambda,\,
\xi_{N_{2}}^{2{\overline{\ell}}_{N_{2}}},\dots,
\xi_{k^{+}}^{2{\overline{\ell}}_{k^{+}}}
\}}}|\prod_{N_{2}\leq n\leq k^{+}}
\xi_{n}^{2{\overline{\ell}}_{n}(1-\alpha_{n})},&\nonumber\\[0.05 cm]
&\hspace{-4.3 cm}
$\emptyset\neq l^{-}\subsetneq K^{+}$, $\theta_{l^{+}}(\sigma)=0$;
\label{istm32}
\\[0.15cm]
\min{\{
\Psi_{K_{0}^{+}}^{\theta_{l^{-}}(\sigma)},\,
\lambda^{\theta_{l^{-}}(\sigma)},
\,{\xi}_{N_{1}}^{2{\overline{\ell}}_{N_{1}}
\theta_{l^{-}}(\sigma)}\}}\prod_{N_{1}\leq n\leq k^{+}}
\xi_{n}^{2{\overline{\ell}}_{n}(1-\alpha_{n})},&
\nonumber\\[0.05 cm]
&\hspace{-4.3 cm}
$\emptyset\neq l^{-}\subsetneq K^{+}$,
$\theta_{l^{+}}(\sigma)>0$\label{istm33}.
\end{numcases}
\end{subequations}
By definition of $l^{-}$, we have
$\theta_{K^{+}}(\sigma)<0$ when
$K^{+}\subseteq l^{-}$.
With $\sigma$ replaced by $2\theta_{K^{+}}(\sigma)$,
$l_{0}$ by $K_{0}^{+}$, and $l$ by $K^{-}$, a
straightforward adaptation of the proof of (\ref{iistm})
yields
\begin{subequations}
\begin{multline}\label{istm31c}
\int \min{\{\Psi_{K_{0}^{+}}^{
\theta_{\!K^{+}}(\sigma)},\,
\lambda^{
\theta_{\!K^{+}}(\sigma)}\}}\,
\prod_{n\in K^{-}} C_{n}(x_{n}) d x_{n}\\
\lesssim
\min{\{
\Phi_{I_{3}}^{
\theta_{\!K^{+}}(\sigma)},\,
\lambda^{
\theta_{\!K^{+}}(\sigma)}\}}
\prod_{n\in K^{-}}
\xi_{n}^{2{\overline{\ell}}_{n}(1-\alpha_{n})}.
\end{multline}
We similarly have
\begin{multline}\label{istm33c}
\int \min{\{
\Psi_{K_{0}^{+}}^{\theta_{l^{-}}(\sigma)},\,
\lambda^{\theta_{l^{-}}(\sigma)},
\,{\xi}_{N_{1}}^{2{\overline{\ell}}_{N_{1}}
\theta_{l^{-}}(\sigma)}\}}\,
\prod_{n\in K^{-}} C_{n}(x_{n}) d x_{n}\\
\lesssim
\min{\{
\Phi_{I_{3}}^{\theta_{l^{-}}(\sigma)},\,
\lambda^{\theta_{l^{-}}(\sigma)},
\,{\xi}_{N_{1}}^{2{\overline{\ell}}_{N_{1}}
\theta_{l^{-}}(\sigma)}\}}
\prod_{n\in K^{-}}
\xi_{n}^{2{\overline{\ell}}_{n}(1-\alpha_{n})}.
\end{multline}
On the other hand, applying (\ref{i3t4}) of Lemma \ref{lmt1}
repeatedly with $(n,I)=(k^{-},K_{0}^{+})$,
$(k^{-}+1,K_{0}^{+}\cup\{k^{-}\})$,
$\dots$,
$(N,K_{0}^{+}\cup (K^{-}\setminus \{N\}))$,
we get
\begin{multline}\label{istm32c}
\int |\log{\max{\{\Psi_{K_{0}^{+}},\,
\lambda,\,
\xi_{N_{2}}^{2{\overline{\ell}}_{N_{2}}},\dots,
\xi_{k^{+}}^{2{\overline{\ell}}_{k^{+}}}
\}}}|\,
\prod_{n\in K^{-}} C_{n}(x_{n}) d x_{n}\\
\lesssim
|\log{\max{\{\Phi_{I_{3}},\,
\lambda,\,
\xi_{N_{2}}^{2{\overline{\ell}}_{N_{2}}},\dots,
\xi_{N}^{2{\overline{\ell}}_{N}}\}}}|\,
\prod_{n\in K^{-}}
\xi_{n}^{2{\overline{\ell}}_{n}(1-\alpha_{n})}.
\end{multline}
\end{subequations}
Combining (\ref{istm31}) with (\ref{istm31c}),
(\ref{istm32}) with (\ref{istm32c}),
and
(\ref{istm33}) with (\ref{istm33c}), we obtain
the estimates
(\ref{istm31f}), (\ref{istm32f}), and (\ref{istm33f}),
respectively, when $K^{+}$, $K^{-}$, and $l^{-}$
are all nonempty. This concludes the proof of Lemma \ref{lmt2}.
\end{proof}
\smallskip\par
In the sequel, let
 $I_{3}=\{n_{*},n_{1},\dots,n_{N}\}$ and suppose that
\begin{equation}\label{rgnn}
\xi_{n_{*}}^{{
\overline{\ell}}_{n_{*}}}\geq
\xi_{n_{1}}^{{
\overline{\ell}}_{n_{1}}}\geq
\dots\geq
\xi_{n_{N}}^{{
\overline{\ell}}_{n_{N}}}.
\end{equation}
In this set up, let
\begin{equation*}
 O_{0}=\{n_{*}\},
\quad O_{j}:=\{n_{*},n_{1},
\dots,n_{j}\},\quad 1\leq j\leq N,
\end{equation*}
and for each $0\leq j<N$. Let $\beta\in
(-2,0)$ and assume that $1-\alpha_{n_{*}}+\beta/2<0$.
In the notations introduced in Section \ref{adobh},
if $\theta_{O_{1}}(\beta)\geq 0$, then
$I^{-}_{3}(\beta)=O_{0}$,
$I^{+}_{3}(\beta)= O_{1}$, $n_{N^{-}(\beta)}=n_{*}$, and
$n_{N_{k}(\beta)}=n_{k}$, $1\leq k\leq N$.
If $\theta_{I_{3}}(\beta)<0$, then $I^{-}_{3}(\beta)=
I_{3}$ and $I^{+}_{3}(\beta)=\emptyset$.
While if
 $\theta_{O_{1}}(\beta)< 0$ and $\theta_{I_{3}}(\beta)\geq 0$, then
$I^{-}_{3}(\beta)=O_{N^{-}(\beta)}$,
where $O_{1}\subseteq O_{N^{-}(\beta)}\subsetneq I_{3}$ is such that
$\theta_{O_{N^{-}(\beta)}}(\beta)<0$ and
$\theta_{O_{N_{1}(\beta)}}(\beta)\geq 0$,
$N_{1}(\beta)=N^{-}(\beta)+1$.
In this case, $I^{+}_{3}(\beta)= O_{N_{1}(\beta)}$ and we are in the region
\begin{equation*}
\xi_{n_{*}}^{{
\overline{\ell}}_{n_{*}}}\geq
\xi_{n_{1}}^{{
\overline{\ell}}_{n_{1}}}\geq
\dots\geq
\xi_{n_{N^{-}(\beta)}}^{{
\overline{\ell}}_{n_{N^{-}(\beta)}}}\geq
\xi_{n_{N_{1}(\beta)}}^{{
\overline{\ell}}_{n_{N_{1}(\beta)}}}\geq
\xi_{n_{N_{2}(\beta)}}^{{
\overline{\ell}}_{n_{N_{2}(\beta)}}}\geq\dots
\xi_{n_{N}}^{{
\overline{\ell}}_{n_{N}}},
\end{equation*}
where $N_{k}(\beta)=N^{-}(\beta)+k$,
$k=1,\dots,N-N^{-}(\beta)$.  Recall that $I_{*}$ denotes either
$I_{3}$ or $I_{3}\cup\{k_{*}\}$. Introduce the function
\begin{equation*}
\Lambda_{I_{3}}(\xi;\beta,I_{*}):=
\begin{cases}
\max{\{
\Phi_{I_{3}}(\xi),
\lambda_{I_{*}}(\xi)\}},&\hspace{0.5 cm}  \text{\bf(c1)}\textnormal{($I_{3}$,$\,I_{*}$,$\,\beta$)};
\\[0.2 cm]
\max{\{\Phi_{I_{3}}(\xi),\,
\lambda_{I_{*}}(\xi),\,
\xi_{n_{N_{2}(\beta)}}^{2{\overline{\ell}}_{n_{N_{2}(\beta)}}}
\}},&\hspace{0.5 cm} \text{\bf(c2)}\textnormal{($I_{3}$,$\,I_{*}$,$\,\beta$)};
\\[0.2 cm]
\max{\{
\Phi_{I_{3}}(\xi),\,
\lambda_{I_{*}}(\xi),
\,{\xi}_{n_{N_{1}(\beta)}}^{2{\overline{\ell}}_{n_{N_{1}(\beta)}}
}\}},&\hspace{0.5 cm}  \text{\bf(c3)}\textnormal{($I_{3}$,$\,I_{*}$,$\,\beta$)},
\end{cases}
\end{equation*}
where
\text{\bf(c1)}\textnormal{($I_{3}$,$\,I_{*}$,$\,\beta$)} -
\text{\bf(c3)}\textnormal{($I_{3}$,$\,I_{*}$,$\,\beta$)} denote
the cases
\begin{description}
\item[\text{\bf(c1)}\textnormal{($I_{3}$,$\,I_{*}$,$\,\beta$)}.]
\hspace{0.1 cm} $I^{-}_{3}(\beta)\neq\emptyset$, $K^{+}(I_{*})\cap I_{3}\subseteq I^{-}_{3}(\beta)$,
\item[\text{\bf(c2)}\textnormal{($I_{3}$,$\,I_{*}$,$\,\beta$)}.]
\hspace{0.1 cm} \textnormal{($I^{-}_{3}(\beta)= \emptyset$, or
$\,\emptyset\neq I^{-}_{3}(\beta)\subsetneq K^{+}(I_{*})\cap I_{3}$), $\,\theta_{I^{+}_{3}(\beta)}(\beta)=0$},
\item[\text{\bf(c3)}\textnormal{($I_{3}$,$\,I_{*}$,$\,\beta$)}.]
\hspace{0.1 cm}\textnormal{($I^{-}_{3}(\beta)= \emptyset$, or
$\,\emptyset\neq I^{-}_{3}(\beta)
\subsetneq K^{+}(I_{*})\cap I_{3}$), $\,
\theta_{I^{+}_{3}(\beta)}(\beta)>0$}.
\end{description}
\begin{corollary}\label{cor43}
Let $-2<\beta<0$, and suppose that $ O_{j}\subseteq
K^{+}(I_{*})\cap I^{-}_{3}(\beta)$,
for some $0\leq j<N$. Then, we have
\begin{equation}\label{uxi}
{\int
\min{\{\Psi^{\theta_{O_{j}}(\beta)}_{O_{j}},\,
\lambda^{\theta_{O_{j}}(\beta)}_{I_{*}}(\xi)
\}}
\prod_{n\in I_{3}\setminus O_{j}}
C_{n}(x_{n})dx_{n}}
\lesssim \mathcal{O}_{I_{3}}(\xi;\beta,I_{*}),
\end{equation}
on the region $S_{r_{*}}(I_{*})\times V$,
where
\begin{equation*}
\mathcal{O}_{I_{3}}(\xi;\beta,I_{*}):=
\begin{cases}
\Lambda_{I_{3}}^{\theta_{\!K^{+}(I_{*})\cap I_{3}}(\beta)}(\xi;\beta,I_{*})\,
\widetilde{P}_{K^{-}(I_{*})\cap I_{3}}(\xi),& \text{\bf(c1)}\textnormal{($I_{3}$,$\,I_{*}$,$\,\beta$)};
\\[0.2 cm]
|\log{\Lambda_{I_{3}}(\xi;\beta,I_{*})}|\,
\widetilde{P}_{I_{3}\setminus I^{+}_{3}(\beta)}(\xi),& \text{\bf(c2)}\textnormal{($I_{3}$,$\,I_{*}$,$\,\beta$)};
\\[0.2 cm]
\Lambda_{I_{3}}^{\theta_{I^{-}_{3}(\beta)}(\beta)}
(\xi;\beta,I_{*})\,\widetilde{P}_{I_{3}\setminus I^{-}_{3}(\beta)}(\xi),&  \text{\bf(c3)}\textnormal{($I_{3}$,$\,I_{*}$,$\,\beta$)}.
\end{cases}
\end{equation*}
\end{corollary}
\begin{proof}
\par With $l_{0}= O_{j}$, $l=I_{3}\setminus O_{j}$, and $\sigma=2\theta_{O_{j}}(\beta)$, the estimate (\ref{uxi}) follows
from the estimates (\ref{abc}) of Lemma \ref{lmt2}
applied to the integral $\mathcal{Z}_{O_{j}}(\xi;\sigma)$.
Observe that
$O_{j}\cup (K^{+}(I_{*})\cap l)=
K^{+}(I_{*})\cap I_{3}$, and $K^{-}(I_{*})\cap l=
K^{-}(I_{*})\cap I_{3}$. Consequently, $
\theta_{K^{+}(I_{*})\cap l}
(\sigma)
=\theta_{K^{+}(I_{*})\cap I_{3}}(\beta)$.
And since $\theta_{l^{-}(\sigma)}
(\sigma)=
\theta_{O_{j}\cup l^{-}(\sigma)}(\beta)$, we have $\theta_{l^{-}(\sigma)}
(\sigma)=\theta_{I^{-}_{3}(\beta)}(\beta)$,
and we also have $O_{j}\cup l^{-}(\sigma)=
I^{-}_{3}(\beta)$, $l\setminus l^{-}(\sigma)=I_{3}\setminus
I^{-}_{3}(\beta)$ and $l\setminus l^{+}(\sigma)=
I_{3}\setminus I^{+}_{3}(\beta)$.
\end{proof}
\smallskip\par
We are ready now to estimate the functions (\ref{etc}).
\subsubsection{Estimates for \texorpdfstring{{\mathversion{bold}$
\mathcal{K}_{I_{1}, I_{2}, I_{3}}$}}{} in the cases \textbf{(b1)} - \textbf{(b3)}}\label{lstp}
We continue to assume that $
I_{3}=\{n_{*},n_{1},n_{2},\dots,n_{N}\}$ and that we are in the region (\ref{rgnn}).
Having restricted to $\zeta_{d}\approx \xi_{n_{*}}^{{
\overline{\ell}}_{n_{*}}}$, we have
\begin{equation*}
\Phi_{I_{3}}(\xi)\,\leq\, \sum_{n\in I_{3}} \xi_{n}^{2{
\overline{\ell}}_{n}}+ |1-\xi_{d}^{2\overline{\ell}_{d}}|\,
\approx\, \xi_{n_{*}}^{2{
\overline{\ell}}_{n_{*}}}+ \zeta^{2}_{d}\,\approx\,
\xi_{n_{*}}^{2{
\overline{\ell}}_{n_{*}}},
\end{equation*}
for $\xi_{d}\in V$, in all the cases \textbf{(b1)} - \textbf{(b3)}. Thus, by (\ref{b123}),
\begin{equation*}
\begin{cases}
\xi_{n_{*}}^{2{
\overline{\ell}}_{n_{*}}}\gtrsim
\max{\{\Phi_{I_{3}}(\xi),\,\lambda_{I_{3}}(\xi)\}},&
\textbf{(b1)},\;\textbf{(b3)};\\[0.05 cm]
\xi_{k_{*}}^{2{
\overline{\ell}}_{k_{*}}}\gg
\xi_{n_{*}}^{2{
\overline{\ell}}_{n_{*}}}\gtrsim
\max{\{\Phi_{I_{3}}(\xi),\,
\lambda_{I_{3}\cup\{k_{*}\}}\}},& \textbf{(b2)}.
\end{cases}
\end{equation*}
In other words
\begin{equation*}
\begin{cases}
n_{*}\in K^{+}(I_{3}),&
 \textbf{(b1)},\;\textbf{(b3)};\\[0.05 cm]
n_{*}\in K^{+}(I_{3}\cup\{k_{*}\}),&
 \textbf{(b2)}.
\end{cases}
\end{equation*}
\par
Let us begin with $\mathcal{A}$,
when $-2<\beta_{0}\leq -1$, in the case
\textbf{(b1)}. As discussed above, the case \textbf{(b1)} puts us
in the region where  $ \xi_{n_{*}}^{2{\overline{\ell}}_{n_{*}}}\gtrsim $
$\max{\left\{\Phi_{I_{3}}(\xi),\right.}$ ${\,\left.\lambda_{I_{3}}(\xi)\right\}}$, which means
that
$n_{*}\in K^{+}(I_{3})$.
If $-2<\beta_{0}\leq-1$, one has
\begin{equation}\label{balf1}
1-\alpha_{n_{*}}<-\beta_{0}/2<1,
\end{equation}
and consequently
$\theta_{\{n_{*}\}}(\beta_{0})=1-\alpha_{n_{*}}+\beta_{0}/2<0$.
This shows that
\begin{equation}\label{nstr}
n_{*}\in K^{+}(I_{3})\cap I^{-}_{3}(\beta_{0}).
\end{equation}
Let $0<\mu\ll 1$,
and let $\alpha\in (1/2,1)$ and $1-\alpha<s<1$.
A more subtle version of the estimate (\ref{stf4}) reads
\begin{multline}\label{tvum}
\int_{\substack{\nu/4<t<4\nu\\
\max{\{|t-\eta|,\,\lambda\}}\gtrsim \mu}}
\dfrac{\min{\{|t-\eta|^{-s},\,\lambda^{-s}\}}}
{|t-\nu|^{\alpha}} d t
\lesssim
\dfrac{\min{\{\nu,\,\max{\{
|\nu-\eta|,\,\lambda,\,\mu\}}\}}^{1-\alpha}}
{\max{\{
|\nu-\eta|,\,\lambda,\,\mu\}}^{s}}\\
\leq\min{\{
|\nu-\eta|^{1-\alpha-s},\,\lambda^{1-\alpha-s},\,\mu^{1-\alpha-s}\}}.
\end{multline}
By (\ref{balf1}), the estimate (\ref{tvum}) can be applied
with $s=-\beta_{0}/2$, $\alpha=\alpha_{n_{*}}$,
$\lambda=\xi_{j_{*}}^{2{
\overline{\ell}}_{j_{*}}}$, $\mu=\xi_{k_{*}}^{2{
\overline{\ell}}_{k_{*}}}$,
$\nu=\xi_{n_{*}}^{2{
\overline{\ell}}_{n_{*}}}$, and
$\eta=1-\xi^{
2{\overline{\ell}}_{d}}_{d}-
\sum_{n\in I_{3}\setminus \{n_{*}\} } x_{n}$, to obtain
\begin{equation*}
\int A(\cdot,\xi)
C_{n_{*}}(x_{n_{*}})dx_{n_{*}} \lesssim
\min{\{\Psi_{\{n_{*}\}}^{1-\alpha_{n_{*}}+\beta_{0}/2},\,
\lambda_{I_{3}}^{1-\alpha_{n_{*}}+\beta_{0}/2}(\xi)\}}.
\end{equation*}
Then, by (\ref{nstr}), it follows from Corollary \ref{cor43}, with
$j=0$, $\beta=\beta_{0}$ and $I_{*}=I_{3}$, that
\begin{equation}\label{ma}
\mathcal{A}(\xi)\lesssim
\mathcal{O}_{I_{3}}(\xi;\beta_{0},I_{3}).
\end{equation}
\par Next, Lemma \ref{lmt2} with $\sigma=\beta_{m}$,
$l_{0}=\emptyset$, $l=I_{3}$ and $I_{*}=I_{3}$ shows that
\begin{equation}\label{mam12}
\mathcal{D}_{m}(\xi)\lesssim
\mathcal{O}_{I_{3}}(\xi;\beta_{m},I_{3}),\quad m=0,1.
\end{equation}
\par
Now, we estimate
$\mathcal{B}$, $-2<\beta_{0}<-1$,
in the case \textbf{(b1)}.
Both (\ref{balf1}) and (\ref{nstr}) continue to hold. Let $1/2<\alpha<1$, $0<\lambda,\mu,\nu\ll 1$,
and $|\eta|\lesssim \nu$. Then
\begin{multline}\label{ii3m}
\int_{\substack{\nu/4<t<4\nu,\\
\max{\{|t-\eta|,\lambda\}}\lesssim\mu\lesssim \nu}}
|t-\nu|^{-\alpha}
\min{\{|t-\eta|^{-s},\,\lambda^{-s}\}} d t\\[0.05 cm]
\lesssim\begin{cases}
{\mu}^{1-s}\min{\{{\mu}^{-\alpha},
\,|\nu-\eta|^{-\alpha}\}},&
s<1-\alpha;
\\[0.05 cm]
|\log{\max{\{
|\nu-\eta|,\,\lambda\}}}|,&
s=1-\alpha,\; \mu\gg  |\nu-\eta|;
\\[0.05 cm]
{\mu}^{1-s}\min{\{{\mu}^{-\alpha},
\,|\nu-\eta|^{-\alpha}\}},&
s=1-\alpha,\;\mu \lesssim |\nu-\eta|;
\\[0.05 cm]
\dfrac{\min{\{\mu,\,\max{\{
|\nu-\eta|,\,\lambda\}}\}}^{1-s}}
{\max{\{
|\nu-\eta|,\,\lambda\}}^{\alpha}},&
1-\alpha<s<1/2,
\end{cases}
\end{multline}
and we have
\begin{equation}\label{i3t4m}
\int_{\substack{\nu/4<t<4\nu,\\
\max{\{|t-\eta|,\lambda\}}\lesssim\mu\lesssim \nu}}
\dfrac{|\log{ \max{\{|t-\eta|,\lambda\}}}|}{|t-\nu|^{\alpha}} dt
\lesssim
\dfrac{\mu|\log{\mu}|}{\max{\{
|\nu-\eta|,\,\lambda,\,\mu\}^{\alpha}}}.
\end{equation}
Let $n\in I_{3}$ and $I\subseteq I_{3}\setminus\{n\}$.
We know that $\theta_{I}(\beta_{1})\geq \beta_{1}/2 >-1/2$. Applying the estimates (\ref{ii3m})
with  $s=-\theta_{I}(\beta_{1})$, $\alpha=\alpha_{n}$,
$\lambda=\lambda_{I_{3}\cup \{k_{*}\}}(\xi)$, $\mu=\xi_{k_{*}}^{2{
\overline{\ell}}_{k_{*}}}$,
$\nu=\xi_{n}^{2{
\overline{\ell}}_{n}}$, and
$\eta=1-\sum_{j\in I_{3}\setminus (I\cup\{n\})} x_{j}
-\sum_{j\in I} \xi^{2\overline{\ell}_{j}}_{j}
-\xi^{2{\overline{\ell}}_{d}}_{d}$,
we find that \medskip\\
$\displaystyle \int \min{\{
\Psi^{\theta_{I}(\beta_{1})}_{I},\,
\lambda^{\theta_{I}(\beta_{1})}_{I_{3}\cup \{k_{*}\}}
\}}
\,\chi_{\{x_{n}\gtrsim\xi_{k_{*}}^{2{
\overline{\ell}}_{k_{*}}}\gtrsim
\max{\{\Psi_{I},\,\lambda_{I_{3}\cup\{k_{*}\}}\}}\}}\,
C_{n}(x_{n})dx_{n}$\\ \vspace{-0.15 cm}
\begin{subequations}\label{glms}
\begin{numcases}
{\lesssim}
\xi_{k_{*}}^{2{\overline{\ell}}_{k_{*}}
(1+\theta_{I}(\beta_{1}))}
\min{\{
\,\Psi_{I\cup\{n\}}^{-\alpha_{n}},\,
\lambda_{I_{3}}^{-\alpha_{n}}\}},&
$\theta_{I\cup \{n\}}(\beta_{1})>0$;\label{glm1}
\\[0.05 cm]
|\log{\max{\{
\Psi_{I\cup\{n\}},\,\lambda_{I_{3}\cup \{k_{*}\}}\}}}|,&
\hspace{-2.5 cm} $\theta_{I\cup \{n\}}(\beta_{1})=0$,\;$\xi_{k_{*}}^{2{
\overline{\ell}}_{k_{*}}}\gg  \Psi_{I\cup\{n\}}$;\label{glm2}
\\[0.05 cm]
\xi_{k_{*}}^{2{
\overline{\ell}}_{k_{*}}(1+\theta_{I}(\beta_{1}))}
\min{\{
\,\Psi_{I\cup\{n\}}^{-\alpha_{n}},\,
\lambda_{I_{3}}^{-\alpha_{n}}\}},&\hspace{-2 cm}
$\theta_{I\cup \{n\}}(\beta_{1})=0$,\;\;$\xi_{k_{*}}^{2{
\overline{\ell}}_{k_{*}}}\lesssim \Psi_{I\cup\{n\}}$;\label{glm3}
\\[0.05 cm]
\dfrac{\min{\{\xi_{k_{*}}^{2{
\overline{\ell}}_{k_{*}}},\,\max{\{
\Psi_{I\cup\{n\}},\,\lambda_{I_{3}\cup \{k_{*}\}}\}}\}^{
1+\theta_{I}(\beta_{1})}}}
{\max{\{
\Psi_{I\cup\{n\}},\,\lambda_{I_{3}\cup \{k_{*}\}}\}}^{
\alpha_{n}}},&
$\theta_{I\cup \{n\}}(\beta_{1})<0$.
\label{glm4}
\end{numcases}
\end{subequations}
The estimates (\ref{glms}) in this general
form will be useful later.
It also follows from (\ref{i3t4m}) that
\begin{multline}\label{plmlg}
 \int |\log{\max{\{
\Psi_{I},\,\lambda_{I_{3}\cup \{k_{*}\}}\}}}|
\,\chi_{\{x_{n}\gtrsim\xi_{k_{*}}^{2{
\overline{\ell}}_{k_{*}}}\gtrsim\max{\{\Psi_{I},\,\lambda_{I_{3}\cup\{k_{*}\}}\}}\}}\,
C_{n}(x_{n})dx_{n}\\
\lesssim
\xi_{k_{*}}^{2{
\overline{\ell}}_{k_{*}}}|\log{\xi_{k_{*}}}|
\min{\{
\Psi_{I\cup\{n\}}^{-\alpha_{n}},
\,\lambda_{I_{3}}^{-\alpha_{n}}\}}.
\end{multline}
For $x\in \Omega^{+}$, one has
\begin{equation*}
\xi_{n_{*}}^{2{\overline{\ell}}_{n_{*}}}\gg \max{\{
\eta_{I_{3}},\,\lambda_{I_{3}\cup \{k_{*}\}}\}},\;\;
\xi_{k_{*}}^{2{\overline{\ell}}_{k_{*}}}\gg \max{\{
\eta_{I_{3}},\,\xi_{j_{*}}^{2{\overline{\ell}}_{j_{*}}}\}}.
\end{equation*}
So, by (\ref{b123}) in the case \textbf{(b1)}, one has
$\xi_{n_{*}}^{2{
\overline{\ell}}_{n_{*}}}\gtrsim\xi_{k_{*}}^{2{
\overline{\ell}}_{k_{*}}} \approx
\lambda_{I_{3}}(\xi)>>\eta_{I_{3}}$. For $n=n_{*}$ and $I=\emptyset$, we have $\theta_{I}(\beta_{1})=\beta_{1}/2$, $\theta_{I\cup\{n\}}(\beta_{1})=1-\alpha_{n_{*}}+\beta_{1}/2$,
and $\psi_{I}=\eta_{I_{3}}$,
and we readily see from (\ref{glms}) that \medskip\\
\begin{subequations}\label{lm14}
$\displaystyle  \int B(\cdot,\xi)
C_{n_{*}}(x_{n_{*}})dx_{n_{*}}$
\begin{numcases}
{\lesssim}
\xi_{k_{*}}^{2{\overline{\ell}}_{k_{*}}(1+\beta_{1}/2)}
\min{\{
\,\Psi_{\{n_{*}\}}^{-\alpha_{n_{*}}},\,
\lambda_{I_{3}}^{-\alpha_{n_{*}}}\}},&
$1-\alpha_{n_{*}}+\beta_{1}/2>0$;\label{lm1}
\\[0.05 cm]
|\log{\max{\{
\Psi_{\{n_{*}\}},\,\lambda_{I_{3}\cup \{k_{*}\}}\}}}|,&
\nonumber\\
&\hspace{-3.2 cm} $1-\alpha_{n_{*}}+\beta_{1}/2=0$,\;\;$\xi_{k_{*}}^{2{
\overline{\ell}}_{k_{*}}}\gg  \Psi_{\{n_{*}\}}$;\label{lm2}
\\[0.05 cm]
\xi_{k_{*}}^{2{
\overline{\ell}}_{k_{*}}(1+\beta_{1}/2)}
\min{\{
\,\Psi_{\{n_{*}\}}^{-\alpha_{n_{*}}},\,
\lambda_{I_{3}}^{-\alpha_{n_{*}}}\}},&\nonumber\\
&\hspace{-3.0 cm}
$1-\alpha_{n_{*}}+\beta_{1}/2=0$,\;\;$\xi_{k_{*}}^{2{
\overline{\ell}}_{k_{*}}}\lesssim \Psi_{\{n_{*}\}}$;\label{lm3}
\\[0.05 cm]
\dfrac{\min{\{\xi_{k_{*}}^{2{
\overline{\ell}}_{k_{*}}},\,\max{\{
\Psi_{\{n_{*}\}},\,\lambda_{I_{3}\cup \{k_{*}\}}\}}\}^{1+\beta_{1}/2}}}
{\max{\{
\Psi_{\{n_{*}\}},\,\lambda_{I_{3}\cup \{k_{*}\}}\}}^{
\alpha_{n_{*}}}},&
$1-\alpha_{n_{*}}+\beta_{1}/2<0$.\label{lm4}
\end{numcases}
\end{subequations}
Starting from (\ref{lm14}), we estimate $\mathcal{B}$ in two cases according to whether $1-\alpha_{n_{*}}+\beta_{1}/2\geq0$ or
$1-\alpha_{n_{*}}+\beta_{1}/2<0$.\par
\textbullet\ First let $1-\alpha_{n_{*}}+\beta_{1}/2\geq0$.
Since $1+\beta_{0}/2>0$ and
$\xi_{k_{*}}^{2{\overline{\ell}}_{k_{*}}}
\approx\lambda_{I_{3}}(\xi)$, we have
\begin{equation}\label{nqmb}
\xi_{k_{*}}^{2{\overline{\ell}}_{k_{*}}(1+\beta_{1}/2)}
=\xi_{k_{*}}^{\overline{\ell}_{k_{*}}}
\xi_{k_{*}}^{2{\overline{\ell}}_{k_{*}}(1+\beta_{0}/2)}
\leq
\xi_{k_{*}}^{\overline{\ell}_{k_{*}}}
\max{\{
\,\Psi_{\{n_{*}\}}^{1+\beta_{0}/2},\,
\lambda_{I_{3}}^{1+\beta_{0}/2}\}}.
\end{equation}
By (\ref{balf1}), $1-\alpha_{n_{*}}+\beta_{0}/2<0$,
and, if $1-\alpha_{n_{*}}+\beta_{1}/2>0$, it follows from (\ref{lm1}) and (\ref{nqmb}) that
\begin{equation}\label{mbpz}
 \int B(\cdot,\xi)
C_{n_{*}}(x_{n_{*}})dx_{n_{*}}\lesssim
\xi_{k_{*}}^{\overline{\ell}_{k_{*}}}
\min{\{
\,\Psi_{\{n_{*}\}}^{1-\alpha_{n_{*}}+\beta_{0}/2},\,
\lambda_{I_{3}}^{1-\alpha_{n_{*}}+\beta_{0}/2}\}}.
\end{equation}
In view of (\ref{nstr}) and the last estimate, it follows from
Corollary \ref{cor43}, with
$j=0$, $\beta=\beta_{0}$ and $I_{*}=I_{3}$, that
\begin{equation}\label{mbpls}
\mathcal{B}(\xi)\lesssim\,
\xi_{k_{*}}^{\overline{\ell}_{k_{*}}}
\mathcal{O}_{I_{3}}(\xi;\beta_{0},I_{3}),
\end{equation}
when $1-\alpha_{n_{*}}+\beta_{1}/2>0$.
Combining the estimates
(\ref{lm3}) and (\ref{nqmb}), we see that the integral
\begin{equation*}
 \int B(\cdot,\xi)\chi_{\{\xi_{k_{*}}^{2{
\overline{\ell}}_{k_{*}}}\lesssim \Psi_{\{n_{*}\}}\}}
C_{n_{*}}(x_{n_{*}})dx_{n_{*}}
\end{equation*}
also satisfies the bound (\ref{mbpz}),
when $1-\alpha_{n_{*}}+\beta_{1}/2=0$.
Therefore, similarly to (\ref{mbpls}),
\begin{equation}\label{mblg0}
 \int B(\cdot,\xi)\chi_{\{\xi_{k_{*}}^{2{
\overline{\ell}}_{k_{*}}}\lesssim \Psi_{\{n_{*}\}}\}}
\prod_{n\in I_{3}}C_{n}(x_{n})dx_{n} \lesssim \xi_{k_{*}}^{\overline{\ell}_{k_{*}}}
\mathcal{O}_{I_{3}}(\xi;\beta_{0},I_{3}),
\end{equation}
when  $1-\alpha_{n_{*}}+\beta_{1}/2=0$, by Corollary \ref{cor43}.
So, when $1-\alpha_{n_{*}}+\beta_{1}/2=0$,
it suffices to estimate the integral
\begin{multline*}
 \int B(\cdot,\xi)\chi_{\{\xi_{k_{*}}^{2{
\overline{\ell}}_{k_{*}}}\gg \Psi_{\{n_{*}\}}\}}
\prod_{n\in I_{3}}C_{n}(x_{n})dx_{n}\\
\lesssim
 \int |\log{\max{\{
\Psi_{\{n_{*}\}},\,\lambda_{I_{3}\cup \{k_{*}\}}\}}}|
\chi_{\{\xi_{k_{*}}^{2{
\overline{\ell}}_{k_{*}}}\gg \Psi_{\{n_{*}\}}\}}
\prod_{n\in I_{3}\setminus\{n_{*}\}}C_{n}(x_{n})dx_{n},
\end{multline*}
by (\ref{lm2}). We consider two cases depending on whether
\text{(i)} $\xi^{\overline{\ell}_{k_{*}}}_{k_{*}}\gg
\xi^{\overline{\ell}_{n_{1}}}_{n_{1}}$ or
 \text{(ii)}  $\xi^{\overline{\ell}_{k_{*}}}_{k_{*}}\lesssim
\xi^{\overline{\ell}_{n_{1}}}_{n_{1}}$.\par
\text{(i)} Assume that $\xi^{\overline{\ell}_{k_{*}}}_{k_{*}}\gg
\xi^{\overline{\ell}_{n_{1}}}_{n_{1}}$. In the intersection of this
 region with the region (\ref{rgnn}), we have
\begin{equation*}
\xi^{2\overline{\ell}_{k_{*}}}_{k_{*}}\gg
\sum_{n_{1}\leq n\leq n_{N}}
\xi^{2\overline{\ell}_{n}}_{n}\gtrsim
\sum_{n\in I_{3}\setminus \{n_{*}\}}{
|x_{n}-\xi^{2\overline{\ell}_{n}}_{n}|},
\end{equation*}
since $\xi^{2\overline{\ell}_{n}}_{n}/4\leq
x_{n}\leq 4\xi^{2\overline{\ell}_{n}}_{n}$, for all
$n\in I_{3}$.
If $\xi_{k_{*}}^{2{
\overline{\ell}}_{k_{*}}}\gg \Psi_{\{n_{*}\}}$, then
\begin{eqnarray*}
\xi_{k_{*}}^{2{
\overline{\ell}}_{k_{*}}}&\gg &|1-\sum_{n\in I_{3}\setminus \{n_{*}\}}{
x_{n}}-{
\xi^{2{\overline{\ell}}_{n_{*}}}_{n_{*}}}-\xi^{
2{\overline{\ell}}_{d}}_{d}|\\
&=&|1-\sum_{n\in I_{3}\setminus \{n_{*}\}}{
(x_{n}-\xi^{2\overline{\ell}_{n}}_{n})}-
\sum_{n\in I_{3}}
{\xi^{2\overline{\ell}_{n}}_{n}}-\xi^{
2{\overline{\ell}}_{d}}_{d}|,
\end{eqnarray*}
whence $\xi_{k_{*}}^{2{
\overline{\ell}}_{k_{*}}}\gg \Phi_{I_{3}}(\xi)$. Since
$\xi_{k_{*}}^{2{
\overline{\ell}}_{k_{*}}}\approx\lambda_{I_{3}}(\xi)$,
we then have
\begin{equation}\label{rgn1}
\xi_{k_{*}}^{2{
\overline{\ell}}_{k_{*}}}\approx
\max{\{
\Phi_{I_{3}}(\xi), \lambda_{I_{3}}(\xi),
\xi^{2\overline{\ell}_{n_{1}}}_{n_{1}}\}}.
\end{equation}
Similarly to the way we obtained (\ref{j0mty1}), starting from
(\ref{lm2}) and then iterating
(\ref{i3t4}) with $(n,I)=$$(n_{1},O_{0})$,
$(n_{2}, O_{1})$, $\dots$, $(n_{N}, O_{N-1})$,
we see that the integral
\begin{multline}\label{mblg1}
 \int B(\cdot,\xi)
\chi_{\{\xi_{k_{*}}^{2{
\overline{\ell}}_{k_{*}}}\gg \Psi_{\{n_{*}\}}\}}
\prod_{n\in I_{3}}C_{n}(x_{n})dx_{n}\\
\lesssim|\log{\max{\{\Phi_{I_{3}}(\xi),\,
\lambda_{I_{3}\cup\{k_{*}\}}(\xi),\,
\xi_{n_{1}}^{2{\overline{\ell}}_{n_{1}}}\}}}|
\widetilde{P}_{I_{3}\setminus \{n_{*}\}}(\xi),
\end{multline}
if $
 \xi_{k_{*}}^{{
\overline{\ell}}_{k_{*}}}\gg
\xi_{n_{1}}^{{\overline{\ell}}_{n_{1}}},
$ and vanishes outside the region (\ref{rgn1}).
Observe that (\ref{mblg1}) holds when $1-\alpha_{n_{*}}+\beta_{1}/2=0$, i.e., when
$I^{-}_{3}(\beta_{1})=\emptyset$ and
$\theta_{I^{+}_{3}(\beta_{1})}(\beta_{1})=0$, which
puts us in the case \textnormal{\bf(c2)}\textnormal{
($I_{3}$,$\,I_{3}\cup\{k_{*}\}$,$\,\beta_{1}$)} and the right side of (\ref{mblg1}) equals
$\mathcal{O}_{I_{3}}(\xi;\beta_{1},I_{3}\cup\{k_{*}\})$
in that case.
\par
\text{(ii)} Assume now that $\xi^{\overline{\ell}_{k_{*}}}_{k_{*}}\lesssim
\xi^{\overline{\ell}_{n_{1}}}_{n_{1}}$. In the region (\ref{rgnn}),
we have $\xi^{2\overline{\ell}_{n_{1}}}_{n_{1}}\gtrsim
\sum_{n\in I_{3}\setminus \{n_{*}\}}{
|x_{n}-\xi^{2\overline{\ell}_{n}}_{n}|}$.
Therefore, if $\xi_{k_{*}}^{2{
\overline{\ell}}_{k_{*}}}\gg \Psi_{\{n_{*}\}}$, then
\begin{equation*}
\xi^{2\overline{\ell}_{n_{1}}}_{n_{1}}\gg
|1-\sum_{n\in I_{3}\setminus \{n_{*}\}}{
(x_{n}-\xi^{2\overline{\ell}_{n}}_{n})}-
\sum_{n\in I_{3}}
{\xi^{2\overline{\ell}_{n}}_{n}}-\xi^{
2{\overline{\ell}}_{d}}_{d}|,
\end{equation*}
which means that $\xi^{2\overline{\ell}_{n_{1}}}_{n_{1}}\gtrsim
\Phi_{I_{3}}(\xi)$. This shows that $\,\xi_{n_{1}}^{2{
\overline{\ell}}_{n_{1}}}$ $\gtrsim\max{\left\{
\Phi_{I_{3}}(\xi)\right.}$, ${\left.\lambda_{I_{3}}(\xi)\right\}}$, i.e.,
$n_{1}\in K^{+}(I_{3})$. Applying the estimate
(\ref{plmlg}) with
$n=n_{1}$ and $I=\{n_{*}\}$, we deduce that
\begin{multline}\label{smlg}
 \int |\log{\max{\{
\Psi_{\{n_{*}\}},\,\lambda_{I_{3}\cup \{k_{*}\}}\}}}|
\chi_{\{x_{n_{1}}\gtrsim\xi_{k_{*}}^{2{
\overline{\ell}}_{k_{*}}}\gg \Psi_{\{n_{*}\}}\}}
C_{n_{1}}(x_{n_{1}})dx_{n_{1}}\\[0.05 cm]
\begin{aligned}[b]
&\lesssim
\xi_{k_{*}}^{2{
\overline{\ell}}_{k_{*}}}|\log{\xi_{k_{*}}}|
\min{\{
\Psi_{O_{1}}^{-\alpha_{n_{1}}},
\,\lambda_{I_{3}}^{-\alpha_{n_{1}}}\}}\\[0.05 cm]
&\leq \xi_{k_{*}}^{{
\overline{\ell}}_{k_{*}}}|\log{\xi_{k_{*}}}|
\min{\{
\Psi_{O_{1}}^{{1}/{2}-\alpha_{n_{1}}},
\,\lambda_{I_{3}}^{{1}/{2}-\alpha_{n_{1}}}\}}\\[0.05 cm]
&=\xi_{k_{*}}^{{
\overline{\ell}}_{k_{*}}}|\log{\xi_{k_{*}}}|
\min{\{
\Psi_{O_{1}}^{\theta_{O_{1}}(\beta_{0})},
\,\lambda_{I_{3}}^{\theta_{O_{1}}(\beta_{0})}\}},
\end{aligned}
\end{multline}
when $1-\alpha_{n_{*}}+\beta_{1}/2=0$. Combining this with (\ref{lm2}),
we see that
\begin{multline*}
 \int  \int B(\cdot,\xi)\chi_{\{\xi_{k_{*}}^{2{
\overline{\ell}}_{k_{*}}}\gg \Psi_{\{n_{*}\}}\}}
C_{n_{*}}(x_{n_{*}})C_{n_{1}}(x_{n_{1}})dx_{n_{*}}dx_{n_{1}}\\ \lesssim
\xi_{k_{*}}^{{
\overline{\ell}}_{k_{*}}}|\log{\xi_{k_{*}}}|
\min{\{
\Psi_{O_{1}}^{\theta_{O_{1}}(\beta_{0})},
\,\lambda_{I_{3}}^{\theta_{O_{1}}(\beta_{0})}(\xi)\}}.
\end{multline*}
When $1-\alpha_{n_{*}}+\beta_{1}/2=0$, we have
\begin{equation*}
\theta_{O_{1}}(\beta_{0})=
2-\alpha_{n_{*}}-\alpha_{n_{1}}+\beta_{0}/2=1/2-\alpha_{n_{1}}<0.
\end{equation*}
So $ O_{1}\subseteq
I^{-}_{3}(\beta_{0})$. We have already noted that
$n_{1}\in K^{+}(I_{3})$ in this case. Thus,
$ O_{1}\subseteq
K^{+}(I_{3})\cap I^{-}_{3}(\beta_{0})$. Using the latter estimate
and then invoking Corollary \ref{cor43}, with $j=1$, $\beta=\beta_{0}$,
and $I_{*}=I_{3}$, we find that
\begin{equation}\label{mblg2}
\int B(\cdot,\xi)
\chi_{\{\xi_{k_{*}}^{2{
\overline{\ell}}_{k_{*}}}\gg \Psi_{\{n_{*}\}}\}}
\prod_{n\in I_{3}}
C_{n}(x_{n})dx_{n}
 \lesssim
\xi_{k_{*}}^{{
\overline{\ell}}_{k_{*}}}|\log{\xi_{k_{*}}}|
\mathcal{O}_{I_{3}}(\xi;\beta_{0},I_{3}),
\end{equation}
in the region $\lambda_{I_{3}}(\xi)\approx\xi^{2\overline{\ell}_{k_{*}}}_{k_{*}}\lesssim
\xi^{2\overline{\ell}_{n_{1}}}_{n_{1}}$.\smallskip\par
\textbullet\ Suppose $1-\alpha_{n_{*}}+\beta_{1}/2<0$.
Then $n_{*}\in I^{-}_{3}(\beta_{1})$.
Since $\lambda_{I_{3}}\geq \lambda_{I_{3}\cup
\{k_{*}\}}$, we have $K^{+}(I_{3})\subseteq
K^{+}(I_{3}\cup\{k_{*}\})$, and it follows from (\ref{nstr}) that
\begin{equation}\label{nstr1}
n_{*}\in K^{+}(I_{3}\cup
\{k_{*}\})\cap I^{-}_{3}(\beta_{1}).
\end{equation}
We have two cases to consider:
\begin{center}
\text{(i)} $\xi^{2{\overline{\ell}}_{k_{*}}}_{k_{*}}\lesssim
\Phi_{I_{3}}(\xi)$, \hspace{3 cm}
\text{(ii)} $\xi^{2{\overline{\ell}}_{k_{*}}}_{k_{*}}\gg
\Phi_{I_{3}}(\xi)$.
\end{center}
\par
\text{(i)} Let $
\xi^{2{\overline{\ell}}_{k_{*}}}_{k_{*}}\lesssim
\Phi_{I_{3}}(\xi)$. Since $1+\beta_{0}/2>0$, we have
\begin{equation*}
\min{\{\xi_{k_{*}}^{2{
\overline{\ell}}_{k_{*}}},\max{\{
\Psi_{\{n_{*}\}},\lambda_{I_{3}\cup
\{k_{*}\}}\}}\}^{1+\beta_{1}/2}} \leq
\xi_{k_{*}}^{\overline{\ell}_{k_{*}}}
\max{\{
\Psi_{\{n_{*}\}},\lambda_{I_{3}\cup
\{k_{*}\}}\}^{1+\beta_{0}/2}}.
\end{equation*}
By (\ref{balf1}), it follows from (\ref{lm4}) that
\begin{equation*}
 \int B(\cdot,\xi)
C_{n_{*}}(x_{n_{*}})dx_{n_{*}}\lesssim
\xi_{k_{*}}^{\overline{\ell}_{k_{*}}}
\min{\{
\Psi_{\{n_{*}\}}^{1-\alpha_{n_{*}}+\beta_{0}/2},
\lambda_{I_{3}\cup
\{k_{*}\}}^{1-\alpha_{n_{*}}+\beta_{0}/2}(\xi)\}}.
\end{equation*}
Using (\ref{nstr1}) and the fact that
$I^{-}_{3}(\beta_{1})
\subseteq I^{-}_{3}(\beta_{0})$, one can apply
 Corollary \ref{cor43} again with $j=0$, $\beta=\beta_{0}$
and $I_{*}=I_{3}\cup\{k_{*}\}$, to obtain
\begin{equation}\label{mb3i}
\mathcal{B}(\xi)\,\lesssim\,
\xi_{k_{*}}^{\overline{\ell}_{k_{*}}}
\mathcal{O}_{I_{3}}(\xi;\beta_{0},I_{3}\cup\{k_{*}\})
\,\approx\,\xi_{k_{*}}^{\overline{\ell}_{k_{*}}}
\mathcal{O}_{I_{3}}(\xi;\beta_{0},I_{3}).
\end{equation}
\smallskip\par
\text{(ii)} Now let $\xi^{2{\overline{\ell}}_{k_{*}}}_{k_{*}}\gg
\Phi_{I_{3}}(\xi)$. By the assumption that
$1-\alpha_{n_{*}}+\beta_{1}/2<0$,
it follows from (\ref{lm4}) that
\begin{equation*}
 \int B(\cdot,\xi)
C_{n_{*}}(x_{n_{*}})dx_{n_{*}}\lesssim
\min{\{
\Psi_{\{n_{*}\}}^{1-\alpha_{n_{*}}+\beta_{1}/2},
\,\lambda_{I_{3}\cup
\{k_{*}\}}^{1-\alpha_{n_{*}}+\beta_{1}/2}(\xi)\}}\}.
\end{equation*}
By virtue of (\ref{nstr1}) and Corollary
\ref{cor43}, with $j=0$, $\beta=\beta_{1}$
and $I_{*}=I_{3}\cup\{k_{*}\}$, it follows from
the latter inequality that
\begin{equation}\label{mb3ii}
\mathcal{B}(\xi)\lesssim
\mathcal{O}_{I_{3}}(\xi;\beta_{1},I_{3}\cup\{k_{*}\}),
\quad \xi^{2\overline{\ell}_{k_{*}}}_{k_{*}}\approx
\lambda_{I_{3}}(\xi)\gg
\Phi_{I_{3}}(\xi).
\end{equation}
This estimate is useful in the region
$\{\xi_{k_{*}}^{2{
\overline{\ell}}_{k_{*}}}\gtrsim
\Lambda_{I_{3}}(\xi;\beta_{1},I_{3}\cup\{k_{*}\})\}$. However, due to the factor
$\xi^{-\alpha_{k_{*}}}_{k_{*}}$ in (\ref{r3s1}),
the estimate (\ref{mb3ii}) is insufficient in the
regions\footnote{\thanks{We refer to the notations introduced
before Corollary \ref{cor43}.}}
\begin{equation*}
\begin{cases}
\xi^{\overline{\ell}_{k_{*}}}_{k_{*}}\ll
\xi_{n_{N_{2}(\beta_{1})}}^{
\overline{\ell}_{n_{N_{2}(\beta_{1})}}},&
\textnormal{\bf(c2)}\textnormal{($I_{3}$,$\,I_{3}\cup\{k_{*}\}$,$\,\beta_{1}$)};\\[0.2 cm]
\xi^{\overline{\ell}_{k_{*}}}_{k_{*}}\ll
\xi_{n_{N_{1}(\beta_{1})}}^{
\overline{\ell}_{n_{N_{1}(\beta_{1})}}},&
\textnormal{\bf(c3)}\textnormal{($I_{3}$,$\,I_{3}
\cup\{k_{*}\}$,$\,\beta_{1}$)},
\end{cases}
\end{equation*}
i.e. in the regions where
\begin{subequations}\label{smptkn}
\begin{numcases}
{}
\label{smptkn1} \xi^{\overline{\ell}_{k_{*}}}_{k_{*}}\ll
\xi^{\overline{\ell}_{n}}_{n},\; \forall
n\in I^{+}_{3}(\beta_{1})\cup\{n_{N_{2}(\beta_{1})}\},&
$\textnormal{\bf(c2)}\textnormal{($I_{3}$,$\,I_{3}\cup\{k_{*}\}$,$\,\beta_{1}$)}$;\\[0.05 cm]
\label{smptkn2} \xi^{\overline{\ell}_{k_{*}}}_{k_{*}}\ll
\xi^{\overline{\ell}_{n}}_{n},\; \forall
n\in I^{+}_{3}(\beta_{1}),&
$\textnormal{\bf(c3)}\textnormal{($I_{3}$,$\,I_{3}\cup\{k_{*}\}$,$\,\beta_{1}$)}$,
\end{numcases}
\end{subequations}
The key to a better estimate for $\mathcal{B}$
in the regions (\ref{smptkn}), under the assumptions that
$\xi^{2{\overline{\ell}}_{k_{*}}}_{k_{*}}\gg
\Phi_{I_{3}}(\xi)$ and $\xi^{2{\overline{\ell}}_{k_{*}}}_{k_{*}}
\approx\lambda_{I_{3}}(\xi)$, is using
the full strength of the estimates (\ref{glms}).
We shall see that
\begin{numcases}
{\mathcal{B}(\xi)\lesssim}
\xi_{k_{*}}^{{
\overline{\ell}}_{k_{*}}}|\log{\xi_{k_{*}}}|
\mathcal{O}_{I_{3}}(\xi;\beta_{0},I_{3}),&\nonumber\\
\label{mb3f1}&\hspace{-2.5 cm}
\textnormal{\bf(c2)}\textnormal{($I_{3}$,$\,I_{3}\cup\{k_{*}\}$,$\,\beta_{1}$)},\;
$\xi^{\overline{\ell}_{k_{*}}}_{k_{*}}\ll
\xi_{n_{N_{2}(\beta_{1})}}^{
\overline{\ell}_{n_{N_{2}(\beta_{1})}}}$;\\
\xi_{k_{*}}^{{
\overline{\ell}}_{k_{*}}}
\mathcal{O}_{I_{3}}(\xi;\beta_{0},I_{3}),&\nonumber\\
&\hspace{-2.5 cm}
\textnormal{\bf(c3)}\textnormal{($I_{3}$,$\,I_{3}\cup\{k_{*}\}$,$\,\beta_{1}$)},\;
$\xi^{\overline{\ell}_{k_{*}}}_{k_{*}}\ll
\xi_{n_{N_{1}(\beta_{1})}}^{
\overline{\ell}_{n_{N_{1}(\beta_{1})}}}$.\label{mb3f2}
\end{numcases}
Observe that when
$\xi^{2{\overline{\ell}}_{k_{*}}}_{k_{*}}\approx\lambda_{I_{3}}(\xi)\gg
\Phi_{I_{3}}(\xi)$, it follows from (\ref{smptkn1}) and (\ref{smptkn2}),
respectively, that
\begin{subequations}\label{trunsm}
\begin{numcases}
{}
\label{trunsm1}
I^{+}_{3}(\beta_{1})\cup\{n_{N_{2}(\beta_{1})}\}
\subseteq
K^{+}(I_{3}),&
$\textnormal{\bf(c2)}\textnormal{($I_{3}$,$\,I_{3}\cup\{k_{*}\}$,$\,\beta_{1}$)}$;\\[0.05 cm]
\label{trunsm2} I^{+}_{3}(\beta_{1})\subseteq
K^{+}(I_{3}),&
$\textnormal{\bf(c3)}\textnormal{($I_{3}$,$\,I_{3}\cup\{k_{*}\}$,$\,\beta_{1}$)}$.
\end{numcases}
\end{subequations}
For each $I\subseteq I_{3}$, let
\begin{eqnarray*}
B^{+}_{I}(x,\xi)&:=&\min{\{\Psi^{\theta_{I}(\beta_{1})}_{I},\,
\lambda^{\theta_{I}(\beta_{1})}_{I_{3}\cup \{k_{*}\}}
\}}
\chi_{\{
\xi^{2{\overline{\ell}}_{k_{*}}}_{k_{*}}\gg \Psi_{I}
\}}, \\
B^{-}_{I}(x,\xi)&:=&
\min{\{
\Psi_{I}^{\theta_{I}(\beta_{0})},
\,\lambda_{I_{3}\cup
\{k_{*}\}}^{\theta_{I}(\beta_{0})}(\xi)\}}
\chi_{\{\xi^{2{\overline{\ell}}_{k_{*}}}_{k_{*}}\lesssim
\Psi_{I}\}}.
\end{eqnarray*}
Since $1+\beta_{0}/2>0$ and $\xi^{2{\overline{\ell}}_{k_{*}}}_{k_{*}}\gtrsim
\lambda_{I_{3}\cup \{k_{*}\}}$,
we have from (\ref{glm4}) the pointwise estimate
\begin{multline}\label{mbb1}
\int
B^{+}_{I}(x,\xi)
\chi_{\{x_{n}\gg \xi^{2{\overline{\ell}}_{k_{*}}}_{k_{*}}
\gtrsim
\max{\{\Psi_{I},\,\lambda_{I_{3}\cup\{k_{*}\}}\}}\}}
C_{n}(x_{n})dx_{n}\\
\lesssim
B^{+}_{I\cup\{n\}}(x,\xi)+
\xi^{\overline{\ell}_{k_{*}}}_{k_{*}}
B^{-}_{I\cup\{n\}}(x,\xi),
\end{multline}
for every $n\in I_{3}$ and $I\subseteq I_{3}\setminus\{n\}$
such that $\theta_{I\cup \{n\}}(\beta_{1})<0$. We  would like to
iterate (\ref{mbb1}) for the pairs
$(n,I)=$ $(n_{*},\emptyset)$,
$(n_{1},\{n_{*}\})$,
$(n_{2},\{n_{*},n_{1}\})$,
$\dots$,$(n_{N^{-}(\beta_{1})},
I^{-}_{3}(\beta_{1})\setminus\{n_{N^{-}(\beta_{1})}\})$.
The restriction that
$x_{n}\gg \xi^{2{\overline{\ell}}_{k_{*}}}_{k_{*}}$
in the support of $C_{n}(x_{n})$, i.e., wherever
$\xi^{2\overline{\ell}_{n}}_{n}/4\leq
x_{n}\leq 4\xi^{2\overline{\ell}_{n}}_{n}$,
is imposed by (\ref{smptkn}) every time we run (\ref{mbb1}).
We will show later that
\begin{equation}\label{mimns}
\int B^{-}_{I}(x,\xi)
\prod_{n\in I_{3}\setminus I}
C(x_{n}) dx_{n}\lesssim \mathcal{O}_{I_{3}}(\xi;\beta_{0},I_{3}).
\end{equation}
for each $I\in \{\{n_{*}\},
\{n_{*},n_{1}\},\dots, I^{-}_{3}(\beta_{1})\}$.
Thus, we only need to estimate
\begin{equation*}
\int
B^{+}_{I^{-}_{3}(\beta_{1})}(x,\xi)
\chi_{\{x_{n_{N_{1}(\beta_{1})}}\gg
\xi^{2{\overline{\ell}}_{k_{*}}}_{k_{*}}\approx\lambda_{I_{3}}\gg
\Psi_{I^{-}_{3}(\beta_{1})}\}}
\prod_{n\in I_{3}\setminus I^{-}_{3}(\beta_{1})}
C(x_{n}) dx_{n}.
\end{equation*}
\par Consider first the case \textnormal{\bf(c2)}\textnormal{($I_{3}$,$\,I_{3}\cup\{k_{*}\}$,$\,\beta_{1}$)}
under the assumption (\ref{trunsm1}).
In this case $\theta_{I^{+}_{3}(\beta_{1})}(\beta_{1})=0$.
Using (\ref{glm2}) and (\ref{glm3})
with $I=I^{-}_{3}(\beta_{1})$ and
$n=n_{N_{1}(\beta_{1})}$, we obtain
\begin{multline}\label{minn1}
\int B^{+}_{I^{-}_{3}(\beta_{1})}(x,\xi)
\chi_{\{x_{n_{N_{1}(\beta_{1})}}\gg
\xi^{2{\overline{\ell}}_{k_{*}}}_{k_{*}}\approx\lambda_{I_{3}}\gg
\Psi_{I^{-}_{3}(\beta_{1})}\}}
C(x_{n_{N_{1}(\beta_{1})}}) dx_{n_{N_{1}(\beta_{1})}}\\
\lesssim
|\log{\max{\{
\Psi_{I^{+}_{3}(\beta_{1})},
\,\lambda_{I_{3}\cup \{k_{*}\}}\}}}|
\chi_{\{\xi_{n_{N_{2}(\beta_{1})}}^{2{
\overline{\ell}}_{n_{N_{2}(\beta_{1})}}}\gg
\xi_{k_{*}}^{2{
\overline{\ell}}_{k_{*}}}\approx\lambda_{I_{3}}\gg
 \Psi_{I^{+}_{3}(\beta_{1})}\}}\\
+\xi_{k_{*}}^{\overline{\ell}_{k_{*}}}
\min{\{
\,\Psi_{I^{+}_{3}(\beta_{1})}^{
\theta_{I^{+}_{3}(\beta_{1})}(\beta_{0})},\,
\lambda_{I_{3}}^{\theta_{I^{+}_{3}(\beta_{1})}
(\beta_{0})}\}}.
\end{multline}
Similarly to (\ref{smlg}), using (\ref{plmlg}) with
$n=n_{N_{2}(\beta_{1})}$ and $I=I^{+}_{3}(\beta_{1})$,
we see that the integral of the first term on the right side of
(\ref{minn1}) in $x_{n_{N_{2}(\beta_{1})}}$ against
$C_{n_{N_{2}(\beta_{1})}}(x_{n_{N_{2}(\beta_{1})}})$ is bounded pointwise
by
\begin{align}
\nonumber &\xi_{k_{*}}^{2{
\overline{\ell}}_{k_{*}}}|\log{\xi_{k_{*}}}|
\min{\{
\Psi_{I^{+}_{3}(\beta_{1})\cup\{n_{N_{2}(\beta_{1})}\}}^{
-\alpha_{n_{N_{2}(\beta_{1})}}},
\,\lambda_{I_{3}}^{-\alpha_{n_{N_{2}(\beta_{1})}}}\}}\\
\label{eqfn}&\leq \xi_{k_{*}}^{{
\overline{\ell}}_{k_{*}}}|\log{\xi_{k_{*}}}|
\min{\{
\Psi_{I^{+}_{3}(\beta_{1})\cup\{n_{N_{2}(\beta_{1})}\}}^{
\theta_{I^{+}_{3}(\beta_{1})\cup
\{n_{N_{2}(\beta_{1})}\}}(\beta_{0})},
\,\lambda_{I_{3}}^{
\theta_{I^{+}_{3}(\beta_{1})\cup
\{n_{N_{2}(\beta_{1})}\}}(\beta_{0})}\}},
\end{align}
since $\xi_{k_{*}}^{2{\overline{\ell}}_{k_{*}}}\approx
\lambda_{I_{3}}$ and
$1/2-\alpha_{n_{N_{2}(\beta_{1})}}=
\theta_{I^{+}_{3}(\beta_{1})\cup
\{n_{N_{2}(\beta_{1})}\}}(\beta_{0})$,
if $\theta_{I^{+}_{3}(\beta_{1})}(\beta_{1})=0$.
This equality also shows that $
\theta_{I^{+}_{3}(\beta_{1})\cup
\{n_{N_{2}(\beta_{1})}\}}(\beta_{0})<0$, whence
$I^{+}_{3}(\beta_{1})\cup\{n_{N_{2}(\beta_{1})}\}\subseteq
I^{-}_{3}(\beta_{0})$.
By this and (\ref{trunsm1}), Corollary \ref{cor43}, with
$j=N_{2}(\beta_{1})$, $\beta=\beta_{0}$ and
$I_{*}=I_{3}$, asserts that
the integral against
$\prod_{n\in I_{3}\setminus (
I^{+}_{3}(\beta_{1})\cup
\{n_{N_{2}(\beta_{1})}\})}C_{n}(x_{n})$ of (\ref{eqfn})
is majorized pointwise by
\begin{equation*}
\xi_{k_{*}}^{{
\overline{\ell}}_{k_{*}}}|\log{\xi_{k_{*}}}|
\mathcal{O}_{I_{3}}(\xi;\beta_{0},I_{3}).
\end{equation*}
We turn to the second term
on the right side of  (\ref{minn1}). We have already seen that
$I^{+}_{3}(\beta_{1})\subseteq I^{-}_{3}(\beta_{0})$.
Once again, by (\ref{trunsm1}),
Corollary \ref{cor43} with $j=N_{1}(\beta_{1})$, $\beta=\beta_{0}$
and $I_{*}=I_{3}$ yields that
\begin{equation*}
\int \min{\{
\,\Psi_{I^{+}_{3}(\beta_{1})}^{\theta_{I^{+}_{3}(\beta_{1})}(\beta_{0})},\,
\lambda_{I_{3}}^{\theta_{I^{+}_{3}(\beta_{1})}
(\beta_{0})}\}}\prod_{n\in I_{3}\setminus
I^{+}_{3}(\beta_{1})}C_{n}(x_{n})dx_{n}\lesssim
\mathcal{O}_{I_{3}}(\xi;\beta_{0},I_{3}).
\end{equation*}
Assuming for the moment the estimates (\ref{mimns}),
this proves (\ref{mb3f1}).
\smallskip
\par Consider next the case
 \textnormal{\bf(c3)}\textnormal{($I_{3}$,$\,I_{3}\cup\{k_{*}\}$,$\,\beta_{1}$)}
under the assumption (\ref{trunsm2}). Then
 $\theta_{I^{+}_{3}(\beta_{1})}(\beta_{1})>0$.
Employing the estimate (\ref{glm1}) with $I=I^{-}_{3}(\beta_{1})$ and
$n=n_{N_{1}(\beta_{1})}$, we find that the integral
on the left side of (\ref{minn1}) is majorized by
\begin{multline}\label{lstfn}
\xi_{k_{*}}^{2{\overline{\ell}}_{k_{*}}
(1+\theta_{I^{-}_{3}(\beta_{1})}(\beta_{1}))}
\min{\{
\,\Psi_{I^{+}_{3}(\beta_{1})}^{-\alpha_{n_{N_{1}(\beta_{1})}}},\,
\lambda_{I_{3}}^{-\alpha_{n_{N_{1}(\beta_{1})}}}\}}\\
\lesssim
\xi_{k_{*}}^{\overline{\ell}_{k_{*}}}
\min{\{
\,\Psi_{I^{+}_{3}(\beta_{1})}^{
\theta_{I^{+}_{3}(\beta_{1})}(\beta_{0})},\,
\lambda_{I_{3}}^{
\theta_{I^{+}_{3}(\beta_{1})}(\beta_{0})}\}}.
\end{multline}
Using Corollary \ref{cor43} as usual,
with $j=N_{1}(\beta_{1})$, $\beta=\beta_{0}$
and $I_{*}=I_{3}$, we see that
the integral of
(\ref{lstfn}) against $\prod_{n\in I_{3}\setminus
I^{+}_{3}(\beta_{1})}C_{n}(x_{n})$ is majorized by $\xi_{k_{*}}^{\overline{\ell}_{k_{*}}}
\mathcal{O}_{I_{3}}(\xi;\beta_{0},I_{3})$. This proves
(\ref{mb3f2}), assuming that (\ref{mimns}) holds.
\smallskip\par
To conclude the proof of (\ref{mb3f1}) and (\ref{mb3f2}),
it remains to prove (\ref{mimns}).
Fix $(n,I)$ in the set of pairs
\begin{equation*}
S
=\left\{(n_{1},\{n_{*}\}), (n_{2},\{n_{*},n_{1}\}),
\dots,(n_{N_{1}(\beta_{1})},I^{-}_{3}(\beta_{1}))\right\}.
\end{equation*}
Since, by definition,
$0>\theta_{I^{-}_{3}(\beta_{1})}(\beta_{1})=
\theta_{I^{-}_{3}(\beta_{1})}(\beta_{0})+1/2$, we have
\begin{equation*}
\theta_{I^{+}_{3}(\beta_{1})}(\beta_{0})=
\theta_{I^{-}_{3}(\beta_{1})}(\beta_{0})+
1-\alpha_{n_{N_{1}(\beta_{1})}}<
1/2-\alpha_{n_{N_{1}(\beta_{1})}}<0.
\end{equation*}
This implies that
\begin{equation}\label{ngxp}
\theta_{\{n_{*}\}}(\beta_{0})<
\theta_{\{n_{*},n_{1}\}}(\beta_{0})<\dots<
\theta_{I^{-}_{3}(\beta_{1})}(\beta_{0})<
\theta_{I^{+}_{3}(\beta_{1})}(\beta_{0})<0.
\end{equation}
We estimate
the integral of $B^{-}_{I}(x,\xi)$ in $x_{n}$ against
$C_{n}(x_{n})$, using (\ref{tvum}). Applying (\ref{tvum}) once with
$s=-\theta_{I}(\beta_{0})$, $\alpha=\alpha_{n}$, $\nu=
\xi_{n}^{2{
\overline{\ell}}_{n}}$, $\eta=1-
\sum_{j\in I}\xi^{
2{\overline{\ell}}_{j}}_{j}-\xi^{
2{\overline{\ell}}_{d}}_{d}-
\sum_{k\in I_{3}\setminus (I\cup\{n\})} x_{k}$,
$\lambda=\lambda_{I_{3}\cup\{k_{*}\}}$,
 and $\mu=\xi_{k_{*}}^{2{
\overline{\ell}}_{k_{*}}}$, we see that
\begin{equation}\label{prvrq}
\int B^{-}_{I}(x,\xi)
C(x_{n}) dx_{n}\lesssim
\min{\{
\Psi_{I\cup\{n\}}^{\theta_{I\cup\{n\}}(\beta_{0})},
\,\lambda_{I_{3}}^{\theta_{I\cup\{n\}}(\beta_{0})}(\xi)\}}.
\end{equation}
The inequalities (\ref{mimns}) follow
then from the latter estimate and Corollary \ref{cor43} with
$j=n$, $\beta=\beta_{0}$ and $I_{*}=I_{3}$. Observe that,
for each $(n,I)\in S$, $I\cup\{n\}\subseteq I^{-}_{3}(\beta_{0})$ by (\ref{ngxp}), and
$I\cup\{n\}\subseteq K^{+}(I_{3})$ by (\ref{trunsm}).\smallskip\par
Collecting the estimates
(\ref{mbpls}), (\ref{mblg0}), (\ref{mblg1}),
(\ref{mblg2}), (\ref{mb3i}), (\ref{mb3ii}),
(\ref{mb3f1}) and (\ref{mb3f2}) for $\mathcal{B}$ we find that, in the case \textbf{(b1)},
\begin{equation}\label{bb1}
\mathcal{B}(\xi)\lesssim\,
\xi_{k_{*}}^{{
\overline{\ell}}_{k_{*}}}|\log{\xi_{k_{*}}}|\,
\mathcal{O}_{I_{3}}(\xi;\beta_{0},I_{3})
+\mathcal{O}^{+}_{I_{3}}(\xi;\beta_{1},I_{3}\cup\{k_{*}\}),
\end{equation}
where
\begin{equation*}
\mathcal{O}^{+}_{I_{3}}(\xi;\beta_{1},I_{3}\cup\{k_{*}\}):=\,
\mathcal{O}_{I_{3}}(\xi;\beta_{1},I_{3}\cup\{k_{*}\})
\chi_{\{\xi_{k_{*}}^{2{
\overline{\ell}}_{k_{*}}}\gtrsim
\Lambda_{I_{3}}(\xi;\beta_{1},I_{3}\cup\{k_{*}\})\}}(\xi).
\end{equation*}
\smallskip\par
We turn to $\mathcal{B}$,
$-2<\beta_{0}<-1$, in the case \textbf{(b2)}.
In this case $J^{+}$ is a singleton.
As pointed out in the
proof of (\ref{fc2}), we necessarily have $J^{+}=\{k_{*}\}=\{r_{*}\}$
in this case. This restricts us to the region where
\begin{equation*}
\xi^{2{\overline{\ell}}_{k_{*}}}_{k_{*}}\gg
\max{\{\Phi_{I_{3}}(\xi),\,
\lambda_{I_{3}\cup\{k_{*}\}},\,
\xi^{2{\overline{\ell}}_{n_{*}}}_{n_{*}}\}}.
\end{equation*}
Using Lemma \ref{lmt2} with
$\sigma=\beta_{1}$, $l_{0}=\emptyset$, $l=I_{3}$ and
$I_{*}=I_{3}\cup\{k_{*}\}$, we see that
the estimate (\ref{mb3ii}) for $\mathcal{B}$
in the case \textbf{(b1)} holds
in the case \textbf{(b2)} as well.\smallskip\par
It only remains to estimate
$\mathcal{C}$, when $\beta_{0}=-1$, in the case \textbf{(b1)}.
Analogously to (\ref{smlg}) and (\ref{eqfn}), using
(\ref{plmlg}) with $n=n_{*}$ and $I=\emptyset$,
we obtain
\begin{multline*}
\int |\log{\max{\{
\eta_{I_{3}},\,\lambda_{I_{3}\cup \{k_{*}\}}\}}}|\,\chi_{\{
x_{n_{*}}\gtrsim
\xi^{2{\overline{\ell}}_{k_{*}}}_{k_{*}}\gg
\max{\{\eta_{I_{3}},\,
\xi_{j_{*}}^{2{\overline{\ell}}_{j_{*}}}\}}\}}
C_{n_{*}}(x_{n_{*}})d x_{n_{*}}\\
\lesssim
\xi_{k_{*}}^{2{
\overline{\ell}}_{k_{*}}}|\log{\xi_{k_{*}}}|
\min{\{\Psi_{\{n_{*}\}}^{-\alpha_{n_{*}}},
\,\lambda_{I_{3}}^{-\alpha_{n_{*}}}(\xi)\}}.
\end{multline*}
Since $\beta_{0}=-1$, we have $\xi^{\overline{\ell}_{k_{*}}}_{k_{*}}\approx
\lambda^{1+\beta_{0}/2}_{I_{3}}(\xi)$, whence
\begin{equation*}
\int C(\cdot,\xi)
C_{n_{*}}(x_{n_{*}})d x_{n_{*}}\lesssim
\xi_{k_{*}}^{{
\overline{\ell}}_{k_{*}}}|\log{\xi_{k_{*}}}|
\min{\{\Psi_{\{n_{*}\}}^{1-\alpha_{n_{*}}+\beta_{0}/2},
\,\lambda_{I_{3}}^{1-\alpha_{n_{*}}+\beta_{0}/2}(\xi)\}}.
\end{equation*}
By Corollary \ref{cor43} one last time, with
$j=0$, $\beta=\beta_{0}$ and $I_{*}=I_{3}$, we have
\begin{equation}\label{mcf}
\mathcal{C}(\xi)
\lesssim \xi_{k_{*}}^{{
\overline{\ell}}_{k_{*}}}|\log{\xi_{k_{*}}}|
\mathcal{O}_{I_{3}}(\xi;\beta_{0},I_{3}).
\end{equation}
\smallskip \par Using the estimates (\ref{ma}), (\ref{mam12}) with $m=0$, (\ref{bb1}) and (\ref{mcf}) in (\ref{misb1}), we deduce that \smallskip \\
\textbullet\ if $\Gamma^{+}_{d}=J^{+}=\emptyset$, then
\begin{subequations}\label{mif123}
\begin{equation}\label{mif1}
\mathcal{K}_{I_{1}, I_{2}, I_{3}}(\xi)\lesssim
\begin{cases}
\widetilde{\mathcal{O}}_{I_{3}}(\xi),&
-2<\beta_{0}<-1,\\[0.05 cm]
|\log{\xi_{k_{*}}}|\,
\mathcal{O}_{I_{3}}(\xi;\beta_{0},I_{3})\tau(\xi),&
\beta_{0}=-1,\\[0.05 cm]
\mathcal{O}_{I_{3}}(\xi;\beta_{0},I_{3})\tau(\xi),&
-1<\beta_{0}<0,
\end{cases}
\end{equation}
where $\widetilde{\mathcal{O}}_{I_{3}}(\xi)$ denotes
\begin{equation*}
|\log{\xi_{k_{*}}}|
\mathcal{O}_{I_{3}}(\xi;\beta_{0},I_{3})
\tau(\xi)+
\mathcal{O}^{+}_{I_{3}}(\xi;\beta_{1},I_{3}\cup\{k_{*}\})
\tau_{k_{*}}(\xi).
\end{equation*}
Inserting the estimate (\ref{mb3ii}) for $\mathcal{B}$
in the case \textbf{(b2)} into (\ref{misb2}), we find that \smallskip\\
\textbullet\ if $\Gamma^{+}_{d}=\emptyset$, $|J^{+}|=1$, and $-2<\beta_{0}<-1$, then
\begin{equation}\label{mif2}
\mathcal{K}_{I_{1}, I_{2}, I_{3}}(\xi)\lesssim
 \mathcal{O}_{I_{3}}(\xi;\beta_{1},I_{3}\cup\{k_{*}\})\tau_{k_{*}}(\xi),
\end{equation}
and we necessarily are in the region $\xi^{2{\overline{\ell}}_{k_{*}}}_{k_{*}}\gg
\Lambda_{I_{3}}(\xi;\beta_{1},I_{3}\cup\{k_{*}\})$.
Finally, it follows from
(\ref{mam12}), with $m=1$, and (\ref{misb3}) that \smallskip\\
\textbullet\ if $|\Gamma^{+}_{d}|=1$, $J^{+}=\emptyset$,
and $-1<\beta_{1}<0$, then
\begin{equation}\label{mif3}
\mathcal{K}_{I_{1}, I_{2}, I_{3}}(\xi)\lesssim
 \mathcal{O}_{I_{3}}(\xi;\beta_{1},I_{3})\tau(\xi).
\end{equation}
\end{subequations}
\par It is only left to substitute the estimates (\ref{mif123}) into formula (\ref{r3s1}). We also
substitute for $\mathcal{O}_{I_{3}}$
and for $\Lambda_{I_{3}}$ in each of the cases
\textbf{(c1)} through \textbf{(c3)}.
For every $I\subseteq \Gamma^{-}_{d}$, let
\begin{equation*}
\lambda^{*}_{I}(\xi)=
\begin{cases}
\lambda_{I}(\xi),& \text{if}\;I_{*}=I_{3};\\
\lambda_{I\cup \{k_{*}\}}(\xi),& \text{if}\;I_{*}=I_{3}\cup \{k_{*}\}.
\end{cases}
\end{equation*}
Take into account that
\begin{equation*}
\begin{cases}
I^{+}_{3}(\beta)=
K^{+}(I_{*})\cap I^{+}_{3}(\beta),&\hspace{0.5 cm} \text{\bf(c2)};
\\[0.1 cm]
I^{-}_{3}(\beta)=
K^{+}(I_{*})\cap I^{-}_{3}(\beta),&\hspace{0.5 cm}  \text{\bf(c3)},
\end{cases}
\end{equation*}
and that
\begin{align*}
\Lambda_{I_{3}}(\xi;\beta,I_{*})&\approx\,
\Phi_{I_{3}}(\xi)+
\begin{cases}
\lambda_{I_{*}}(\xi),&\hspace{0.5 cm}  \text{\bf(c1)};
\\[0.1 cm]
\lambda^{*}_{I_{3}^{+}(\beta)}(\xi),&\hspace{0.5 cm} \text{\bf(c2)};
\\[0.1 cm]
\lambda^{*}_{I_{3}^{-}(\beta)}(\xi),&\hspace{0.5 cm}  \text{\bf(c3)}
\end{cases}\\
&\approx\, \Phi_{I_{3}}(\xi)+
\begin{cases}
\lambda^{*}_{K^{+}(I_{*})\cap I_{3}}(\xi),&\hspace{0.5 cm}  \text{\bf(c1)};
\\[0.1 cm]
\lambda^{*}_{K^{+}(I_{*})\cap I^{+}_{3}(\beta)}(\xi),&\hspace{0.5 cm} \text{\bf(c2)};
\\[0.1 cm]
\lambda^{*}_{K^{+}(I_{*})\cap I_{3}^{-}(\beta)}(\xi),&\hspace{0.5 cm}  \text{\bf(c3)}.
\end{cases}
\end{align*}
With these remarks in mind, using (\ref{mif1}) in (\ref{r3s1}) we find that, in the region $\Gamma_{d}^{+}=J^{+}=\emptyset$
and $\sum_{n\in I_{3}}
\xi_{n}^{\overline{\ell}_{n}}\approx \zeta_{d}$,
we have $\mathcal{R}_{-2<\beta_{0}(I_{1})<0}(\xi)$ bounded by
\begin{multline}\label{w3b1b2}
\sum_{\substack{I\subseteq\Gamma_{d}^{-}\cap J^{-}\\
-2<\beta_{0}(I)<0}}
\sum_{J\subseteq (\Gamma_{d}^{-}\setminus I)\cap J^{-}}
\sum_{\substack{K\subseteq J\\
-1<\theta_{K}(\beta_{0}(I))\leq 0}}
(1+\sum_{k\in \Gamma_{d}^{-}\setminus (I\cup J)}|\log{\xi_{k}}|)\;
\mathcal{N}_{I,J,K}(\xi;\theta_{K}(\beta_{0}(I)))\\+
\sum_{j\in \Gamma_{d}^{-}}
\sum_{\substack{I\subseteq \Gamma_{j, d}^{-}\cap J^{-}\\
-2<\beta_{0}(I)<-1}}
\sum_{J\subseteq (\Gamma_{j, d}^{-}\setminus I)\cap J^{-}}
\sum_{\substack{K\subseteq J\\
-1<\theta_{K}(\beta_{0}(I))\leq -1/2}}
Q_{J, K\cup \{j\}}(\xi_{j})\\
\mathcal{N}_{I,J,K\cup \{j\},K}(\xi;\theta_{K}(\beta_{0}(I))+1/2).
\end{multline}
Using (\ref{mif2}) in (\ref{r3s1}), we deduce that $\mathcal{R}_{-2<\beta_{0}(I_{1})<-1}(\xi)$
is bounded by the second quadruple sum in (\ref{w3b1b2})
in the region $\Gamma_{d}^{+}=\emptyset$, $|J^{+}|=1$ and $\sum_{n\in I_{3}}
\xi_{n}^{\overline{\ell}_{n}}\approx \zeta_{d}$. Similarly, it follows from (\ref{mif3}) and (\ref{r3s1}) that, wherever $|\Gamma_{d}^{+}|=1$, $J^{+}=\emptyset$ and $\sum_{n\in I_{3}}
\xi_{n}^{\overline{\ell}_{n}}\approx \zeta_{d}$, we have
\begin{equation}\label{w3b3}
\mathcal{R}_{-1<\beta_{1}(I_{1})<0}(\xi)\lesssim
\sum_{\substack{I\subseteq\Gamma_{d}^{-} \cap J^{-}\\
-1<\beta_{1}(I)<0}}
\sum_{J\subseteq (\Gamma_{d}^{-}\setminus I)\cap J^{-}}
\sum_{\substack{K\subseteq J\\
-1/2<\theta_{K}(\beta_{1}(I))\leq 0}}\mathcal{N}_{I,J,K}(\xi;\theta_{K}(\beta_{1}(I))).
\end{equation}
\subsubsection{Concluding the estimate for \texorpdfstring{{\mathversion{bold}$
W_{3}(\xi)$}}{} on \texorpdfstring{{\mathversion{bold}$
U\times V$}}{}}
As discussed in Section \ref{prmlst},
$W_{3}$ has the bound (\ref{hr1}), when
$|\Gamma^{+}_{d}|\geq 2$. We put together the estimates obtained
in Sections \ref{471lgs} and  \ref{lstp} for
(portions of) the sum (\ref{r3s1}) when $|\Gamma^{+}_{d}|\leq 1$.\smallskip\par
Let $\Gamma_{d}^{+}=\emptyset$. Combining (\ref{lge2}), (\ref{lq1jnp1}), (\ref{lq1jnp2}), (\ref{lge3}) and the first two inequalities in (\ref{xinsmzt}), we discover that if we extend the first quadruple sum in (\ref{w3b1b2})
to $\beta_{0}(I)=0$, extend the second quadruple sum in (\ref{w3b1b2}) to $\beta_{0}(I)=-1$, and add the sum $\mathcal{S}_{-2<\beta_{0}(I)<0}(\xi)$, then we end up with a pointwise upper bound
for $\mathcal{R}_{-2<\beta_{0}(I_{1})\leq 0}(\xi)$. Observe that
when $\beta_{0}(I)=0$, the restriction $-1<\theta_{K}(\beta_{0}(I))\leq 0$ is possible only if
$K=\emptyset$, in which case $\mathcal{N}_{I,J,K}(\xi;\theta_{K}(\beta_{0}(I)))$ coincides with the summand in (\ref{lge2}). Similarly, if $\beta_{0}(I)=-1$ and
$-1<\theta_{K}(\beta_{0}(I))\leq -1/2$, then $K=\emptyset$ and
the summand of the second quadruple sum in (\ref{w3b1b2})
coincides with the summand in (\ref{lge3}). Finally, by
(\ref{lge121}), adding $\widetilde{\mathcal{S}}(\xi)$ to the above-described bound for $\mathcal{R}_{-2<\beta_{0}(I_{1})\leq 0}(\xi)$ yields a bound for $W_{3}(\xi)$ when $\Gamma_{d}^{+}=\emptyset$.
\smallskip\par
Assume now that $|\Gamma_{d}^{+}|=1$. According to
(\ref{lq1jnp1}) and the third inequality in (\ref{xinsmzt}),
if we add the sum
$\mathcal{S}_{-1<\beta_{1}(I)< 0}(\xi)$ to
the sum (\ref{w3b3}), we obtain a bound
for $\mathcal{R}_{-1<\beta_{1}(I_{1})< 0}(\xi)$. If, in addition, we
extend the sum (\ref{w3b3}) to $\beta_{1}(I)=0$, then, by (\ref{lge2}), we recover a bound for
$\mathcal{R}_{-1<\beta_{1}(I_{1})\leq 0}(\xi)$. Adding $\widetilde{\mathcal{S}}(\xi)$ to this bound implies an estimate for $W_{3}(\xi)$ when $|\Gamma_{d}^{+}|=1$. \smallskip\par
With these remarks in mind, we see that
\begin{equation}\label{w3f}
W_{3}(\xi)\lesssim
\mathcal{S}_{d,\Gamma_{d}^{-}}(\xi)+
\mathcal{R}_{d,\Gamma_{d}^{-}}(\xi)+\widetilde{\mathcal{S}}(\xi).
\end{equation}
\subsection{Concluding the proof}
Adding the estimates (\ref{hr1}), (\ref{w2f}) and (\ref{w3f}),
 gives an estimate for $W(\xi)$ on $U\times (d^{-\ell_{d}}/2,2)$
 whenever $\Gamma_{d}^{-}\neq \emptyset$.
Observe that
\begin{equation*}
P_{\Gamma^{-}_{d}\cap E_{2}}(\xi)\gtrsim\widetilde{\mathcal{S}}(\xi)
,\quad P_{\Gamma^{-}_{d}\cap E_{2}}(\xi)\gtrsim\mathcal{S}_{\beta_{|\Gamma_{d}^{+}|}(I)\geq 0}(\xi).
\end{equation*}
Additionally, the sum that defines
$\mathcal{R}_{d,\Gamma_{d}^{-}}$ (with $\mathcal{L}_{\Gamma^{-}_{d}\setminus (I\cup J)}$ replaced by 1) restricted to $K=J=\emptyset$ coincides with the
 sum (\ref{w2f}) restricted to $\beta_{|\Gamma_{d}^{+}|}(I)\leq 0$.
Therefore, by (\ref{smsm1}), we have
\begin{equation*}
W(\xi)\lesssim (P_{\Gamma^{-}_{d}\cap E_{2}}(\xi)+
\mathcal{S}_{d,\Gamma_{d}^{-}}(\xi)+
\mathcal{R}_{d,\Gamma_{d}^{-}}(\xi))Q_{\Gamma_{d}^{+}}(\xi),\quad \xi\in U\times (d^{-\ell_{d}}/2,2).
\end{equation*}
This agrees with (\ref{sh0}) when $\Gamma_{d}^{-}=\emptyset$.
Adding the estimates (\ref{hyg1}) and (\ref{hyg2}), we obtain
\begin{equation*}
G_{1}(\xi)+G_{2}(\xi)\lesssim\,
P_{\Gamma^{-}_{d}\cap E_{2}}(\xi)(
\chi_{(0,2]}(\xi_{d})+\xi_{d}^{-\alpha_{d}}
\chi_{[1,\infty)}(\xi_{d}))Q_{\Gamma_{d}^{+}}(\xi),\quad \xi_{d}^{\prime}\in U.
\end{equation*}
In view of (\ref{hyg0}), adding the latter two estimates gives a bound for $L^{d}(\xi)$ on $U\times \mathbb{R}_{+}$. An analogous bound can be inferred for $L^{k}(\xi)$ for each $k\in \Gamma$.
Using each of these estimates in (\ref{kl2}), and then summing in $k$ implies an estimate for $L_{\epsilon}(\xi)$ by (\ref{lps}).
This completes the proof of (\ref{mnk}).
%-----------------------------------------------------
%-----------------------------------------------------

\bibliographystyle{plain}
\bibliography{refpstm}

\end{document}